\documentclass[11pt]{amsart}
\usepackage{amssymb}
\usepackage{latexsym}
\usepackage{graphicx}
\usepackage{subfigure}
\usepackage{tikz}
\usepackage[linktocpage=true]{hyperref}
\usepackage{ulem}
\usepackage[left=3cm, right=3cm, top=3cm, bottom=3cm]{geometry}
\usepackage{color}

\hypersetup{ colorlinks   = true, 
urlcolor  = blue, 
linkcolor    = blue, 
citecolor   = blue 
}

\def\dist{{\rm dist}}
\def\sl{{\rm sl}(2,\mathbb{R})}
\def\SL{{\rm SL}(2,\mathbb{R})}

\def\om{{\omega}}

\def\lm{{\lambda}}

\def\vp{{\varphi}}

\def\beq{\begin{equation}}
\def\eeq{\end{equation}}

\def\bm{\begin{matrix}}
\def\em{\end{matrix}}

\newcommand{\Z}{{\mathbb Z}}
\newcommand{\R}{{\mathbb R}}

\newcommand{\T}{{\mathbb T}}

\newcommand{\N}{{\mathbb N}}

\newcommand{\CF}{{\mathcal F}}
\newcommand{\CG}{{\mathcal G}}
\newcommand{\CB}{{\mathcal B}}
\newcommand{\CH}{{\mathcal H}}
\newcommand{\CI}{{\mathcal I}}
\newcommand{\CL}{{\mathcal L}}
\newcommand{\CM}{{\mathcal M}}

\newcommand{\CP}{{\mathcal P}}

\newcommand{\CZ}{{\mathcal Z}}
\newcommand{\CU}{{\mathcal U}}

\newcommand{\CK}{{\mathcal K}}

\newcommand{\CJ}{{\mathcal J}}
\newcommand{\CR}{{\mathcal R}}

\newtheorem{Theorem}{Theorem}[section]
\newtheorem{remark}{Remark}[section]
\newtheorem{lemma}{Lemma}[section]

\newtheorem{Definition}{Definition}[section]
\newtheorem{prop}{Proposition}[section]

\newcommand{\la}{\langle}
\newcommand{\ra}{\rangle}

\begin{document}

\title[]{Almost reducibility and oscillatory growth of Sobolev norms}

\author{Zhenguo Liang}
\address{School of Mathematical Sciences and Key Lab of Mathematics for Nonlinear Science, Fudan University, Shanghai 200433, China}
\email{zgliang@fudan.edu.cn}
\thanks{Z. Liang was partially supported by
NSFC grant (12071083) and Natural Science Foundation of Shanghai (Grants No. 19ZR1402400).}

\author{Zhiyan Zhao}
\address{Universit\'e C\^ote d'Azur, CNRS, Laboratoire J. A. Dieudonn\'{e}, 06108 Nice, France}
\email{zhiyan.zhao@univ-cotedazur.fr}
\thanks{Z. Zhao was partially supported by the French government through the National Research Angency (ANR) grant for the project KEN ANR-22-CE40-0016, and by
NSFC grant (11971233).}

\author{Qi Zhou}
\address{Chern Institute of Mathematics and LPMC, Nankai University, Tianjin 300071, China
} \email{qizhou@nankai.edu.cn}
\thanks{Q. Zhou was partially supported by National Key R\&D Program of China (2020 YFA0713300), NSFC grant (12071232), the Science Fund for Distinguished Young Scholars of Tianjin (No. 19JCJQJC61300) and Nankai Zhide Foundation}

\

\begin{abstract}
For 1D quantum harmonic oscillator perturbed by a time quasi-periodic quadratic form of $(x,-{\rm i}\partial_x)$, we show its almost reducibility.
The growth of Sobolev norms of solution is described based on the scheme of almost reducibility. In particular, an $o(t^s)-$upper bound is shown for the $\CH^s-$norm if the equation is non-reducible.
Moreover, by Anosov-Katok construction, we also show the optimality of this upper bound,
i.e., the existence of quasi-periodic quadratic perturbation for which the growth of ${\mathcal H}^s-$norm of the solution is $o(t^s)$ as $t\to\infty$ but arbitrarily ``close" to $t^s$ in an oscillatory way.

\

\noindent
{\bf Keywords:}
1D quantum harmonic oscillator; time quasi-periodic; quadratic quantum Hamiltonian; Metaplectic representation; almost reducibility; growth of Sobolev norms

\

\noindent
{\bf MSC 2020:} 35Q40, \ 37C60, \ 35Q41, \ 47G30

\end{abstract}

\maketitle
\underline{}
\section{Introduction and main results}

The quadratic quantum Hamiltonian plays an important role in partial differential equations, since it gives non-trivial examples of wave propagation phenomena in quantum mechanics.
A typical example is the quantum Hamiltonian given by the Schr\"odinger operator $-\hbar^{2} \triangle+V$, with $V$ a harmonic potential, which has been well studied since the last century. In this paper, we consider the time quasi-periodic Schr\"odinger equation
\begin{equation}\label{eq_original}
{\rm i}\partial_t\psi=(\CJ+{\CP}(t))\psi,
	\end{equation}
where, we assume that
\begin{itemize}
  \item the operator $\CJ :=D^2+X^2$ is the one-dimensional quantum harmonic oscillator, with
  $$(Du)(x):=-{\rm i}u'(x), \quad  (Xu)(x):= x u(x), \qquad \forall \ u\in L^2(\R),$$
  \item the perturbation $\CP(t)$ is a time quasi-periodic self-adjoint quadratic form of $(X,D)$:
  \begin{equation}\label{SA_qua_form}
{\CP}(t)= \frac12\left(p_{20}(\om t)X^2+p_{11}(\om t)XD+ p_{11}(\om t) DX + p_{02}(\om t)D^2\right),
  \end{equation}
(or more precisely, the symbol of ${\CP}(t)$ is a quadratic form)
 with the frequency vector $\omega\in \R^d$, $d\geq 1$, satisfying the {\it Diophantine} condition (denoted by $\omega\in {\rm DC}_d(\gamma,\tau)$ for $\gamma>0$, $\tau>d-1$):
      \begin{equation}\label{dio}
      \inf_{j\in\Z}|\la n,\omega\ra-j|>\frac{\gamma}{|n|^\tau},\qquad \forall \ n\in\Z^d\setminus\{0\},
      \end{equation}
and with the coefficients $p_{20},\, p_{11}, \, p_{02}$ real-analytic functions on $\T^d=(\R/2\pi\Z)^d$, extending to a bounded complex analytic function on $\{|{\rm Im}z|<r\}$ for some $r>0$.
\end{itemize}
The aim of this paper is to exploit the behavior of solutions to Eq. (\ref{eq_original}) in Sobolev space
\begin{equation}\label{SobolevSpace}
{\CH}^s:=
\left\{\psi\in L^2(\R):\CJ^{\frac{s}2}\psi \in L^2(\R)\right\}, \quad s>0,
\end{equation}
through almost reducibility of Eq. (\ref{eq_original}). In particular, for the non-trivial solution to Eq. (\ref{eq_original}), an oscillatory behavior of Sobolev norm (${\CH}^s-$norm) $\|u(t)\|_{{\CH}^s}:=\|\CJ^{\frac{s}2} u(t)\|_{L^2}$ is observed, with an optimal $o(t^s)$ growth upper bound.

\subsection{Almost reducibility}

For quadratic Hamiltonians, it is well known that the correspondance between classical and quantum mechanics is exact. If the coefficients of the quadratic form are constant, there are explicit formulas for solutions, such as Melher formula (see e.g. \cite{Fol1989}) for the quantum harmonic oscillator.
But for time dependent quadratic Hamiltonians, explicit formulas for solutions are usually not expectable. In this sense, it is natural to wonder whether the time dependent quantum Hamiltonian is {\it reducible}, i.e., conjugated to an equation independent of time via a time dependent unitary transformation. This unitary transformation is sometimes the limit of a convergent sequence (in a suitable sense) of time dependent transformations constructed by a KAM scheme.

Many previous works are devoted to reducibility for time dependent Hamiltonian PDEs. Regarding 1D quantum harmonic oscillator, the time periodic smooth perturbations were firstly considered \cite{Com87, DLSV2002, EV83, HLS86, Kuk1993}.
In the time quasi-periodic case, bounded perturbations \cite{GreTho11, LiangWZQ2022, Wang08, WLiang17} and
unbounded ones \cite{Bam2018, Bam2017, BM2018, Liangluo19, LZZ2021, LLZ2022} are both investigated.
It is worth  mentioning that, in studying reducibility problems, KAM theory for 1D PDEs with unbounded perturbations has been well developed by Bambusi-Graffi \cite{BG2001}, Kuksin \cite{Kuk1997} and Liu-Yuan \cite{LY2010}.
For PDEs in higher dimension, reducibility was initiated by Eliasson-Kuksin \cite{EK2009, EK2010} for the quasi-periodic Schr\"odinger equation.
As for higher-dimensional quantum harmonic oscillator, we can refer to \cite{GrePat16, LiangWang19, LiangWZQ2022plus} for bounded potentials, and to \cite{BGMR2018} for unbounded perturbations which include the quadratic ones.
Through reducibility, all above results give the indications to the long-time behaviors of solutions (at least upper bounds of Sobolev norms), or the spectral feature of Floquet operator.

However, reducibility is not always feasible since the problem of small divisors are unavoidable in the KAM construction.
As an alternative concept, {\it almost reducibility} was introduced by Eliasson \cite{Eli1992} for quasi-periodic linear systems, and appears subsequently in many other settings of dynamical systems \cite{A2015, AJ2010, AK2022, E2017}. The main principle of almost reducibility is to conjugate the time dependent system to constant with a possibly divergent sequence of transformations.



\medskip

For the time dependent Schr\"odinger equation (\ref{eq_original}), noting that $\CP$ is uniquely determined by the coefficients $p_{20}(\cdot)$, $p_{11}(\cdot)$, $p_{02}(\cdot)$, one can endow the topology of $\CP$ by the topology of the function triplet $(p_{20}, p_{11}, p_{02})$. For example, one can define
$$\|\CP\|_r:= |p_{20}|_r+|p_{11}|_r+|p_{02}|_r=\sup_{|{\rm Im}z|<r} |p_{20}(z)|+\sup_{|{\rm Im}z|<r} |p_{11}(z)|+\sup_{|{\rm Im}z|<r} |p_{02}(z)|.$$
The main result about almost reducibility is the following:
\begin{Theorem}\label{thm_al_re} There exists $\varepsilon_*=\varepsilon_*(\gamma, \tau, d, r)>0$ such that if $\|\CP\|_r <\varepsilon_* $,
then Eq. (\ref{eq_original}) is almost reducible, i.e., for every $j\in\N^*$, there exists an equation
\begin{equation}\label{eq_ar_j-thm}
{\rm i}\partial_t\psi_j=({\CL}_j+{\CP}_{j+1}(t))\psi_j,
  \end{equation}
where both ${\CL}_j$ and ${\CP}_{j+1}(t)$ are self-adjoint quadratic forms of $(X,D)$ as (\ref{SA_qua_form}), with the coefficients of ${\CL}_j$ constant, and with those of ${\CP}_{j+1}(t)$ time quasi-periodic and tending to $0$ as $j\to\infty$, such that Eq. (\ref{eq_original}) is conjugated to Eq. (\ref{eq_ar_j-thm}) via a time quasi-periodic $L^2-$unitary transformation $\psi(t)={\CU}_j(t)  \psi_j(t)$.
\end{Theorem}

The reducibility is already studied by the authors for Eq. (\ref{eq_original}), with the coefficients smoothly depending on a well interpolated parameter $E\in {\CI}\subset\R$ (see Theorem 1 in \cite{LZZ2021}). It is shown that, for almost every $E$, Eq. (\ref{eq_original}) is reducible, hence almost reducible. Through Theorem \ref{thm_al_re}, we see ``merely" almost reducibility for the zero-measure subset of parameters not mentioned in \cite{LZZ2021}.
However, it will be shown that, there exist plenty of time quasi-periodic perturbations $\CP(t)$ such that the corresponding sequence of unitary transformations $\{{\CU}_j(t)\}$ obtained in Theorem \ref{thm_al_re} ``diverges" in a suitable sense. Such a phenomenon does not occur if the coefficients are time periodic \cite{HLS86} since the reducibility is guaranteed by Floquet theory.

\subsection{Growth of Sobolev norm}

As an important problem in Mathematical Physics, the behavior of solutions to Hamiltonian PDEs in Sobolev space, proposed by Bourgain \cite{Bou96} and followed by other pioneering works (e.g., \cite{CKSTT2010, GG2016, GK2015, HPTV2015}), has attracted plenty of interest during the past decades. Besides the boundedness of Sobolev norms and the upper bound of growth,
it is more challenging to investigate the existence of unbounded trajectories in Sobolev space, naturally related to weak turbulent effects and energy cascades, and precise the growth rate with time of Sobolev norms. During the last years, fruitful achievements have been made on this theme for various Hamiltonian PDEs.

As mentioned in the previous subsection, it is effective to study the long-time behavior of solutions through reducibility of the original quantum Hamiltonian, especially for the linear Hamiltonians. Following \cite{GY00}, Bambusi-Gr\'ebert-Maspero-Robert \cite{BGMR2018} considered 1D quantum harmonic oscillator with a time periodic linear potential, and got $t^{s}-$polynomial growth for ${\CH}^s-$norm, by reducibility to a transport equation. This strategy indeed also works for higher-dimensional quantum harmonic oscillator with perturbations linear in $(X,D)$ (via Theorem 3.3 of \cite{BGMR2018}). For 1D quantum harmonic oscillator with a time quasi-periodic perturbation which is a polynomial of $(X,D)$ of degree $2$, various growth rates of ${\CH}^s-$norm have been shown according to different normal forms of reducibility \cite{LZZ2021, LLZ2022}:
$t^{s}-$polynomial growth for parabolic normal form, exponential growth for hyperbolic normal form, $t^{2s}-$polynomial growth for reducibility to Stark Hamiltonian.

For Hamiltonian PDEs, almost reducibility is also useful to see the long-time behaviors of solutions.
Eliasson \cite{E2017} showed almost reducibility for the quasi-periodic linear wave equation, through which a log-log-bound on Sobolev norms of solutions is obtained.
Similar idea was employed by Bambusi and his collaborators \cite{BGMR2021, BLM20, BL2022} for time dependent Schr\"odinger equations, with a more and more regularized remaining part along with the iteration step, instead of the usual asymptotic smallness assumption. A $t^\epsilon-$upper bound, for arbitrary $\epsilon>0$, on Sobolev norms of solutions is obtained through such an argument.

Another way usually used to get the unbounded trajectories, is to construct specific perturbations towards the growth to infinity.
For 1D quantum harmonic oscillator,
Delort \cite{Del2014} (followed by a refined proof of Maspero \cite{Mas2018}) constructed a time periodic order-zero pseudo
differential operators as the perturbation such that some solution exhibits $t^{\frac{s}{2}}-$polynomial growth of ${\CH}^s-$norms, and in an abstract setting, Maspero \cite{Mas2021, Mas2022} exploited further time periodic perturbations and gave sufficient conditions such that some solution exhibits such polynomial growth.
For 2D quantum harmonic oscillator, by embedding Arnold diffusion into the infinite-dimensional quantum system, Faou-Rapha\"el \cite{FR20} constructed a potential decaying in time such that some solution presents logarithmic growth with time in Sobolev space, and, based on the study in \cite{ST20} for linear Lowest Landau Level equations with a time dependent potential, Thomann \cite{Thomann20} constructed the perturbation being the projection onto Bargmann-Fock space such that some travelling wave whose Sobolev norm presents polynomial growth with time.
In other settings, the logarithmic growth of Sobolev norms was also shown by Bourgain \cite{Bou99} for 1D and 2D linear Schr\"odinger equations with quasi-periodic potential, and by Haus-Maspero \cite{HM2020} for semiclassical anharmonic oscillators with regular time dependent potentials.



\medskip

Let us come back to Eq.(\ref{eq_original}). We consider the behavior of its solution in Sobolev space defined as in (\ref{SobolevSpace}). In view of Theorem 1.2 of \cite{MR2017}, the solution to the quadratic quantum Hamiltonian (\ref{eq_original}) is globally well-posed in ${\CH}^s$ space with a general exponential upper bound. According to Theorem 1 of \cite{LZZ2021}, Sobolev norms always present normal behaviors (boundedness, polynomial growth, exponential growth) if the quadratic Hamiltonian (\ref{eq_original}) is reducible.
It is interesting to consider the behavior of Sobolev norms when the quadratic quantum Hamiltonian is almost reducible instead of reducible, especially in the case that Eq. (\ref{eq_original}) is {\bf non-reducible}, i.e., the sequence of unitary transformations $\{{\CU}_j(t)\}$ obtained in Theorem \ref{thm_al_re} does not converge to a $L^2-$unitary transformation uniformly in $t$.

{The behavior of Sobolev norms of solutions to the non-reducible equation is described in the following theorems.
\begin{Theorem}\label{thm_growth_upper}
If $\|\CP\|_{r}<\varepsilon_*$ as in Theorem \ref{thm_al_re}, and the quadratic Hamiltonian (\ref{eq_original}) is non-reducible,
then for the solution $\psi(t)$ to Eq. (\ref{eq_original}) with $\psi(0)\in \CH^{s+2}$, we have
\begin{equation} \label{upper-thm}
\lim_{t\to +\infty}\frac{\|\psi(t)\|_{\CH^s}}{t^s}= 0.
\end{equation}
\end{Theorem}


\medskip

With the $o(t^s)-$upper bound of ${\CH}^s-$norm stated in Theorem \ref{thm_growth_upper}, one may naturally wonder if it is optimal. Based on Anosov-Katok construction from dynamical systems, we will construct
explicit ``\textit{dense}" perturbations ${\CP}(t)$, such that the corresponding quantum Hamiltonian exhibits an optimal $o(t^s)-$upper bound in $\CH^s-$space:

\begin{Theorem}\label{thm_sub_polynome}
For $s>0$ and non-vanishing $\psi(0)\in {\CH}^{s+2}$, given $f:\R_+^*\to \R_+^*$ with
$$f(t)\to\infty,\quad f(t)=o(t^s), \qquad  t\to\infty,$$
there exists time quasi-periodic quadratic perturbation ${\CP}(t)$ such that for the solution $\psi$ to Eq. (\ref{eq_original}), we have
\begin{eqnarray}
\limsup_{t\to+\infty}\frac{\|\psi(t)\|_{{\CH}^s}}{f(t)}&=& \infty, \label{limsup-thm}\\
\liminf_{t\to+\infty}\|\psi(t)\|_{{\CH}^s}&<& \infty.\label{liminf-thm}
\end{eqnarray}
Moreover, for any $\varepsilon>0$, there exists ${\CP}(t)$ with  $\|\CP\|_r<\varepsilon $ such that (\ref{limsup-thm}) and (\ref{liminf-thm}) are satisfied.
\end{Theorem}

One may interpret the result as follows. For any $f(t)=o(t^s)$ which may grow with a rate arbitrarily ``close" to $t^s$, e.g.,
$$f(t)=\frac{t^s}{\log\log\log(e^e+t)},$$
according to (\ref{limsup-thm}) and (\ref{liminf-thm}), there exist two sequences of moments $\{T_j\}$, $\{t_j\}$, both of which tending to $\infty$, and a constant $c$, depending on $s$ and $\psi(0)$, such that
\begin{equation}\label{rem_seq_moments}
\|\psi(T_j)\|_{{\CH}^s}\geq f(T_j)= \frac{T_j^s}{\log\log\log(e^e+T_j)},\quad \|\psi(t_j)\|_{{\CH}^s}\leq c,\qquad \forall \ j\in\N.
\end{equation}
Note that Theorem \ref{thm_sub_polynome} gives an oscillatory growth of Sobolev norms for the solution to Eq. (\ref{eq_original}). This kind of result was also obtained by G\'erard-Grellier \cite{GG2016} for the cubic Szeg\H{o} equation.

\subsection{Ideas of proof}

 It is known that the classical-quantum correspondence is realized by Weyl quantization:
 given any symbol $f=f(x,\xi)$, with $x,\xi\in\R^n$, $n\geq 1$, the Weyl operator $f^{W}$ of $f$ is defined as
 \begin{equation}\label{WeylQuantization}
 \left(f^{W} u\right)(x)=\frac{1}{(2\pi)^n}\int_{y, \, \xi\in\R^n} e^{{\rm i}(x-y)\xi} \,  f\left(\frac{x+y}{2},\xi\right) u(y) \, dy  \, d\xi,\qquad \forall \ u\in L^2(\R^n).
 \end{equation}
In particular, if $f$ is a polynomial of degree at most $2$ in $(x,\xi)$, then $f^W$ is a polynomial of degree at most $2$ in $(X,D)$ after the symmetrization. In this sense, in order to well understand Eq. (\ref{eq_original}), it is crucial to study the corresponding finite-dimensional Hamiltonian system. This idea has been employed in \cite{BGMR2018} to obtain the reducibility and the almost preservation of Sobolev norms for the time quasi-periodic perturbed quantum harmonic oscillator, by establishing a classical KAM scheme for ``most" frequency vector. With the same idea, in 1D case, the reducibility and various behaviors of Sobolev norms have been described in \cite{LZZ2021, LLZ2022}, by adapting the KAM scheme of reducibility for time quasi-periodic $\sl-$linear systems developed by Eliasson \cite{Eli1992}, in which the reducibility transformation is not necessarily close to identity, to the quantum Hamiltonian.

Different from the above achievements, the present work deals with the quantum Hamiltonian for which the reducibility is not expectable.
In this situation, it is not possible to simply conjugate the original time dependent Hamiltonian PDE to a constant one by a unitary transformation, but the theory of Eliasson \cite{Eli1992} is still helpful to realize the almost reducibility of Eq. (\ref{eq_original}), by asymptotically conjugating the original Hamiltonian to constant via a sequence of unitary transformations.
The possible divergence of this sequence  makes it useless to exploit the long-time behavior of solutions by simply calculating with the constant quantum Hamiltonian.
It is more reasonable to focus on a certain time dependent state after the application of finitely many unitary transformations, with the time dependent part arbitrarily small according to the prescribed time.

To show the optimality of  $o(t^s)-$upper bound in Theorem \ref{thm_growth_upper}, the proof of Theorem \ref{thm_sub_polynome} is based on fast approximation conjugation method or the well-known Anosov-Katok construction. This is the first novelty of this paper. Anosov-Katok \cite{AK1970} first developed this method to construct mixing diffeomorphisms
of the unit disc arbitrarily close to Liouvillean rotations.  We give a sketch of the constructions as follows. The irrational number $\rho$ is Liouvillean, if there exist sequences $p_{n}, q _{n} \in \N ^{*}$ such that
\begin{equation} \label{liou}
\left| \rho - \frac{p_{n}}{q_{n}}\right| < |q_{n}|^{-n}.
\end{equation}
It is
tempting to study diffeomorphisms with rotation number equal to $\rho$, by
approximating them from diffeomorphisms that are conjugate to    periodic rotation $T_{p_{n}/q_{n}}.$
A
diffeomorphism $f_{n}$ is constructed such that
\[
f_{n} = H_{n}\circ T_{p_{n}/q_{n}} \circ H^{-1}_{n}
\]
for some smooth conjugation $H_{n}$, satisfying some finitary version of mixing at a scale
of iteration $q_{n}$. A $q_{n}$-periodic conjugation $h_{n}$ is constructed,
such that the diffeomorphism
\begin{equation}\label{hn}
f_{n+1} = H_{n} \circ h_{n} \circ T_{p_{n+1}/q_{n+1}} \circ h_{n} ^{-1}  \circ H_{n} ^{-1}
\end{equation}
satisfies some improved finitary version of mixing at a scale of iteration $q_{n+1}$.
The Liouvillean character of $\rho$ is necessary in order to assure both the convergence
of $f_{n}$ and the divergence of $H_{n}$.

If one wants to apply Anosov-Katok construction to quadratic quantum Hamiltonian, the basic question is what the corresponding periodic rotation $T_{p_{n}/q_{n}}$ is. Actually, in our case, we will develop a fibred Anosov-Katok construction, and the periodic quantum Hamiltonian takes the form
$${\CL}_n=(\varphi_n+\lambda_n)D^2+(\varphi_n-\lambda_n)X^2 \   \  {\rm with} \  \ \sqrt{\varphi_n^2-\lambda_n^2}=\la k_n,\om\ra, $$
i.e., rational with respect to the base frequency $\om$.
At a scale of iteration $k_n$, we will show that the approximated quantum Hamiltonian ${\CL}_n$ presents large but oscillatory growth of Sobolev norms, if $\la k_n,\om\ra$ shrinks fast enough, and the limit Liouvillean quantum Hamiltonian will present optimal $o(t^s)$ bound. It is worth pointing out that classical Anosov-Katok construction is only valid in $C^{\infty}$ topology \cite{FK2005}, while we are able to carry out the  construction in analytic topology. One may consult Section \ref{sec_al_ODE} for detailed constructions and more comments.

Similar fibered Anosov-Katok construction has also been applied to quasi-periodic $\SL-$cocycles \cite{KXZ2021}.
One can consult a nice survey of Fayad-Katok \cite{FK2005} on  results and references of classical Anosov-Katok construction. In principle,
Anosov-Katok construction has been useful in producing examples of dynamics incompatible with quasi-periodicity, in the vicinity of quasi-periodic dynamics.
It is in some sense the counterpart of the KAM method.
In our case, we develop almost reducibility (KAM method) to show the quantum Hamiltonian \eqref{eq_original} has $o(t^s)-$upper bound, while we develop Anosov-Katok construction method to prove this kind of growth is optimal.

To carry out the above strategies on Sobolev norms of solutions with almost reducibility via the classical and quantum Hamiltonians, more quantitative computations are required.
Noting the homogeneity of Eq. (\ref{eq_original}), we make use of the Metaplectic representation \cite{Fol1989} instead of Weyl quantization (\ref{WeylQuantization}) which works in more general settings. This is another novelty of this paper.
The Metaplectic representation is a ``fascinating double-valued unitary representation of the symplectic group", and in particular, it is explicitly defined by two integral operators for $\SL$ group (corresponding to 1D quadratic quantum Hamiltonian).
Moreover, the conclusions in \cite{LZZ2021} can also be shown through the Metaplectic representation, which simplifies the proof.

\subsection{Notations}\label{sec_Nota}

For the convenience of readers, we introduce several notations frequently used through this paper.

\begin{enumerate}
  \item The inequality with `` $\lesssim_s$" (or `` $\gtrsim_s$") means boundedness from above (or below) by a positive constant depending only on $s\in[0,\infty[$ but independent of other factors, i.e., for $a,b\geq 0$, the inequality $a \lesssim_s b$ (or $a \gtrsim_s b$) means that there is a constant $c_s$ depending only on $s$ such that $a \leq c_s b$ (or $a \geq c_s b$). Moreover, $a\sim_s b$ means that $a \lesssim_s b$ and $a \gtrsim_s b$.

\smallskip

  \item For $0\leq s'\leq s$, let ${\CB}({\CH}^{s}\to{\CH}^{s'})$ be the space of bounded operators from ${\CH}^s$ to ${\CH}^{s'}$ with $\|\cdot\|_{{\CB}({\CH}^s\to{\CH}^{s'})}$ the operator norm on ${\CH}^s$, simplified to ${\CB}({\CH}^s)$ and  $\|\cdot\|_{{\CB}({\CH}^s)}$ if $s'=s$. Given non-vanishing $u\in{\CH}^s$, let us define
  $$\Gamma_s(u):=\frac{\|u\|_{{\CH}^s}}{\|D^s u\|_{L^2}}.$$

\smallskip

  \item Let ${\rm SL}(2,\R)$ (or $\sl$) be the special linear group (or Lie algebra) of $2\times2$ real matrices with the matrix norm defined by
  $$\|A\|=\sum_{i,j=1,2}|a_{ij}|,\quad A=(a_{ij})_{i,j=1,2}\in\SL \  ({\rm or} \   \sl).$$
  Moreover, we define the identity, symplectic and rotation matrices $I$, $J$, $R_\theta$, by
  $$I:=\left(\begin{array}{cc}
               1 & 0 \\
               0 & 1
             \end{array}\right),\quad J:=\left(\begin{array}{cc}
               0 & 1 \\
               -1 & 0
             \end{array}\right),\quad R_\theta:=\left(\begin{array}{cc}
               \cos(\theta) & \sin(\theta) \\
               -\sin(\theta) & \cos(\theta)
             \end{array}\right), \;\  \theta\in\R.$$

\smallskip

\item For given $A=\left(\begin{array}{cc}
               a_{11} & a_{02} \\
               -a_{20} & -a_{11}
             \end{array}\right)\in\sl$, let
$$\frac{-1}2 \left\la \left(\begin{array}{c}
X\\
D
 \end{array}\right) , J A \left(\begin{array}{c}
X\\
D
 \end{array}\right) \right\ra$$
be the self-adjoint quadratic pseudo-differential operator
$$\frac{1}2\left(a_{20} X^2 + a_{11}(XD+DX)+a_{02}D^2\right).$$

\smallskip

  \item With $\om\in {\rm DC}_d$ the Diophantine vector defined as above and fixed, we denote $$\la k \ra:=\la k,\om\ra,\quad k\in\Z^d.$$

\smallskip

\item For $r>0$, let $C_r^\om(\T^d,\sl)$ be the space of real-analytic $\sl-$valued functions which extends to the complex neighbourhood $\{|{\rm Im} z|<r\}$ of $\T^d$, with the norm
$$\|A\|_r:= \sup_{|{\rm Im} z|<r}\|A(z)\|,\quad  A\in C_r^\om(\T^d,\sl).$$
In particular, $C^\om(\T^d,\sl)$ denotes the space of real-analytic $\sl-$valued functions on an imprecise complex neighbourhood of $\T^d$, and we define
$$\|A\|_{\T^d}:= \sup_{\theta\in\T^d}\|A(\theta)\|,\quad  A\in C^\om(\T^d,\sl).$$
The above spaces and norms can be defined similarly for $\SL$ instead of $\sl$.
\end{enumerate}

\subsection{Description of the remaining of paper}

Section \ref{sec_Meta} contains the precise definition of the Metaplectic representation of $\SL$ group in Definition \ref{def_meta}, its basic properties formulated in Proposition \ref{prop_elem_meta} (with the proof provided in Appendix \ref{app_pr_prop}), and its application on the estimates of Sobolev norms formulated mainly in Proposition \ref{cor_metaplectic_upper} and \ref{cor_metaplectic_lower}. The estimates are based on the elementary calculations formulated in Proposition \ref{lemma_Hs} for the function of the form $U_{{\CG}, \, {\CK}}(x)=e^{{\rm i}{\CG} x^2} u({\CK} x)$, ${\CG}\in\R$, ${\CK}\in\R\setminus\{0\}$.

Section \ref{sec_AR} is devoted to the proof of almost reducibility of Eq. (\ref{eq_original}) formulated in Proposition \ref{prop_upper},
as well as the concrete construction of time quasi-periodic perturbation formulated in Proposition \ref{prop_lower} prepared for the specific oscillatory growth. These two propositions for quantum Hamiltonians will be shown through their parallel statements in classical Hamiltonians formulated in Proposition \ref{prop_al_ODE} and \ref{prop_AK_ODE}. Proposition \ref{prop_al_ODE} is indeed another version of the almost reducibility theory of Eliasson \cite{Eli1992}, and we give a sketch of proof in Appendix \ref{App_proof_ar}.

Section \ref{sec_pr_thm} is devoted to the proofs of Theorem \ref{thm_growth_upper} and Theorem \ref{thm_sub_polynome}. The growth of Sobolev norms is computed with the example constructed in Proposition \ref{prop_lower}. The $o(t^s)-$upper bound (\ref{upper-thm}) is formulated in Proposition \ref{prop_growth_upper}, and the optimality of upper bound, as well as the oscillatory growth of Sobolev norms, (\ref{limsup-thm}), (\ref{liminf-thm}), is formulated in Proposition \ref{prop_limsup} and \ref{prop_liminf}.

\medskip

\noindent
{\bf Acknowledgement.} The authors are grateful to D. Bambusi and A. Maspero for fruitful discussions during these years, especially to A. Maspero for helping formulate Lemma \ref{lemma_wellposed} in this manuscript by simplifying Theorem 1.2 of \cite{MR2017}.
Z. Zhao would like to thank J. Zheng for an inspiring discussion about Mehler's formula as well as the Metaplectic representation, and thank the Key Lab of Mathematics for Nonlinear Science of Fudan University (China) for its hospitality during his visits in 2021.

\section{Metaplectic representation and Sobolev norm}\label{sec_Meta}

Let us introduce in this section the Metaplectic representation on $\SL$ and its application on the calculation of Sobolev norms. It is well known that, for $s\in \N$, the definition of $\CH^s-$norm is equivalent to
\begin{equation}\label{norm_equiv}
\sum\limits_{\alpha+\beta\leq s\atop{\alpha,\beta\in\N} }\|x^{\alpha}\cdot\partial^{\beta} \psi\|_{L^2}.
\end{equation}
In view of Proposition 2.3 of \cite{BM2018},
we see that, for a given $\psi\in {\CH}^s$,
\begin{equation}\label{norm_equiv-1}
\|\psi\|_{\CH^s}\simeq \|\psi\|_{H^s}+ \|x^{s} \psi\|_{L^2},
\end{equation}
where `` $H^s$ " means the standard Sobolev space and $\|\cdot\|_{H^s}$ is the corresponding norm. Hence, to calculate the norm $\|\psi\|_{\CH^s}$, $s\geq 0$, it is sufficient to focus on $\|X^{s} \psi\|_{L^2}$ for $s\geq0$ and $\|D^s\psi\|_{L^2}$ for $s\in\N$.

Define the {\bf Fourier transform} ${\CF}$ and its inverse on $L^2(\R)$ as
\begin{eqnarray*}
({\CF}u)(\xi) \, = \, \widehat{u}(\xi)&=&\frac{1}{(2\pi)^{\frac{1}2}}\int_{\R} e^{-{\rm i} x\xi } \, u(x) \,  dx, \quad u\in L^2(\R), \\
({\CF}^{-1}\widehat{u})(x) \, = \, u(x)&=&\frac{1}{(2\pi)^{\frac{1}2}}\int_{\R} e^{{\rm i}x\xi} \, \widehat{u}(\xi) \,  d\xi, \quad \widehat{u}\in L^2(\R).
\end{eqnarray*}
For given $s\geq 0$, it is well known that $\CF^{\pm1}:\CH^s\to \CH^s$ and
$$\|D^s u\|_{L^2}\sim_s \|X^s \CF u\|_{L^2}, \quad \|X^s u\|_{L^2}\sim_s \|D^s \CF u\|_{L^2}, \quad \|u\|_{\CH^s}\sim_s \|\CF u\|_{\CH^s}, \qquad u\in\CH^s.$$

\subsection{Metaplectic representation on $\SL$ -- Introduction}\label{sec_def_meta}

The Metaplectic representation, also called oscillator representation, is originally defined through the Heisenberg group and Schr\"odinger representation on the symplectic group ${\rm Sp}(n,\R)$ (refer to Chapter 4 of \cite{Fol1989}). In this paper, we focus on the case $n=1$, i.e., ${\rm Sp}(n,\R)=\SL$, and give the simplified definition of the Metaplectic representation by explicit formulas with integral operators.

\begin{Definition}\label{def_meta}
For given $A=\left(\begin{array}{cc}
               a & b \\
               c & d
             \end{array}\right)\in \SL$, the {\bf Metaplectic representation}, denoted by ${\CM}(A)$, is the $L^2-$unitary transformation defined as
\begin{itemize}
  \item if $a\neq 0$, then, for $u\in L^2(\R)$,
  \begin{equation}\label{Meta_a}
  ({\CM}(A)u)(x) = \frac{1}{(2\pi a)^{\frac{1}2}}\int_{\R} e^{{\rm i} \left(\frac12 ca^{-1}x^2+a^{-1} x\xi-\frac12 a^{-1}b\xi^2\right)} \widehat{u}(\xi) d\xi;
  \end{equation}
  \item if $b\neq 0$, then, for $u\in L^2(\R)$,
  \begin{equation}\label{Meta_b}
  ({\CM}(A)u)(x) = \frac{{\rm i}^{\frac{1}{2}}}{(2\pi b)^{\frac{1}2}}\int_{\R} e^{{\rm i} \left(\frac12 db^{-1}x^2+ b^{-1} xy+\frac12b^{-1}ay^2\right)} u(y) dy.
  \end{equation}
\end{itemize}
We shall call ${\CM}(A)$ as Meta representation for short.
\end{Definition}
\begin{remark}\label{rmk_def}
The above formulas are indeed Theorem (4.51) and (4.53) of \cite{Fol1989} in 1D case.
It can be shown similarly as in \cite{Fol1989} that the formulas (\ref{Meta_a}) and (\ref{Meta_b}) agree if both $a$ and $b$ are non-vanishing.
We note that the formulas (\ref{Meta_a}) and (\ref{Meta_b}) are not exactly the same with \cite{Fol1989} because of the slightly different definition of Fourier transform and its inverse, as well as the adaption to the solution to Eq. (\ref{eq_original}).

\end{remark}

\begin{remark}\label{rmk_signs}
It has been shown in \cite{Fol1989} that ${\CM}(A)$ is determined up to factors $\pm1$ and ${\CM}$ is a double-valued unitary representation of $\SL$. With (\ref{Meta_a}) and (\ref{Meta_b}), we see that such ambiguity of sign appears in the square root. Through the paper, we will neglect the ambiguity of sign since the Meta representation will be used on the computations of Sobolev norms.
\end{remark}

Briefly speaking, the Meta representation given in Definition \ref{def_meta} transforms any $\SL-$matrix to a $L^2-$unitary transformation on $L^2(\R)$, by which more information of the quadratic quantum Hamiltonian can be read off from the corresponding finite-dimensional linear systems.

\begin{prop}\label{prop_elem_meta}(Properties of the Meta representation.)
 The following assertions hold.
\begin{enumerate}
\item [(i)] (Mehler's formulas) For $\theta\in\R$ with $\cos(\theta)\neq0$,
$${\CM}(R_\theta) u(x)= \frac{1}{\sqrt{2\pi \cos(\theta)}}\int_{\R} e^{-\frac{\rm i}{\cos(\theta)}\left(\frac{x^2+\xi^2}{2}\sin(\theta)-x\xi\right)}\widehat u(\xi) \, d\xi,\quad u\in L^2(\R),$$
and for $\theta\in\R$ with $\sin(\theta)\neq0$,
$${\CM}(R_\theta) u(x)= \frac{{\rm i}^{\frac12}}{\sqrt{2\pi \sin(\theta)}}\int_{\R} e^{\frac{\rm i}{\sin(\theta)}\left(\frac{x^2+y^2}{2}\cos(\theta)+xy\right)} u(y) \, dy,\quad u\in L^2(\R).$$
In particular, ${\CM}(I)={\CM}(R_0)={\rm Id}.$, ${\CM}(J)={\CM}(R_{\frac{\pi}{2}})={\rm i}^{\frac12} {\CF}^{-1}$.

\smallskip

\item [(ii)] For $\kappa\in\R$,
$$ {\CM}\left(\begin{array}{cc}
               1 & 0 \\
               \kappa & 1
             \end{array}\right) u (x)= e^{\frac{\rm i}{2}\kappa x^2} u(x),\quad u\in L^2(\R). $$

\smallskip

\item [(iii)] For $\lambda\in\R\setminus\{0\}$,
$$ {\CM}\left(\begin{array}{cc}
               \lambda & 0 \\
               0 & \lambda^{-1}
             \end{array}\right) u (x)= \lambda^{-\frac12} u(\lambda^{-1}x),\quad u\in L^2(\R).  $$

\smallskip

\item [(iv)] For $A,B\in\SL$, ${\CM}(AB)={\CM}(A)\CM(B)$ \footnote{According to the double-value ambiguity mentioned in Remark \ref{rmk_signs}, it should be stated as ${\CM}(AB)=\pm{\CM}(A)\CM(B)$. Here we neglect the sign.}. In particular, ${\CM}(A^{-1})={\CM}(A)^{-1}$.

\smallskip

\item [(v)] Given any $Y\in \SL$ and $L\in\sl$, we have
\begin{equation}\label{conj_quadratic}
{\CM}(Y)^{-1} \, \frac{-1}2 \left\la \left(\begin{array}{c}
X\\
D
 \end{array}\right) , JL\left(\begin{array}{c}
X\\
D
 \end{array}\right) \right\ra {\CM}(Y)=\frac{-1}2 \left\la \left(\begin{array}{c}
X\\
D
 \end{array}\right) , JY^{-1}LY \left(\begin{array}{c}
X\\
D
 \end{array}\right) \right\ra.
\end{equation}

\smallskip

\item [(vi)] Given $L\in \sl$, for the quadratic quantum Hamiltonian
$${\CL}={\CL}(X,D):=\frac{-1}2 \left\la \left(\begin{array}{c}
X\\
D
 \end{array}\right) , JL \left(\begin{array}{c}
X\\
D
 \end{array}\right) \right\ra,$$
we have that $e^{-{\rm i}t{\CL}}={\CM}(e^{tL})$ for every $t\in\R$.

\smallskip

\item [(vii)]
Given two $\sl-$linear systems
\begin{equation}\label{linear_systems}
\dot{y}_j(t)=L_j(t)y_j(t), \qquad L_j(t)\in \sl, \quad t\in\R, \quad j=1,2,
\end{equation}
and two quadratic quantum Hamiltonians
\begin{equation}\label{quantum_Hamiltonians}
{\rm i}\partial_t\psi_j=\CL_j(t)\psi_j=\frac{-1}2 \left\la \left(\begin{array}{c}
X\\
D
 \end{array}\right) , J L_j(t) \left(\begin{array}{c}
X\\
D
 \end{array}\right) \right\ra \psi_j,\quad j=1,2,\end{equation}
if the change of variables
\begin{equation}\label{conj_ode}
y_1(t)= Y(t) \,  y_2(t),\qquad Y(t)\in C^1(\R, \SL),
\end{equation}
gives the conjugation between linear systems in (\ref{linear_systems}), then
the  $L^2-$unitary transformation defined by the Meta representation
\begin{equation}\label{conj_pde}
\psi_1(t)={\CM}(Y(t)) \psi_2(t), \quad t\in\R ,
\end{equation}
 gives the conjugation between the two quantum Hamiltonians in (\ref{quantum_Hamiltonians}).
\end{enumerate}
\end{prop}

With the general definition of Meta representation, some of the above properties have been stated and shown in \cite{Fol1989}. Only the assertions of Proposition \ref{prop_elem_meta} will be used in the sequel. For completeness, we give a proof of this proposition in Appendix \ref{app_pr_prop}.

\subsection{Estimates of ${\CH}^s-$norms for quantum propagator}

With the Meta representation, we have a correspondence between the classical linear systems and quadratic quantum Hamiltonians.
We are interested in the behavior of Sobolev norms for the quadratic quantum Hamiltonians through the Meta representation of a given $\SL-$matrix. In the remaining of this section, we give several estimates on Sobolev norms, according to the elements of a given $\SL-$matrix.

For $L\in{\rm sl}(2,\R)$,
let us consider the 1D partial differential equation
\begin{equation}\label{PDE_quadric}
{\rm i}\partial_t u(t,x)=({\CL}u)(t,x)=\frac{-1}2 \left\la \left(\begin{array}{c}
X\\
D
 \end{array}\right) , JL \left(\begin{array}{c}
X\\
D
 \end{array}\right) \right\ra u(t,x),\quad u(0,\cdot)=u_0(\cdot)\in {\CH}^s.
\end{equation}
According to Proposition \ref{prop_elem_meta}-(vi), we see that, $e^{-{\rm i}t{\CL}}=\CM(e^{tL})$, i.e., the propagator $e^{-{\rm i}t{\CL}}$ of Eq. (\ref{PDE_quadric}) is the Meta representation of the classical flow $e^{tL}$.
Moreover, the well-posedness of quadratic quantum Hamiltonians implies that $e^{-{\rm i}t{\CL}}\in{\CB}(\CH^s)$ for any $s\geq 0$.

Let us assume that the determinant of $L=\left(
\begin{array}{cc}
  a_{11} & a_{02} \\
  -a_{20} & -a_{11}
\end{array}
\right)$ satisfies that $a_{20}a_{02}-a_{11}^2>0$ and define
$$\varrho:=\sqrt{a_{20}a_{02}-a_{11}^2}.$$ We have the following upper bound of Sobolev norms on the propagator $e^{-{\rm i}t{\CL}}$.
\begin{prop}\label{cor_metaplectic_upper}
For $s\geq 0$, we have
\begin{equation}\label{esti_eitL_ellip}
\|e^{-{\rm i}t{\CL}}\|_{{\CB}({\CH}^s)}\lesssim_s 1+ \left(|a_{20}|+|a_{11}|+|a_{02}|\right)^s\left|\frac{\sin(t\varrho)}{\varrho}\right|^s,\quad t\in\R.
\end{equation}
\end{prop}
\begin{remark} The factor $\left|\frac{\sin(\varrho t)}{\varrho}\right|^s$ in the upper bound (\ref{esti_eitL_ellip}) can be bounded by $t^s$ or $\varrho^{-s}$ as we need.
\end{remark}

Assume further that $a_{11}=0$ in $L$, i.e., $L=\left(
\begin{array}{cc}
  0 & a_{02} \\
  -a_{20} & 0
\end{array}
\right)$. Recall that, for non-vanishing $u\in{\CH}^s$, we defined in Section \ref{sec_Nota},
$$
\Gamma_s(u)=\frac{\|u\|_{{\CH}^s}}{\|D^su\|_{L^2}}.
$$
We have the following lower bound of Sobolev norms.

\begin{prop}\label{cor_metaplectic_lower}
For $s\geq 0$ and non-vanishing $u\in {\CH}^s$, if $\frac{|a_{02}|\sin(\varrho t)}{\varrho (1+\Gamma_s(u))^2}$ is large enough (depending only on $s$), then
\begin{equation}\label{esti_eitL_ellip-lower}
\|e^{-{\rm i}t\CL}u\|_{{\CH}^s}\gtrsim_{s}  |a_{02}|^s\left|\frac{\sin(\varrho t)}{\varrho}\right|^s \|D^su\|_{L^2}+ |a_{02}|^{-s}\left|\frac{\sin(\varrho t)}{\varrho}\right|^{-s}\|X^s u\|_{L^2}.
\end{equation}
\end{prop}

Proposition \ref{cor_metaplectic_upper} and \ref{cor_metaplectic_lower} will be proved in the remaining of this section, by direct calculations on the general Meta representation of $A\in \SL$.

\subsection{Bounds of Sobolev norms -- basic calculations}

Before studying the ${\CH}^s-$norm by the Meta representation, let us give the following basic lemma about upper and lower bounds of Sobolev norms for the function
$$U_{{\CG}, \, {\CK}}(x):=e^{{\rm i}{\CG} x^2} u({\CK} x),\quad {\CG}\in\R, \quad {\CK}\in\R\setminus\{0\},\quad u\in{\CH}^s,$$
since, in view of (\ref{Meta_a}) and (\ref{Meta_b}), every Meta representation can be presented in this form.

\begin{prop}\label{lemma_Hs}
For every ${\CG}\in\R$ and ${\CK}\in\R\setminus\{0\}$, we have
 \begin{equation}\label{upper}
\|U_{{\CG}, \, {\CK}}\|_{\CH^s}\lesssim_s |{\CK}|^{-\frac12} \left(|{\CK}|+|{\CG}{\CK}^{-1}|+|{\CK}|^{-1}\right)^{s} \left\|u\right\|_{{\CH}^s},\quad \forall \ u\in{\CH}^s.
 \end{equation}
 Moreover, for non-vanishing $u\in{\CH}^s$,
\begin{itemize}
  \item If $|{\CG}{\CK}^{-2}|(1+\Gamma_s(u))^2$ is small enough (depending only on $s$),
  then we have
 \begin{equation}\label{lower-j}
\|U_{{\CG}, \, {\CK}}\|_{{\CH}^s} \gtrsim_{s} |{\CK}|^{s-\frac12} \|D^su\|_{L^2}+ |{\CK}|^{-s-\frac12}\|X^s u\|_{L^2}.
\end{equation}
  \item If $|{\CG}{\CK}^{-2}|\max\{1, |{\CK}|^{4s}\}$ is small enough (depending only on $s$),
   then we have
  \begin{equation}\label{lower_transf}
  \|U_{{\CG}, \, {\CK}}\|_{\CH^s} \gtrsim_{s}  |{\CK}|^{-\frac12}\min\{|{\CK}|^{s}, |{\CK}|^{-s}\} \left\|u\right\|_{\CH^s}.
  \end{equation}
  In particular, if $|\CG|$ is small enough (depending only on $s$), then
 \begin{equation}\label{lower-K=1}
\|U_{{\CG}, \, 1}\|_{\CH^s} \gtrsim_{s}\left\|u\right\|_{\CH^s}.
\end{equation}
\end{itemize}
\end{prop}

At first, we have the following property for derivatives of the function $e^{{\rm i}{\CG}x^2}$, $\CG\in\R$.

\begin{lemma}\label{derivative_e_x2} For $\alpha\in \N^*$, we have
\begin{equation}\label{expre_deriv}
\frac{d^\alpha}{dx^\alpha}\left(e^{{\rm i}{\CG}x^2}\right)=e^{{\rm i}{\CG}x^2} \sum_{k=0}^{\lfloor\frac{\alpha}{2}\rfloor} p_{\alpha, k} \cdot  (2{\rm i}{\CG})^{\alpha-k}x^{\alpha-2k}
\end{equation}
with the coefficients satisfying $p_{\alpha, 0}=1$ and for $k\geq 1$ (if existing),
\begin{equation}\label{esti_coeff}
p_{\alpha, k}\geq 1,\quad \sum_{k=1}^{\lfloor\frac{\alpha}{2}\rfloor} p_{\alpha, k}\leq \alpha!.
\end{equation}
\end{lemma}
\proof We can compute directly the derivatives until the third order:
\begin{eqnarray*}
\frac{d}{dx}\left(e^{{\rm i}{\CG}x^2}\right)(x)&=&2{\rm i}{\CG}x e^{{\rm i}{\CG}x^2}, \\
\frac{d^2}{dx^2}\left(e^{{\rm i}{\CG}x^2}\right)(x)&=&((2{\rm i}{\CG} x)^2+2{\rm i}{\CG}) e^{{\rm i}{\CG}x^2}, \\
\frac{d^3}{dx^3}\left(e^{{\rm i}{\CG}x^2}\right)(x)&=&((2{\rm i}{\CG} x)^3+3(2{\rm i}{\CG})^2x) e^{{\rm i}{\CG}x^2}.
\end{eqnarray*}
Then we have (\ref{expre_deriv}) and (\ref{esti_coeff}) for $\alpha=1,2,3$.
Assume that for some $\alpha\geq 3$, we have
$$\frac{d^\alpha}{dx^\alpha}\left(e^{{\rm i}{\CG}x^2}\right)(x)=e^{{\rm i}{\CG}x^2}\left((2{\rm i}{\CG})^\alpha x^\alpha + \sum_{k=1}^{\lfloor\frac{\alpha}{2}\rfloor} p_{\alpha, k} \cdot (2{\rm i}{\CG})^{\alpha-k}x^{\alpha-2k}\right).$$
Then, through a direct calculation, we have
\begin{eqnarray*}
\frac{d^{\alpha+1}}{dx^{\alpha+1}}\left(e^{{\rm i}{\CG}x^2}\right)(x)&=&\frac{d}{dx}\left(\frac{d^\alpha}{dx^\alpha}e^{{\rm i}{\CG}x^2}\right)\\
&=&e^{{\rm i}{\CG}x^2}\left(2{\rm i}{\CG}x\left((2{\rm i}{\CG})^\alpha x^\alpha + \sum_{k=1}^{\lfloor\frac{\alpha}{2}\rfloor} p_{\alpha, k}  \cdot (2{\rm i}{\CG})^{\alpha-k}x^{\alpha-2k}\right)\right.\\
& & \  \  \  \  \  \  \  \  \  \ + \left. \, \left(\alpha (2{\rm i}{\CG})^\alpha x^{\alpha-1} +\sum_{k=1}^{\lfloor\frac{\alpha}{2}\rfloor}  p_{\alpha, k}  \cdot  (\alpha-2k) (2{\rm i}{\CG})^{\alpha-k}x^{\alpha-2k-1}\right)\right)\\
&=&e^{{\rm i}{\CG}x^2}\left((2{\rm i}{\CG})^{\alpha+1} x^{\alpha+1} + \sum_{k=1}^{\lfloor\frac{\alpha+1}{2}\rfloor} p_{\alpha+1, k} (2{\rm i}{\CG})^{\alpha+1-k}x^{\alpha+1-2k}\right)
\end{eqnarray*}
where the coefficients $p_{\alpha+1, k}$, $k=1,\cdots \left\lfloor\frac{\alpha+1}{2}\right\rfloor$, are
\begin{itemize}
  \item $p_{\alpha+1, 1}=\alpha+p_{\alpha, 1}$
  \item if $\alpha$ is odd, then
  $$p_{\alpha+1, k}=\left\{\begin{array}{cl}
                           p_{\alpha, k}+(\alpha-2k+2)p_{\alpha, k-1},   & 2\leq k \leq \frac{\alpha-1}{2}=\left\lfloor\frac{\alpha}{2}\right\rfloor \\
                             (\alpha-2k+2)p_{\alpha, k-1}, & k=\frac{\alpha+1}{2}=\left\lfloor\frac{\alpha+1}{2}\right\rfloor
                           \end{array}
  \right. ,$$
  and if $\alpha$ is even, then
  $$p_{\alpha+1, k}=p_{\alpha, k}+(\alpha-2k+2)p_{\alpha, k-1}, \quad 2\leq k \leq \frac{\alpha}{2}=\left\lfloor\frac{\alpha+1}{2}\right\rfloor.$$
\end{itemize}
Hence we obtain (\ref{expre_deriv}) for any $\alpha\geq 3$. It is easy to see from the above recurrence that all coefficients $p_{\alpha,k}\geq 1$.
Moreover, if $\alpha$ is odd, then
  $$   \sum_{k=1}^{\frac{\alpha+1}{2}} p_{\alpha+1, k}=\alpha+ \sum_{k=1}^{\frac{\alpha-1}{2}}p_{\alpha, k}+ \sum_{k=2}^{\frac{\alpha+1}{2}} (\alpha-2k+2)p_{\alpha, k-1}
    \leq   \alpha+\alpha!+ (\alpha-2) \alpha! <(\alpha+1)!, $$
and if $\alpha$ is even, then
   $$   \sum_{k=1}^{\frac{\alpha+1}{2}} p_{\alpha+1, k}=\alpha+ \sum_{k=1}^{\frac{\alpha}{2}}p_{\alpha, k}+ \sum_{k=2}^{\frac{\alpha}{2}} (\alpha-2k+2)p_{\alpha, k-1}
    \leq   \alpha+\alpha!+ (\alpha-2) \alpha! <(\alpha+1)! .$$
The inequality (\ref{esti_coeff}) is shown for $\alpha\geq 3$.\qed

\medskip

\noindent
{\bf Proof of Proposition \ref{lemma_Hs}.} In view of the equivalent definition of the ${\CH}^s-$norm in (\ref{norm_equiv-1}), to estimate $\|U_{{\CG}, \, {\CK}}\|_{\CH^s}$, it is sufficient to calculate the $L^2-$norms
$\left\|X^sU_{{\CG}, \, {\CK}}\right\|_{L^2}$ and $\left\|D^s U_{{\CG}, \, {\CK}}\right\|_{L^2}$ assuming that $s\in \N^*$.

By a direct calculation, we have that
\begin{equation}\label{norm_xsf}
\|X^s U_{{\CG}, \, {\CK}}\|_{L^2}=|{\CK}|^{-s-\frac12}\|X^s u\|_{L^2}.
\end{equation}
For the $s^{\rm th}$-order derivative of $U_{{\CG}, \, {\CK}}$, by Lemma \ref{derivative_e_x2}, we have that
\begin{eqnarray}
& &\left(D^s  U_{{\CG}, \, {\CK}}\right)(x)\nonumber\\
&=&(-{\rm i})^s\sum_{\alpha=0}^s C_s^{\alpha}  \left(e^{{\rm i}{\CG}x^2}\right)^{(\alpha)}\left(u({\CK}x)\right)^{(s-\alpha)}\label{form_upper}\\
&=&{\CK}^s e^{{\rm i}{\CG}x^2}( D^su)({\CK} x) +  (-{\rm i})^s e^{{\rm i}{\CG}x^2} \sum_{\alpha=1}^s C_s^{\alpha} {\CK}^{s-\alpha} \sum_{k=0}^{\lfloor\frac{\alpha}{2}\rfloor} p_{\alpha, k}  (2{\rm i}{\CG})^{\alpha-k}x^{\alpha-2k}u^{(s-\alpha)}({\CK} x).\label{form_lower}
\end{eqnarray}
According to (\ref{form_upper}), we obtain
\begin{eqnarray*}
\left\|D^s  U_{{\CG}, \, {\CK}}\right\|_{L^2}&=&\left\|\sum_{\alpha=0}^s C_s^{\alpha} {\CK}^{s-\alpha} \sum_{k=0}^{\lfloor\frac{\alpha}{2}\rfloor} p_{\alpha, k} \cdot  (2{\rm i}{\CG})^{\alpha-k}x^{\alpha-2k}u^{(s-\alpha)}({\CK} x) \right\|_{L^2}\\
&\leq&\sum_{\alpha=0}^s C_s^{\alpha} |{\CK}|^{s-\alpha} \sum_{k=0}^{\lfloor\frac{\alpha}{2}\rfloor} p_{\alpha, k}|2{\CG}|^{\alpha-k} \|x^{\alpha-2k}u^{(s-\alpha)}({\CK} x)\|_{L^2},
\end{eqnarray*}
where, similar to (\ref{norm_xsf}), we have
$$\|x^{\alpha-2k}u^{(s-\alpha)}({\CK} x)\|_{L^2}=|{\CK}|^{-(\alpha-2k)-\frac12}\|x^{\alpha-2k}u^{(s-\alpha)}( x)\|_{L^2}. $$
Hence, by the equivalent definition (\ref{norm_equiv}) of ${\CH}^s-$norm, we have
\begin{eqnarray}
\left\|D^s U_{{\CG}, \, {\CK}}\right\|_{L^2}
&\leq&|{\CK}|^{s-\frac12}\sum_{\alpha=0}^s C_s^{\alpha} \sum_{k=0}^{\lfloor\frac{\alpha}{2}\rfloor} p_{\alpha, k}|2{\CG}|^{\alpha-k}|{\CK}|^{-2(\alpha-k)} \|x^{\alpha-2k}  u^{(s-\alpha)}(x)\|_{L^2} \label{norm_parsf}\\
&\lesssim_s&|{\CK}|^{s-\frac12} \left\|u\right\|_{{\CH}^s}\sum_{\alpha=0}^s C_s^{\alpha} \left(1+2|{\CG}{\CK}^{-2}|\right)^{\alpha}\nonumber\\
&=&|{\CK}|^{s-\frac12} \left\|u\right\|_{{\CH}^s}\left(2+2|{\CG}{\CK}^{-2}|\right)^s\nonumber\\
&=& 2^s |{\CK}|^{-\frac12}\left(|{\CK}|+|{\CG}{\CK}^{-1}|\right)^{s} \left\|u\right\|_{{\CH}^s}.\nonumber
\end{eqnarray}
Combining (\ref{norm_xsf}) and (\ref{norm_parsf}), we obtain (\ref{upper}).

According to (\ref{form_lower}), we see that
\begin{eqnarray*}
& & \left\|(D^s U_{{\CG}, \, {\CK}})(x)- {\CK}^s e^{{\rm i}{\CG}x^2} (D^s u)( {\CK} x) \right\|_{L^2}\\
&=&\left\|\sum_{\alpha=1}^s C_s^{\alpha} {\CK}^{s-\alpha} \sum_{k=0}^{\lfloor\frac{\alpha}{2}\rfloor} p_{\alpha, k} (2{\rm i}{\CG})^{\alpha-k}x^{\alpha-2k}u^{(s-\alpha)}({\CK}x) \right\|_{L^2}\\
&\leq&|{\CK}|^{s-\frac12}\sum_{\alpha=1}^s C_s^{\alpha} \sum_{k=0}^{\lfloor\frac{\alpha}{2}\rfloor} p_{\alpha, k}|2{\CG}|^{\alpha-k}|{\CK}|^{-2(\alpha-k)}\|x^{\alpha-2k} u^{(s-\alpha)}(x)\|_{L^2},
\end{eqnarray*}
where, by the smallness assumption on $|{\CG}{\CK}^{-2}|$, we have that
$$\sum_{k=0}^{\lfloor\frac{\alpha}{2}\rfloor} p_{\alpha, k}|2{\CG}|^{\alpha-k}|{\CK}|^{-2(\alpha-k)} \leq (1+\alpha!) |2{\CG}{\CK}^{-2}|^{\alpha-\lfloor\frac{\alpha}{2}\rfloor}\leq (s!+1) |2{\CG}{\CK}^{-2}|^{\frac{\alpha}{2}}.$$
Then we obtain
\begin{eqnarray}
\left\|(D^s U_{{\CG}, \, {\CK}})(x)- {\CK}^s e^{{\rm i}{\CG}x^2} (D^s u)( {\CK} x) \right\|_{L^2}
&\lesssim_s&  |{\CK}|^{s-\frac12} \left\|u\right\|_{{\CH}^s}\sum_{\alpha=1}^s C_s^{\alpha}|2{\CG}{\CK}^{-2}|^{\frac\alpha 2}\label{esti_JK-}\\
&=&  |{\CK}|^{s-\frac12} \left(\left(1+|2{\CG}{\CK}^{-2}|^\frac12\right)^{s}-1 \right) \left\|u\right\|_{{\CH}^s}\nonumber\\
&\lesssim_s&  |{\CK}|^{s-\frac12} |{\CG}{\CK}^{-2}|^\frac12\left\|u\right\|_{{\CH}^s}.\nonumber
\end{eqnarray}
Since $|{\CG}{\CK}^{-2}|(1+\Gamma_s(u))^2$ is sufficiently small, through (\ref{esti_JK-}), we have that
$$ \left\|(D^s U_{{\CG}, \, {\CK}})(x)-{\CK}^s e^{{\rm i}{\CG}x^2}(D^s u)({\CK} x)\right\|_{L^2}\lesssim_s
|{\CK}|^{s-\frac12} \cdot \frac{\|u\|_{{\CH}^s}}{10 \Gamma_s(u)}=
\frac{1}{10}|{\CK}|^{s-\frac12} \|D^su\|_{L^2}.$$
Then, we have
\begin{eqnarray*}
\|(D^s U_{{\CG}, \, {\CK}})\|_{L^2}
&\geq&   \left\|{\CK}^s e^{{\rm i}{\CG}x^2} (D^s u)( {\CK} x) \right\|_{L^2} -  \left\|(D^s U_{{\CG}, \, {\CK}})(x)-{\CK}^s e^{{\rm i}{\CG}x^2} (D^s u)( {\CK} x)\right\|_{L^2}\\
 &\gtrsim_{s}&\frac{9}{10}|{\CK}|^{s-\frac12} \left\|D^s u\right\|_{L^2}.
\end{eqnarray*}
Combining with (\ref{norm_xsf}), we obtain (\ref{lower-j}).

Moreover, if $|{\CG}{\CK}^{-2}|^{\frac12}\max\{1 , |{\CK}|^{2s} \}$ is sufficiently small,
then (\ref{esti_JK-}) implies that
$$ \left\|(D^s U_{{\CG}, \, {\CK}})(x)- {\CK}^s e^{{\rm i}{\CG}x^2} (D^s u)( {\CK} x) \right\|_{L^2} \lesssim_s \frac{1}{10}|{\CK}|^{s-\frac12}\min\{|{\CK}|^{-2s}, 1\} \left\|u\right\|_{{\CH}^s}.$$
Then, we have
\begin{eqnarray*}
& & \|X^s U_{{\CG}, \, {\CK}}\|_{L^2} +\|D^s U_{{\CG}, \, {\CK}}\|_{L^2}\\
&\geq& |{\CK}|^{-s-\frac12}\|X^s u\|_{L^2}+  \left\|{\CK}^s e^{{\rm i}{\CG}x^2} (D^s u)( {\CK} x) \right\|_{L^2} - \,  \left\|(D^s U_{{\CG}, \, {\CK}})(x)-{\CK}^s e^{{\rm i}{\CG}x^2} (D^s u)( {\CK} x)\right\|_{L^2}\\
 &\gtrsim_{s} &|{\CK}|^{-s-\frac12}\|X^s u\|_{L^2}+ |{\CK}|^{s-\frac12} \left\|D^s u\right\|_{L^2} - \frac{1}{10}|{\CK}|^{s-\frac12}\min\{|{\CK}|^{-2s}, 1\} \left\|u\right\|_{{\CH}^s} \\
 &\gtrsim_s & \frac{9}{10} |{\CK}|^{-\frac12}\min\{|{\CK}|^{-s}, |{\CK}|^{s}\} \left\|u\right\|_{{\CH}^s},
\end{eqnarray*}
which implies (\ref{lower_transf}). In particular, if ${\CK}=1$, then (\ref{lower_transf}) implies (\ref{lower-K=1}) provided that $|{\CG}|$ is sufficiently small.
\qed

\subsection{Upper bound of Meta representation}

\begin{lemma}\label{lemma_esti_metaplectic}
Given $s\geq 0$ and $A\in \SL$, we have that ${\CM}(A):{\CH}^s\to{\CH}^s$ with
\begin{equation}\label{esti_meta_upper}
\|{\CM}(A)\|_{{\CB}({\CH}^s)}\lesssim_s \|A\|^s.
\end{equation}
\end{lemma}

\proof Let $A=\left(
\begin{array}{cc}
  a & b \\
  c & d
\end{array}
\right)\in {\rm SL}(2,\R)$. There are four possibly overlapped cases about the matrix elements:
\begin{itemize}
  \item [1)] $|ad|\geq |bc|$ with $|a|\geq \frac{1}{\sqrt{2}}$,
  \item [2)] $|ad|\geq |bc|$ with $|d|\geq \frac{1}{\sqrt{2}}$,
  \item [3)] $|ad|\leq |bc|$ with $|b|\geq \frac{1}{\sqrt{2}}$,
  \item [4)] $|ad|\leq |bc|$ with $|c|\geq \frac{1}{\sqrt{2}}$.
\end{itemize}

In {\rm Case 1)}, according to (\ref{Meta_a}), we have, for any $u\in L^2(\R)$,
$$({\CM}(A)u)(x)=(2\pi a)^{-\frac12}e^{\frac{\rm i}2 ca^{-1}x^2}\int_\R e^{{\rm i}a^{-1}x\xi}e^{-\frac{\rm i}2 ba^{-1}\xi^2} \widehat{u}(\xi)  d\xi.$$
By a change of variable $\xi\mapsto a\xi$ in the above integral, we obtain
\begin{equation}\label{meta-a_cv}
({\CM}(A)u)(x)
= (2\pi )^{-\frac12} a^{\frac12} e^{\frac{\rm i}2 ca^{-1}x^2} \int_\R e^{{\rm i} x \xi} e^{-\frac{\rm i}2  ab {\xi}^2} \widehat{u}(a\xi) \, d\xi.
\end{equation}
Applying (\ref{upper}) in Proposition \ref{lemma_Hs} with ${\CG}=\frac{1}2 ca^{-1}$ and ${\CK}=1$, we obtain, for any $u\in {\CH}^s$,
\begin{eqnarray}
\|{\CM}(A)u\|_{{\CH}^s}&\lesssim_s&  |a|^{\frac12} \left(1+|ca^{-1}|\right)^s\left\|\int_\R e^{{\rm i} x \xi} e^{-\frac{\rm i}2 ab {\xi}^2}\widehat{u}(a\xi) \, d\xi\right\|_{{\CH}^s} \label{esti_Hs_app-1} \\
&\sim_s&  |a|^{\frac12} \left(1+|ca^{-1}|\right)^s \|  e^{-\frac{\rm i}2 ab {\xi}^2} \widehat{u}(a\xi)\|_{{\CH}^s}.\nonumber
\end{eqnarray}
Then, applying (\ref{upper}) with ${\CG}=-\frac12 ab$ and ${\CK}=a$, we get
\begin{equation}\label{esti_Hs_app-2}
\| e^{-\frac{\rm i}2 ab {\xi}^2} \widehat u(a\xi)\|_{{\CH}^s}\lesssim_s  |a|^{-\frac12}\left(|a|+|a|^{-1}+|b|\right)^s\|u\|_{{\CH}^s}\lesssim_s |a|^{-\frac12}\left(|a|+|b|\right)^s\|u\|_{{\CH}^s},
\end{equation}
in view of the assumption that $|a|\geq \frac{1}{\sqrt{2}}$.
Combining (\ref{esti_Hs_app-1}) and (\ref{esti_Hs_app-2}) and recalling that $|ad|\geq |bc|$, we get
$$\|{\CM}(A)u \|_{{\CH}^s}\lesssim_s   \left(1+|ca^{-1}|\right)^s\left(|a|+|b|\right)^s \|u\|_{{\CH}^s}
\leq\left(|a|+|b|+|c|+|d|\right)^s \|u\|_{{\CH}^s},$$
which implies (\ref{esti_meta_upper}).

\smallskip

In {\rm Case 3)}, according to (\ref{Meta_b}), we have, for any $u\in L^2(\R)$,
$$({\CM}(A)u)(x)={\rm i}^{\frac12}(2\pi)^{-\frac12}b^{-\frac12}e^{\frac{\rm i}2 db^{-1} x^2} \int_\R e^{{\rm i}b^{-1}xy}e^{\frac{\rm i}2b^{-1}a y^2} u(y)  dy.$$
By a change of variable $y\mapsto b y$ in the above integral, we get
\begin{equation}\label{meta-b_cv}
({\CM}(A)u)(x)
={\rm i}^{\frac12}(2\pi)^{-\frac12}b^{\frac12}e^{\frac{\rm i}2 db^{-1} x^2}\int_\R e^{{\rm i}xy}e^{\frac{\rm i}2ab y^2} u(by)  dy
\end{equation}
Applying (\ref{upper}) with ${\CG}=\frac12 db^{-1}$ and ${\CK}=1$, we obtain, for any $u\in {\CH}^s$,
\begin{eqnarray}
\|{\CM}(A)u\|_{{\CH}^s}&\lesssim_s&  |b|^{\frac12} \left(1+|db^{-1}|\right)^s\left\|\int_\R e^{{\rm i} x y} e^{\frac{\rm i}2 ab y^2}u(by) \, dy\right\|_{{\CH}^s}\nonumber  \\
&\lesssim_s&  |b|^{\frac12} \left(1+|db^{-1}|\right)^s \left\|  e^{\frac{\rm i}2 ab y^2} u(by)\right\|_{{\CH}^s}.\label{esti_Hs_app-3}
\end{eqnarray}
Then, applying (\ref{upper}) with ${\CG}=\frac12 ab$ and ${\CK}=b$, we get
\begin{equation}\label{esti_Hs_app-4}
\left\|  e^{\frac{\rm i}2 ab y^2} u(by)\right\|_{{\CH}^s}\lesssim_s  |b|^{-\frac12}\left(|b|+|b|^{-1}+|a|\right)^s\|u\|_{{\CH}^s}\lesssim_s |b|^{-\frac12}\left(|a|+|b|\right)^s\|u\|_{{\CH}^s},
\end{equation}
in view of the assumption that $|b|\geq \frac{1}{\sqrt{2}}$.
Combining (\ref{esti_Hs_app-3}) and (\ref{esti_Hs_app-4}) and recalling that $|ad|< |bc|$, we get
$$\|{\CM}(A)u\|_{{\CH}^s}\lesssim_s   \left(1+|db^{-1}|\right)^s\left(|a|+|b|\right)^s \|u\|_{{\CH}^s}
\leq\left(|a|+|b|+|c|+|d|\right)^s \|u\|_{{\CH}^s},$$
which implies (\ref{esti_meta_upper}).

\smallskip

As for {\rm Case 2)} and {\rm Case 4)}, since we see that
$$A= \left(
\begin{array}{cc}
  a & b \\
  c & d
\end{array}
\right)=-J \left(
\begin{array}{cc}
  c & d \\
   -a & -b
\end{array}
\right),$$
and, by Proposition \ref{prop_elem_meta},
$${\CM}(-J)={\CM}(J)^{-1}=\frac{1}{{\rm i}^{\frac12}}{\CF},$$
we have
$${\CM}(A)= {\CM}(-J) {\CM}\left(
\begin{array}{cc}
  c & d \\
   -a & -b
\end{array}
\right)=\frac{1}{{\rm i}^{\frac12}}{\CF} {\CM}\left(
\begin{array}{cc}
  c & d \\
   -a & -b
\end{array}
\right).$$
Hence we obtain, through {\rm Case 1)} and {\rm Case 3)},
$$\|{\CM}(A)\|_{{\CB}({\CH}^s)}\sim_s \left\|{\CM}\left(
\begin{array}{cc}
  c & d \\
   -a & -b
\end{array}
\right)\right\|_{{\CB}({\CH}^s)}\lesssim_s \|A\|^s.\qed$$

Now we are ready to show the upper bound of Sobolev norms for the propagator $e^{-{\rm i}t\CL}$.

\noindent
{\bf Proof of Proposition \ref{cor_metaplectic_upper}.} Recall that $L=\left(
\begin{array}{cc}
  a_{11} & a_{02} \\
  -a_{20} & -a_{11}
\end{array}
\right)$ with the determinant
$$a_{02}a_{20}-a_{11}^2>0.$$
 If $a_{02}>0$, by a direct computation, we have the normalization of $L$:
\begin{equation}\label{normal_L}
L=\frac{1}{\sqrt{a_{02}\varrho}}\left(\begin{array}{cc}
               a_{02} & 0\\[1mm]
               -a_{11} & \varrho
             \end{array}
\right) \cdot \left(\begin{array}{cc}
             0  & \varrho \\[1mm]
             -\varrho  & 0
             \end{array}
\right)\cdot \frac{1}{\sqrt{a_{02}\varrho}}\left(\begin{array}{cc}
               \varrho & 0\\[1mm]
               a_{11} & a_{02}
             \end{array}
\right).
\end{equation}
Then we see that
\begin{eqnarray}
e^{tL} &=& \frac{1}{\sqrt{a_{02}\varrho}}\left(\begin{array}{cc}
               a_{02} & 0\\[1mm]
               -a_{11} & \varrho
             \end{array}
\right)\cdot  \left(\begin{array}{cc}
             \cos(\varrho t)  & \sin(\varrho t) \\[1mm]
             -\sin(\varrho t)  & \cos(\varrho t)
             \end{array}
\right)  \cdot \frac{1}{\sqrt{a_{02}\varrho}}\left(\begin{array}{cc}
               \varrho & 0\\[1mm]
               a_{11} & a_{02}
             \end{array}
\right)\label{e_tL}\\
&=&\left(\begin{array}{cc}
             \cos(\varrho t) + \frac{a_{11}}{\varrho}\sin(\varrho t)  & \frac{a_{02}}{\varrho} \sin(\varrho t) \\[1mm]
             -\frac{a_{20}}{\varrho} \sin(\varrho t)  & \cos(\varrho t)-\frac{a_{11}}{\varrho}\sin(\varrho t)
             \end{array}
\right).\nonumber
\end{eqnarray}
Otherwise, if $a_{02}<0$, we have the normalization of $-L$ as in (\ref{normal_L}) with $-a_{02}>0$, and we obtain $e^{tL}=e^{(-t)(-L)}$, which is the same form with (\ref{e_tL}).
Then the upper bound (\ref{esti_eitL_ellip}) follows immediately from Lemma \ref{lemma_esti_metaplectic}.\qed

\subsection{Lower bounds of Meta representation}

We shall give two lower bounds of Sobolev norms through the Meta representation of
 $A=\left(
\begin{array}{cc}
  a & b \\
  c & d
\end{array}
\right)\in {\rm SL}(2,\R)$, under suitable assumptions on the matrix elements.

\begin{lemma}\label{lemma_metaplectic_lower-1}
Assume that $a\neq 0$. If $|ca^{-1}|$ and $|ba^{-1}|\max\{|a|^{4s}, 1\}$ are small enough (depending only on $s$), then
\begin{equation}\label{meta-lower-1}
\|{\CM}(A)u\|_{{\CH}^s} \gtrsim_{s} \min\{|a|^{s}, |a|^{-s} \}\|u\|_{{\CH}^s},\quad u\in{\CH}^s.
\end{equation}
Moreover, if $b=c=0$, then
\begin{equation}\label{meta-lower-1-deriv}
\|D^s{\CM}(A)u\|_{L^2}=|a|^{-s} \|D^s u\|_{L^2},\quad \|X^s{\CM}(A)u\|_{L^2}=|a|^{s} \|X^s u\|_{L^2}.
\end{equation}
\end{lemma}

\proof As in the proof of Lemma \ref{lemma_esti_metaplectic}, we have (\ref{meta-a_cv}) for ${\CM}(A)u$ under the assumption that $a\neq 0$.
Since $|ca^{-1}|$ is sufficiently small, we apply (\ref{lower-K=1}) in Proposition \ref{lemma_Hs} with ${\CG}=\frac12 ca^{-1}$ and ${\CK}=1$, and obtain
$$\|{\CM}(A)u\|_{{\CH}^s}\gtrsim_{s}|a|^{\frac12}\left\|\int_\R e^{{\rm i} x \xi} e^{-\frac{\rm i}2 ab {\xi}^2} \widehat{u}(a\xi) \, d\xi\right\|_{{\CH}^s}\sim_{s}|a|^{\frac12}\| e^{-\frac{\rm i}2 ab {\xi}^2} \widehat{u}(a\xi) \|_{{\CH}^s}.$$
Applying (\ref{lower_transf}) with ${\CG}=-\frac12 ab$ and ${\CK}=a$, we obtain
$$\left\|  e^{-\frac{\rm i}2 ab {\xi}^2} \widehat{u}(a\xi) \right\|_{{\CH}^s} \gtrsim_{s} |a|^{-\frac12} \min\{|a|^{s}, |a|^{-s} \}\|\widehat{u}\|_{{\CH}^s} \sim_{s} |a|^{-\frac12} \min\{|a|^{s}, |a|^{-s} \}\|u\|_{{\CH}^s},$$
since
$|ba^{-1}|\max\{|a|^{4s}, 1\}$ is small enough.
Then (\ref{meta-lower-1}) is shown.

If $A$ is, in particular, diagonal, i.e., $A=\left(
\begin{array}{cc}
  a & 0 \\
  0 & a^{-1}
\end{array}
\right)$, then, by the assertion iii) of Proposition \ref{prop_elem_meta}, we have ${\CM}(A)u (x)= a^{-\frac12} u(a^{-1} x)$.
We obtain (\ref{meta-lower-1-deriv}) through direct computations.\qed

\begin{lemma}\label{lemma_metaplectic_lower-2}
Given $s\geq 0$ and $u\in{\CH}^s$ non-vanishing, assume that $|ab^{-1}|(1+\Gamma_s(u))^2$ and $|db^{-1}|$ are sufficiently small. Then
\begin{equation}\label{meta-lower-2}
\|{\CM}(A)u\|_{{\CH}^s} \gtrsim_{s}  |b|^{s} \|D^su\|_{L^2}+ |b|^{-s}\|X^s u\|_{L^2}.
\end{equation}
\end{lemma}
\proof Under the assumption that $b\neq 0$, we have (\ref{meta-b_cv}) for $\CM(A)u$ as in the proof of Lemma \ref{lemma_esti_metaplectic}. Since $|db^{-1}|$ is small enough, we
apply (\ref{lower-K=1}) in Proposition \ref{lemma_Hs} with ${\CG}=\frac12 db^{-1}$ and ${\CK}=1$, and obtain
$$\|{\CM}(A)u\|_{{\CH}^s}\gtrsim_{s}|b|^{\frac12}\left\|\int_\R e^{{\rm i} x y} e^{\frac{\rm i}2 ab y^2}u(by) \, dy\right\|_{{\CH}^s}\sim_{s}|b|^{\frac12}\|  e^{\frac{\rm i}2 ab y^2} u(by)\|_{{\CH}^s}.$$
Applying (\ref{lower-j}) with ${\CG}=\frac12 ab$ and ${\CK}=b$, we have
$$
\|e^{\frac{\rm i}2 ab y^2} u(by)\|_{{\CH}^s}\gtrsim_{s}  |b|^{s-\frac12} \|D^s u\|_{L^2}+ |b|^{-s-\frac12} \|X^s u\|_{L^2}.
$$
since $|ab^{-1}|(1+\Gamma_s(u))^2$ is sufficient small.
Hence (\ref{meta-lower-2}) is shown.\qed

\

\noindent{\bf Proof of Proposition \ref{cor_metaplectic_lower}.} Under the assumptions of Proposition \ref{cor_metaplectic_lower}, we have $$L=\left(\begin{array}{cc}
  0 & a_{02} \\
  -a_{20} & 0\end{array}\right) \quad  {\rm with}  \  \  a_{20}a_{02}>0.$$
By (\ref{e_tL}), we obtain
$$
e^{tL}=\left(\begin{array}{cc}
             \cos(\varrho t) & \frac{a_{02}}{\varrho} \sin(\varrho t) \\[1mm]
             \frac{-a_{20}}{\varrho} \sin(\varrho t)  & \cos(\varrho t)
             \end{array}
\right).$$
Since $\left|\frac{a_{02}\sin(\varrho t)}{\varrho(1+\Gamma_s(u))^2}\right|$ is large enough, we have that
$\left|\frac{\varrho\cos(\varrho t)}{a_{02}\sin(\varrho t)}\right|(1+\Gamma_s(u))^2 $ is sufficiently small.
Then the lower bound (\ref{esti_eitL_ellip-lower}) follows from Lemma \ref{lemma_metaplectic_lower-2}.\qed

\section{Almost reducibility in quantum Hamiltonian systems}\label{sec_AR}

In this section, we shall give more descriptions about the almost reducibility scheme of Eq. (\ref{eq_original}) (formulated in Proposition \ref{prop_upper}), which is closely related to that of corresponding linear system. Such an argument leads to the proof of Theorem \ref{thm_al_re}.
We will also construct a concrete time quasi-periodic quadratic perturbation ${\CP}(t)$, for which the sequence of unitary transformation $\{{\CU}_j(t)\}$ in the almost reducibility scheme of Eq. (\ref{eq_original}) does not converge in some suitable sense (formulated in Proposition \ref{prop_lower}). This concrete perturbation leads to the oscillatory growth of Sobolev norm stated in Theorem \ref{thm_sub_polynome}, whose proof will be completed in Section \ref{sec_pr_thm}.
The proof of Proposition \ref{prop_upper} and \ref{prop_lower} will be given in Section \ref{sec_al_ODE} through the exact correspondence between classical and quantum quadratic Hamiltonians formulated by the Meta representation.

\subsection{General scheme of almost reducibility for Eq. (\ref{eq_original})}



Let  $(\om, J+P(\cdot))$ denote the quasi-periodic linear system $\dot{y}(t)=(J+P(\om t))y(t)$
with
$$J=\left(
\begin{array}{cc}
  0 & 1 \\
  -1 & 0
\end{array}
\right) ,\qquad P(\cdot)=\left(
\begin{array}{cc}
  p_{11}(\cdot) & p_{02}(\cdot) \\
  -p_{20}(\cdot) & -p_{11}(\cdot)
\end{array}
\right)\in C_r^\om(\T^d, {\sl}),$$
where the frequency vector $\omega$, as well as the coefficients  $p_{11}(\cdot)$, $p_{02}(\cdot)$, $p_{20}(\cdot)$, is the same with that in Eq. (\ref{eq_original}). Theorem \ref{thm_al_re} is shown through the following almost reducibility scheme.

\begin{prop}\label{prop_upper}
There exists $\varepsilon_*=\varepsilon_*(\gamma,\tau,r,d)>0$ such that, if $\|P\|_{r}<\varepsilon_*$, then the following holds.
\begin{enumerate}
\item Eq. (\ref{eq_original}) is almost reducible, i.e., there exist sequences $\{L_j\}_{j\in\N}\subset\sl$ and $\{P_j\}_{j\in\N^*}\subset C^\omega(\T^d,\sl)$ with $\|P_j\|_{\T^d} \to 0$
such that, Eq. (\ref{eq_original}) is conjugated to the equation
\begin{equation}\label{eq_conj_j}
{\rm i}\partial_t\psi_j = ({\CL}_j+{\CP}_{j+1}(t))\psi_j,
\end{equation}
via a time quasi-periodic unitary transformation $\psi(t)={\CU}_j(t)\psi_{j}(t)$ with ${\CU}_j(t)={\CM}(U_j(t))$ the Meta representation of
some $U_j(t)\in \SL$,
 where ${\CL}_j$ and ${\CP}_{j+1}(t)$ are defined as
 \begin{eqnarray}
{\CL}_j(X,D)&=& \frac{-1}2 \left\la \left(\begin{array}{c}
X\\
D
 \end{array}\right) , J L_j \left(\begin{array}{c}
X\\
D
 \end{array}\right) \right\ra,\label{L_j}\\
{\CP}_{j+1}(t, X,D)&=& \frac{-1}2 \left\la \left(\begin{array}{c}
X\\
D
 \end{array}\right) , J P_{j+1}(\om t) \left(\begin{array}{c}
X\\
D
 \end{array}\right) \right\ra.\label{P_j+1}
\end{eqnarray}

\item If Eq. (\ref{eq_original}) is non-reducible, then $\sup_{t\in\R} \|U_j(t)\|\to \infty$ as $j\to \infty$, and there is a sequence $\{\varepsilon_j\}_{j\in\N}$ with
\begin{equation}\label{sequence_epsilon}
0<\varepsilon_{j+1}< \varepsilon_j \exp\left\{-\frac{r}{2}\varepsilon_j^{-\frac{1}{18\tau}}\right\},
\end{equation}
such that
 \begin{equation}\label{l_small_rho}
 \sup_{t\in\R}\|U_j(t)\|\leq |\ln\varepsilon_{j}|^{2\tau}, \quad   \|L_j\| < \varepsilon_{j}^{\frac1{16}},\quad {\rm det}(L_j)>\frac{\gamma^2}{ |\ln\varepsilon_{j+1}|^{2\tau}}, \quad \|P_{j+1}\|_{\T^d} < \varepsilon_{j+1}.
  \end{equation}
\end{enumerate}

\end{prop}

\subsection{Concrete construction for non-reducible Eq. (\ref{eq_original})}\label{sec_concrete_PDE}

Recall $s> 0$ and $f:\R_+^*\to \R_+^*$ with $f(t)=o(t^s)$ given in Theorem \ref{thm_sub_polynome}.
Let us define $$g(t):=1-\frac{\ln(f(t))}{s\ln(t)}.$$
\begin{lemma}\label{lemma_gt}
We have $g(t)\in (0,1)$ for $t$ large enough and $t^{g(t)}\to \infty$ as $t\to\infty$.
\end{lemma}
\proof Calculating $s\ln(t)\cdot g(t)$ with
$$s\ln(t)\cdot g(t)=s\ln(t)- \ln(f(t))=\ln\left(\frac{t^s}{f(t)}\right)\to\infty,$$
we see that $t^{g(t)}=e^{\frac1s \cdot s\ln(t)g(t)}\to\infty$ as $t\to\infty$.\qed

\

Let us construct sequences $\{k_j\}_{j\in\N}\subset \Z^d$ and $\{T_j\}_{j\in\N}\subset \R_+$ according to the above $g(\cdot)$.

\begin{itemize}
\item Let $k_0, k_1\in \Z^d$ satisfy that
\begin{equation}\label{defi_k0k1}
\frac34<\la k_0\ra < 1,\quad 0<\la k_1\ra<\frac14.
\end{equation}
Let $T_0:=\frac{5\pi}{2\la k_0+k_1\ra}$, and, with defined $\{k_n\}_{0\leq n\leq j+1}$, $j\geq 1$, such that $\la k_n\ra>0$, let
\begin{equation}\label{defi_T_j}
T_j:=\frac{5\pi}{2\la k_{j+1} \ra}+\Xi_{j,j}+\cdots+\Xi_{j,1},
\end{equation}
where $-\frac{\pi}{\la k_0 + k_1 \ra}\leq\Xi_{j,1}\leq \frac{\pi}{\la k_0 + k_1\ra}$ and $-\frac{\pi}{\la k_{n} \ra}\leq\Xi_{j,n}\leq \frac{\pi}{\la k_{n}\ra}$ for $2\leq n \leq j$ such that
\begin{equation}\label{defi_T_j-error}
\frac{5\pi}{2\la k_{j+1}\ra}+\Xi_{j,j}+\cdots+\Xi_{j,n} \in \frac{2\pi}{\la k_n \ra} \Z,\quad T_j\in \frac{2\pi}{\la k_0 + k_1 \ra} \Z.
\end{equation}
  \item With defined $\{k_n\}_{0\leq n\leq j}$ and $\{T_n\}_{0\leq n\leq j-1}$, let $k_{j+1}\in \Z^d$ satisfy
  \begin{equation}\label{seq_k-1}
   |k_{j+1}|>e^{|k_j|}+10,
  \end{equation}
  and, with $r>0$ given in Eq. (\ref{eq_original}), with $T_{j}$ uniquely defined as in (\ref{defi_T_j}),
    \begin{equation}\label{seq_k-2}
    \la k_{j+1} \ra^{\frac{g(T_j)}{2}} < \la k_j\ra^{\frac{g(T_{j-1})}{2}} e^{-33 r  |k_j|} \exp\left\{-\la k_j\ra^{-\left(1+\frac{1}{36\tau}\right)}\right\}.
  \end{equation}
\end{itemize}
The Diophantine condition (\ref{dio}) of $\om$ means that $\{\la k\ra\}_{k\in\Z^d}$ is dense in $\R$, then there exists a sequence $\{k_j\}\subset\Z^d$ such that (\ref{seq_k-1}) is satisfied. In view of (\ref{defi_T_j}) and (\ref{defi_T_j-error}), we see that
$T_j$ is uniquely defined through $\{k_n\}_{1\leq n\leq j+1}$ and
$\la k_{j+1}\ra< 8T_j^{-1}$. Hence, the fact $t^{g(t)}\to\infty$ shown in Lemma \ref{lemma_gt} implies that
$$\la k_{j+1} \ra^{\frac{g(T_j)}{2}}<T_j^{-\frac{g(T_j)}{2}}\to 0,$$
which guarantees the existence of above sequences $\{k_j\}$ and $\{T_j\}$.

\

Now, for given $s> 0$, let us come back to Eq. (\ref{eq_original}):
$$
{\rm i}\partial_t \psi=(\CJ+{\CP}(t))\psi,\quad \psi(0)\in\CH^{s+2} \ {\rm non-vanishing}.
$$
For any $\varepsilon>0$, we choose $k_0, k_1\in\Z^d\setminus\{0\}$ with $\la k_0\ra, \la k_1\ra>0$ such that
 \begin{equation}\label{k0k1_small}
 |1-\la k_0+k_1\ra| e^{2 r |k_0|}+  \la k_1\ra<\frac{\varepsilon}{64},
\end{equation}
$g(T_0)\in (0,1)$ and $\la k_1\ra^{\frac18}(1+\Gamma_s(\psi(0)))$ sufficiently small (depending on $s$).
With $\{k_j\}$ and $\{T_j\}$ generated as in (\ref{defi_T_j}) -- (\ref{seq_k-2}),
let us define sequences $\{\varphi_j\}_{j\in\N^*}, \{\lambda_j\}_{j\in\N^*}\subset (0,1)$ and $\{z_j\}_{j\in\N^*}\subset (1,\infty)$ by
\begin{equation}\label{seq_vp_lm}
 \vp_j= \la k_{j+1}\ra^{\frac34 g(T_j)}, \quad \lm_j:=\sqrt{\vp_j^2-\la k_{j+1}\ra^2},\quad z_j=\sqrt{\frac{\vp_j+\lm_j}{\la k_{j+1}\ra}}.
\end{equation}
The above choice of $\{\vp_j\}$ and $\{\lm_j\}$ implies that
$$ \frac{\left(\frac{\varphi_j+\lambda_j}{\la k_{j+1} \ra}\right)^s}{f(T_j)}\sim_s \frac{\varphi_j^s\cdot  T_j^s}{T_j^{(1- g(T_j))s}}\sim_s \frac{T_j^{(1-\frac34 g(T_j))s}}{T_j^{(1- g(T_j))s}} =T_j^{\frac{g(T_j)s}{4}}\to \infty,\quad j\to\infty,$$
which provides the ingredient for the proof of (\ref{limsup-thm}) (see (\ref{eq_rmk_limsup}) in Remark \ref{rmk_limsup} and the proof of Proposition \ref{prop_limsup}).

\smallskip

We have the following concrete almost reducibility argument for the non-reducible Eq. (\ref{eq_original}).
\begin{prop}\label{prop_lower}
For any $r>0,$ $\varepsilon>0$, there exists $P\in C_r^\omega(\T^d,\sl)$  satisfying $\|P\|_{r}<\varepsilon$, such that, for $j\in\N^*$, Eq. (\ref{eq_original}) is conjugated to the equation
\begin{equation}\label{eq_conj_j-lower}
{\rm i}\partial_t\psi_j = ({\CL}_j+{\CP}_{j+1}(t )) \psi_j,
\end{equation}
via a time quasi-periodic unitary transformation $\psi(t)={\CU}_j(t)\psi_{j}(t)$, with
${\CU}_j(t)={\CM}(U_j(t))$ the Meta representation of
\begin{equation}\label{U_j_prop}
U_1(t)=R_{\la k_0+k_1\ra t}, \qquad
U_j(t)=R_{\la k_0+k_1\ra t} \prod_{n=1}^{j-1}\left(\left(\begin{array}{cc}
                                  z_n & 0 \\
                                  0  & z^{-1}_n
                                \end{array}
\right) R_{\la k_{n+1} \ra t}\right), \quad j\geq  2,
\end{equation}
where ${\CL}_j$ and ${\CP}_{j+1}(\cdot)$ are defined as in (\ref{L_j}) and (\ref{P_j+1}) with
\begin{equation}\label{LjPj+1}
L_j=\left(
\begin{array}{cc}
  0 & \varphi_j+\lambda_j \\
  -(\varphi_j-\lambda_j) & 0
\end{array}
\right),\quad \|P_{j+1}\|_{\T^d}<\la k_{j+2}\ra^{\frac{g(T_{j+1})}2}.
\end{equation}
Moreover, we have $\{\sup_{t\in\R} \|U_j(t)\|\}\to\infty$ as $j\to\infty$.
\end{prop}

\begin{remark}
With the fast-decaying positive sequence $\{\la k_{j}\ra\}$ chosen according to the prescribed $o(t^s)-$growth rate $f(t)$, the construction of $L_j$ in (\ref{LjPj+1}), as well as the corresponding classical propagator
$$
e^{tL_j}=\left(\begin{array}{cc}
              \cos(\la k_{j+1} \ra t) & (\varphi_j+\lambda_j)\frac{\sin(\la k_{j+1} \ra t)}{\la k_{j+1} \ra} \\
              -(\varphi_j-\lambda_j)\frac{\sin(\la k_{j+1} \ra t)}{\la k_{j+1} \ra} & \cos(\la k_{j+1} \ra t)
            \end{array} \right),
$$
incites the oscillatory growth of Sobolev norms for the quantized Hamiltonian via the computations with Meta representation. More precisely, the matrix element $(\varphi_j+\lambda_j)\frac{\sin(\la k_{j+1} \ra t)}{\la k_{j+1} \ra}$ of $e^{tL_j}$ satisfies that
\begin{equation}\label{low-up}
0\leq \left|(\varphi_j+\lambda_j)\frac{\sin(\la k_{j+1} \ra t)}{\la k_{j+1} \ra}\right|\leq  \frac{\varphi_j+\lambda_j}{\la k_{j+1}\ra},
\end{equation}
with the upper bound $\frac{\varphi_j+\lambda_j}{\la k_{j+1} \ra}$ large enough and tending to $\infty$ as $j\to \infty$.
Moreover, both bounds in (\ref{low-up}) can be reached since $\left|\sin(\la k_{j+1} \ra t)\right|$ oscillates slowly between $0$ and $1$. In this sense, it provides the clue to the oscillatory orbit.

It is crucial that $L_j\in {\rm sl}(2,\R)$ given in (\ref{LjPj+1}) is of the elliptic type, but not in the standard elliptic form. Otherwise, for the standard elliptic form $\la k_{j+1} \ra J$ and its propagator
$$ e^{t \la k_{j+1} \ra J}=\left(\begin{array}{cc}
              \cos(\la k_{j+1} \ra t) & \sin(\la k_{j+1} \ra t) \\
              -\sin(\la k_{j+1} \ra t) & \cos(\la k_{j+1} \ra t)
            \end{array} \right),$$
the uniform boundedness of matrix elements will not give rise to the growth of Sobolev norms through Meta representation.
\end{remark}

\subsection{Almost reducibility and Anosov-Katok construction for linear systems}\label{sec_al_ODE}

Since the quadratic quantum Hamiltonians in Eq. (\ref{eq_original}) and (\ref{eq_conj_j}) satisfy that
\begin{eqnarray*}
\CJ+{\CP}(t)&=&\frac{-1}2 \left\la \left(\begin{array}{c}
X\\
D
 \end{array}\right) , J (J+P(\om t)) \left(\begin{array}{c}
X\\
D
 \end{array}\right) \right\ra,\\
\CL_j+{\CP}_{j+1}(t)&=&\frac{-1}2 \left\la \left(\begin{array}{c}
X\\
D
 \end{array}\right) , J (L_j+P_{j+1}(\om t)) \left(\begin{array}{c}
X\\
D
 \end{array}\right) \right\ra,
\end{eqnarray*}
according to Proposition \ref{prop_elem_meta}, Proposition \ref{prop_upper} and \ref{prop_lower} can be shown through the Meta representaion with the almost reducibility of the linear system $(\om, J+P(\cdot))$:
\begin{equation}\label{diagram}
\begin{array}{ccccc}
(\om, J+P(\cdot)) &{\longrightarrow} &     {\rm i}\partial_t\psi=(\CJ+{\CP}(t))\psi  \\
    &  &        \\
  \  \  \   \  \    \big\downarrow  \,  U_j(t)  &    & \  \  \  \  \  \  \  \  \  \  \  \  \  \     \big\downarrow  \,  {\CU}_j(t)={\CM}(U_j(t))   \\
    &  &        \\
(\om, L_j+P_{j+1}(\cdot))& {\longrightarrow} &     {\rm i}\partial_t\psi_j=({\CL}_j+{\CP}_{j+1}(t))\psi_j
  \end{array}   \  \  \   .
\end{equation}

In the remaining of this section, let us focus on the quasi-periodic linear system $(\om, J+P(\cdot))$:
\begin{equation}\label{lin_sys_original}
\dot{y}(t)= (J + P(\om t)) y(t), \quad J=\left(
\begin{array}{cc}
  0 & 1 \\
  -1 & 0
\end{array}
\right) ,\quad P(\cdot)=\left(
\begin{array}{cc}
  p_{11}(\cdot) & p_{02}(\cdot) \\
  -p_{20}(\cdot) & -p_{11}(\cdot)
\end{array}
\right).
\end{equation}
It is well known from Theorem A of Eliasson \cite{Eli1992} that  $(\om, J+P(\cdot))$ is reducible if the rotation number ${\rm rot}(\omega, J+P)$ is Diophantine or rational w.r.t. $\om$. However, for the case that ${\rm rot}(\omega, J+P(\cdot))$ is Liouville w.r.t. $\omega$, the reducibility is sometimes not expected, since, in KAM scheme aiming for reducibility, it is possible to meet the resonances infinitely many times, and the corresponding renormalization, which is not close to identity, makes the convergence of sequence of transformation unrealizable.

Nevertheless, almost reducibility is always true for the system $(\om, J+P(\cdot))$. We state in the following the almost reducibility argument with the KAM scheme developed in \cite{LYZZ}, where all the close-to-identity transformations in two consecutive resonant steps are combined into one close-to-identity transformation and, as in (\ref{sequence_epsilon}), the size of time dependent perturbation decays much faster than that in \cite{Eli1992}.

\begin{prop}\label{prop_al_ODE}
There exists $\varepsilon_*=\varepsilon_*(\gamma,\tau,r,d)>0$ such that, if $\|P\|_{r}<\varepsilon_*$, then the following holds.
\begin{enumerate}
\item The linear system $(\omega, J+P(\cdot))$ is almost reducible, i.e., there exist sequences $\{L_j\}\subset \sl$ and $\{P_j\}\subset C^\om(\T^d, \sl)$ with $\|P_j\|_{\T^d} \to 0$, such that the system (\ref{lin_sys_original}) is conjugated to
$$\dot{y}_j(t)=(L_j+P_{j+1}(\om t))y_j(t)$$
via a time quasi-periodic change of variables $y=U_j(t)y_j$.
\item If $(\omega, J+P(\cdot))$ is reducible, then there exists $N\in\N^*$ such that, for $j\geq N$,
\begin{equation}\label{CV_lin_sys_Uj}
U_j(t)=U_N(t) \ with \ \sup_{t\in\R}\|U_N(t)\|<\infty, \quad L_j = L_N, \quad P_{j+1}=0.
\end{equation}
Otherwise, $\displaystyle \sup_{t\in\R} \|U_j(t)\|\to \infty$ as $j\to\infty$, and there exists $\{\varepsilon_j\}_{j\in\N}\subset (0,1)$ satisfying (\ref{sequence_epsilon}) such that
 (\ref{l_small_rho}) is satisfied.
\end{enumerate}
\end{prop}

Since almost reducibility of quasi-periodic $\sl$ linear system has been well studied in previous works \cite{Eli1992, LYZZ}, we give a sketch of proof in Appendix \ref{App_proof_ar}.

\smallskip

\noindent
{\bf Proof of Proposition \ref{prop_upper}.} The assertion (1), i.e., the almost reducibility of Eq. (\ref{eq_original}), follows immediately from that of the linear system (\ref{lin_sys_original}) through (\ref{diagram}).

If Eq. (\ref{eq_original}) is non-reducible, i.e., the sequence of unitary transformations $\{{\CU}_j(t)\}$ obtained in the assertion (1) does not converge to a $L^2-$unitary transformation uniformly in $t$, then the linear system $(\omega, J+P(\cdot))$ is not reducible neither. Indeed, if the system $(\omega, J+P(\cdot))$ is reducible, then, according to (\ref{CV_lin_sys_Uj}), $\{{\CU}_j(t)\}=\{\CM(U_j(t))\}$ converges to ${\CU}_N(t)=\CM(U_N(t))$. Then the assertion (2) is shown.\qed

\smallskip

Let us turn to the proof of Proposition \ref{prop_lower}, which is based on the fibred Anosov-Katok construction \cite{KXZ2021}. As we mentioned,  Anosov-Katok construction is in some sense the counterpart of the KAM method: KAM tends to prove rigidity in the Diophantine world, while Anosov-Katok is used in order to prove non-rigidity in the Liouvillean world.
In our context, the concept of almost reducibility, obtained via KAM, allows us to study the rigidity results if the rotation number ${\rm rot}(\omega, J+P)$ is Diophantine or rational w.r.t. $\om$, while the fibered Anosov-Katok construction will be an efficient method to study the non-rigidity results if ${\rm rot}(\omega, J+P(\cdot))$ is Liouville w.r.t. $\omega$. Different from classical Anosov-Katok construction that is only valid in $C^{\infty}$ topology \cite{FK2005}, our fibred construction works in the analytic topology.
The main reason is that our method is a combination of the almost reducibility scheme and Anosov-Katok construction.
Roughly speaking, the almost reducibility scheme developed in \cite{LYZZ} mainly consists of two ingredients:
a close-to-identity transformation to remove all non-resonant terms of the perturbation and a rotation to eliminate the truncated resonant terms (see also Appendix \ref{App_proof_ar}).
In our fibred Anosov-Katok construction, we assume, for simplification, that there is no non-resonant terms at all. Hence there is no close-to-identity transformation appearing in the almost reducibility scheme, and we only need to construct the resonant terms according to the prescribed resonant sites, to keep the fibered rotation number of the conjugated normal form rational. This scheme works clearly in the analytic topology.
However, the conjugation $h_{n}$ in the classical  Anosov-Katok construction (obtained in \eqref{hn}) is $q_{n}$-periodic, and then is extended $1-$periodic, thus in general is not analytic.

\smallskip

With sequences $\{k_j\}$, $\{T_j\}$, $\{\vp_j\}$, $\{\lm_j\}$, $\{z_j\}$ given as in Section \ref{sec_concrete_PDE}, we have the fibred Anosov-Katok construction  for the linear system $(\om, J + P(\cdot))$ in (\ref{lin_sys_original}).

\begin{prop}\label{prop_AK_ODE}
For any $r>0,$ $\varepsilon>0$,    there exists $P\in C_r^\omega(\T^d,\sl)$  satisfying $\|P\|_{r}<\varepsilon$,  such that, for every $j\in\N^*$, the quasi-periodic linear system $(\om, J + P(\cdot))$ is conjugated to the linear system
$\dot{y}_j(t)= (L_j + P_{j+1}(\om t)) y_j(t)
$,
via the time quasi-periodic transformation $y(t)=U_j(t) y_j(t)$ with $L_j$, $P_{j+1}$ and $U_j$ the same as in (\ref{U_j_prop}) and (\ref{LjPj+1}).
\end{prop}

Before entering the proof of Proposition \ref{prop_AK_ODE}, let us introduce several estimates (thanks to the fast decay property of $\{\la k_j\ra\}$) on sequences given in Section \ref{sec_concrete_PDE}.
These estimates will be also applied in Section \ref{sec_proof_osci} to describe the oscillatory growth of Sobolev norms.

\begin{lemma}\label{lemma_seq_kT} For $\la k_1\ra$ sufficiently small,
the sequences $\{\la k_j\ra\}$ and $\{T_j\}$ satisfy that, for $j\geq 1$,
\begin{eqnarray}
\la k_{j+1} \ra^{\frac{g(T_j)}{33^{j+1-n}}}&<& { e^{-2 r |k_n|} \la k_{n} \ra}, \quad 1\leq n \leq j \label{esti_chain_k}\\
\la k_{j+1}\ra^{\frac{g(T_j)}{32}}&<&\prod_{n=1}^{j}\la k_{n}\ra, \label{seqest_1}\\
\left|\la k_{j+1}\ra T_j -\frac{5\pi}{2}\right|   &<& \la k_{j+1}\ra^{1-\frac{g(T_j)}{32}}, \label{seqest_2}\\
{\rm dist}\left(\la k_{n}\ra T_j,2\pi\Z\right) &<& \la k_{n}\ra^{1-\frac{g(T_{n-1})}{32}}, \quad 2\leq n\leq j. \label{seqest_3}
\end{eqnarray}
\end{lemma}
\proof In view of (\ref{seq_k-2}), we have immediately that $\la k_{j+1} \ra^{\frac{g(T_j)}{33}}< { e^{-2 r |k_j|} \la k_{j} \ra} $. Hence, for $1\leq n\leq j$, we have (\ref{esti_chain_k}) by seeing that
$$
e^{-2 r |k_n|} \la k_{n}\ra> \la k_{n+1}\ra^{\frac{g(T_n)}{33}}>\cdots >\la k_{j+1}\ra^{\frac{g(T_j)\cdots g(T_{n+1})g(T_n)}{33^{j+1-n}}}> \la k_{j+1}\ra^{\frac{g(T_j)}{33^{j+1-n}}},
$$
which means that
$$
\prod_{n=1}^{j}\la k_{n}\ra > \la k_{j+1}\ra^{g(T_j)\sum_{n=1}^{j}\frac{1}{33^{j+1-n}}}> \la k_{j+1}\ra^{\frac{g(T_j)}{32}}.
$$
By (\ref{defi_k0k1}) -- (\ref{defi_T_j-error}), we have that
$$\left|\la k_{j+1} \ra T_j -\frac{5\pi}{2}\right|\leq 2\pi\left(\frac{\la k_{j+1}\ra}{\la k_0+k_1\ra}+\sum_{n=2}^j\frac{\la k_{j+1}\ra}{\la k_n\ra}\right)< 2\pi\la k_{j+1}\ra\left(2+\sum_{n=2}^j\la k_{j+1}\ra^{-\frac{g(T_j)}{33^{j+1-n}}}\right),$$
and for $2\leq n\leq j$,
$${\rm dist}\left(\la k_{n}\ra T_j,2\pi\Z\right) <2\pi\left(\frac{\la k_{n}\ra}{\la k_0+k_1\ra}+ \sum_{l=2}^{n-1}\frac{\la k_{n}\ra}{\la k_l\ra}\right)<2\pi\la k_{n}\ra \left(2+\sum_{l=1}^{n-1}\la k_{n}\ra^{-\frac{g(T_{n-1})}{33^{n-l}}}\right).$$
Since the sequence $\{\la k_{j+1}\ra^{g(T_j)}\}$ decays super-exponentially w.r.t. $j$,
the above estimates imply (\ref{seqest_2}) and (\ref{seqest_3}).\qed

\begin{lemma}\label{lemma_seq_z}
The sequence $\{z_j\}$ satisfies that
\begin{eqnarray}
\frac{1}{\sqrt{2}}\la k_{j+1} \ra^{\frac38 g(T_{j})-\frac12} < z_j < \sqrt{2}\la k_{j+1} \ra^{\frac38 g(T_{j})-\frac12},\label{seq_z0}\\
|z_n\sin(\la k_{n+1} \ra T_j)|< z_n^{-1} \la k_{n+1} \ra^{\frac{5}8 g(T_n)}, \quad 1\leq n\leq j-1,\label{seq_z1}\\
\prod_{n=1}^{j-1} z_n < \la k_{j+1} \ra^{-\frac{g(T_{j})}{2200}}, \quad \la k_{j+1} \ra^{\frac{g(T_{j})}{1100}} \prod_{n=1}^{j-1} z_n < \prod_{n=1}^{j-1} z^{-1}_n.\label{seq_z2}
\end{eqnarray}
\end{lemma}
\proof According to (\ref{seq_vp_lm}), we have
$$z_j=\sqrt{\frac{\vp_j+\lm_j}{\la k_{j+1}\ra}}\in \left(\frac{1}{\sqrt{2}}\la k_{j+1} \ra^{\frac38 g(T_{j})-\frac12}, \sqrt{2}\la k_{j+1} \ra^{\frac38 g(T_{j})-\frac12}\right).$$
In view of (\ref{seqest_3}), we have, for $1\leq n\leq j-1$,
\begin{eqnarray*}
|z_n\sin(\la k_{n+1} \ra T_j)|&<& \sqrt{2}\la k_{n+1} \ra^{\frac38g(T_n)-\frac12}\cdot
 \la k_{n+1}\ra^{1-\frac{g(T_{n})}{32}} \\
   &<& \frac{1}{\sqrt{2}}\la k_{n+1} \ra^{\frac12+\frac{g(T_n)}4} \\
   &<& z_n^{-1} \la k_{n+1} \ra^{\frac{5}8 g(T_n)}.
\end{eqnarray*}
According to (\ref{esti_chain_k}), we obtain (\ref{seq_z2}) by
$$ \prod_{n=1}^{j-1} z_n< \sqrt{2}^{j-1}\prod_{n=1}^{j-1} \la k_{n+1} \ra^{-\frac12+\frac38g(T_n)}
<\sqrt{2}^{j-1}  \la k_{j+1} \ra^{-\frac{g(T_{j})}2\sum_{n=1}^{j-1}\frac{1}{33^{j+1-n}}}
< \la k_{j+1} \ra^{-\frac{g(T_{j})}{2200}},$$
$$\la k_{j+1} \ra^{\frac{g(T_{j})}{1100}} \prod_{n=1}^{j-1} z_n < \la k_{j+1} \ra^{\frac{g(T_{j})}{2200}} < \prod_{n=1}^{j-1} z^{-1}_n.\qed $$

\noindent
{\bf Proof of Proposition \ref{prop_AK_ODE}.} Let us construct the time dependent part $P(\cdot)\in C^\om(\T^d,\sl)$ such that Proposition \ref{prop_AK_ODE} holds for the linear system $(\om, J + P(\cdot))$.
Recall the given $k_0,k_1\in \Z^d$ and the sequences $\{k_j\}$, $\{T_j\}$ defined by (\ref{seq_k-1}) and (\ref{seq_k-2}) with (\ref{k0k1_small}) satisfied, as well as $\{\vp_j\}$, $\{\lm_j\}$, $\{z_j\}$ defined in (\ref{seq_vp_lm}). We also define the diagonal matrix
\begin{equation}\label{diagonal_Zj}
Z_j:= \left(\begin{array}{cc}
  z_j & 0 \\
  0 & z_j^{-1}
\end{array}
\right).
\end{equation}
By (\ref{seq_z0}) in Lemma \ref{lemma_seq_z}, we have
\begin{equation}\label{esti_Zj}
\|Z_j\|, \  \|Z_j^{-1}\|\leq z_j +z_j^{-1} \leq  2\la k_{j+1}\ra^{\frac38 g(T_{j})-\frac12}.
\end{equation}
The construction is composed of the following three steps.

\smallskip

\noindent{\bf Step 1.} Preliminary change of variables $y(t)= R_{\la k_0+k_1 \ra t} y_1(t)$.

\smallskip

Let us define $P(\cdot)\in C^\om(\T^d, \sl)$ as
\begin{equation}\label{p-pre}
P(\cdot)=R_{\la k_0,\cdot \ra}\left(
\begin{array}{cc}
  0 & -1+\la k_0+k_1\ra \\[2mm]
  1-\la k_0+ k_1\ra & 0
\end{array}
\right)R^{-1}_{\la k_0,\cdot \ra} + \tilde P_1(\cdot):= F_0(\cdot) + \tilde P_1(\cdot),
\end{equation}
with the error term $\tilde P_1(\cdot)$, to be determined later, satisfying
 $\|\tilde P_1\|_{r}<\la k_2 \ra^{\frac{g(T_1)}{4}}$.
By the change of variables
\begin{equation}\label{pre_change-0}
y(t)= R_{\la k_0 \ra t} y_0(t),
\end{equation}
we obtain the quasi-periodic linear system
$$\dot{y}_0=\left(-R^{-1}_{\la k_0 \ra t}\dot{R}_{\la k_0 \ra t}+R^{-1}_{\la k_0 \ra t}\left( J + P(\om t)\right)R_{\la k_0 \ra t}\right)y_0
= \left(\la k_1\ra J
+P_1(\om t) \right) y_0$$
with the time quasi-periodic perturbation
$P_1(\cdot):= R^{-1}_{\la k_0,\cdot \ra} \tilde P_1(\cdot)R_{\la k_0,\cdot \ra} $.

Now let us define $P_1(\cdot)$ as
\begin{equation}\label{p-F1}
P_1(\cdot)=\left(
\begin{array}{cc}
  \lambda_1 \sin(2\la k_1, \cdot \ra) & \varphi_1+\lambda_1 \cos(2\la k_1, \cdot \ra) \\[2mm]
  -\varphi_1+\lambda_1 \cos(2\la k_1, \cdot \ra) & -\lambda_1 \sin(2\la k_1, \cdot \ra)
\end{array}
\right) + \tilde P_2(\cdot):= F_1(\cdot) + \tilde P_2(\cdot),
\end{equation}
with the coefficients in the leading term $F_1$ of perturbation defined as in (\ref{seq_vp_lm}):
$$\varphi_1= \la k_2 \ra^{\frac34 g(T_1)}, \quad \lm_1= \sqrt{\varphi_1^2-\la k_2 \ra^2},$$
and the error term $\tilde P_2(\cdot)$, to be determined later, satisfying
 $\|\tilde P_2\|_{r}<\la k_3 \ra^{\frac{g(T_2)}{4}}$.
By the change of variables
\begin{equation}\label{pre_change-1}
y_0(t)= R_{\la k_1 \ra t} y_1(t),
\end{equation}
we obtain the first quasi-periodic linear system
$$\dot{y}_1=\left(-R^{-1}_{\la k_1 \ra t}\dot{R}_{\la k_1 \ra t}+R^{-1}_{\la k_1 \ra t}\left(\la k_1 \ra J + P_1(\om t)\right)R_{\la k_1 \ra t}\right)y_1
= (L_1 + P_2(\om t)) y_1 $$
with the matrices of coefficients
$$L_1=\left(
\begin{array}{cc}
  0 & \varphi_1+\lambda_1 \\
  -(\varphi_1-\lambda_1) & 0
\end{array}
\right),\quad P_2(\cdot)= R^{-1}_{\la k_1,\cdot \ra} \tilde P_2(\cdot)R_{\la k_1,\cdot \ra} \ .$$

\smallskip

\noindent{\bf Step 2.} Iterative change of variables $y_j(t)=Z_j R_{\la k_{j+1} \ra t} y_{j+1}(t)$.

\smallskip

Suppose that we arrive at the $j^{\rm th}$ quasi-periodic linear system, $j\geq 1$,
$$\dot{y}_j(t)= (L_j + P_{j+1}(\om t)) y_j(t),$$
where the matrices of coefficients are
$$L_j=\left(
\begin{array}{cc}
  0 & \varphi_j+\lambda_j \\
  -(\varphi_j-\lambda_j) & 0
\end{array}
\right),\quad P_{j+1}(\cdot)=  F_{j+1}(\cdot) + \tilde P_{j+2}(\cdot),$$
with the leading term of perturbation
\begin{equation}\label{F_j+1}
F_{j+1}(\cdot)=Z_j\left(
\begin{array}{cc}
  \lambda_{j+1} \sin(2\la k_{j+1}, \cdot \ra) & \varphi_{j+1}+\lambda_{j+1}\cos(2\la k_{j+1}, \cdot \ra) \\[2mm]
  -\varphi_{j+1}+\lambda_{j+1} \cos(2\la k_{j+1}, \cdot \ra) & -\lambda_{j+1}\sin(2\la k_{j+1}, \cdot \ra)
\end{array}
\right) Z_j^{-1},
\end{equation}
and the error term $\tilde P_{j+2}(\cdot)$, to be determined in the next steps, satisfying
\begin{equation}\label{small_tildeP}
\|\tilde P_{j+2}\|_{r}< \la k_{j+3} \ra^{\frac{g(T_{j+2})}{4}}.
\end{equation}
By the change of variables
\begin{equation}\label{iter_change}
y_{j}(t)=Z_j R_{\la k_{j+1} \ra t} y_{j+1}(t) \  \  {\rm with}    \  \  Z_j= \left(
\begin{array}{cc}
  z_j & 0 \\
  0 & z_j^{-1}
\end{array}
\right),\quad  z_j=\sqrt{\frac{\varphi_j+\lambda_j}{\la k_{j+1} \ra}},\end{equation}
we obtain, noting that $Z_j$ normalize $L_j$,
\begin{eqnarray*}
\dot{y}_{j+1}&=&\left(-R^{-1}_{\la k_{j+1} \ra t} \dot{R}_{\la k_{j+1} \ra t}+R^{-1}_{\la k_{j+1} \ra t}Z_j^{-1}\left(L_j + P_{j+1}(\om t)\right)Z_j R_{\la k_{j+1} \ra t}\right)y_{j+1}\\
&=& (L_{j+1} + P_{j+2}(\om t)) y_{j+1}
\end{eqnarray*}
with the new matrices of coefficients
$$L_{j+1}=\left(
\begin{array}{cc}
  0 & \varphi_{j+1}+\lambda_{j+1} \\
  -(\varphi_{j+1}-\lambda_{j+1}) & 0
\end{array}
\right),\quad P_{j+2}(\cdot)= R^{-1}_{\la k_{j+1} ,\cdot \ra}Z_j^{-1} \tilde P_{j+2}(\cdot)Z_j R_{\la k_{j+1},\cdot \ra }.$$

\smallskip

\noindent{\bf Step 3.} Estimates on changes of variables and perturbations.

\smallskip

The above construction shows the conjugation between two linear systems $(\om, L_j+P_{j+1}(\cdot))$ and $(\om, L_{j+1}+P_{j+2}(\cdot))$.  For $F_{j+1}$, $j\geq 0$, given in (\ref{F_j+1}), according to (\ref{seq_vp_lm}) and Lemma \ref{lemma_seq_kT}, \ref{lemma_seq_z}, we see that $F_{j+1}\in C_{r}^{\omega}(\T^d, \sl)$, with
\begin{eqnarray}
\|F_{j+1}\|_{r} &\leq& 2 \|Z_j\| \left(\varphi_{j+1}+ 2 \lambda_{j+1}e^{2 r |k_{j+1}|}\right) \|Z_j^{-1}\|  \label{esti_Fj+1}\\
  &\leq& 8\la k_{j+1}\ra^{\frac34 g(T_{j})-1} \left(\la k_{j+2}\ra^{\frac34g(T_{j+1})}+2 e^{2 r |k_{j+1}|} \la k_{j+2}\ra^{\frac34g(T_{j+1})}  \right)\nonumber\\
  &\leq& 24 \la k_{j+1} \ra^{-1} \cdot e^{2 r |k_{j+1}|}\la k_{j+2}\ra^{\frac34g(T_{j+1})}  \nonumber\\
  &\leq&\la k_{j+2}\ra^{\frac58g(T_{j+1})}.\nonumber
\end{eqnarray}
Combining with (\ref{small_tildeP}), we obtain
$$\|P_{j+1}\|_{r}\leq \|F_{j+1}\|_{r} + \|\tilde P_{j+2}\|_{r}\leq \la k_{j+2} \ra^{\frac58g(T_{j+1})} + \la k_{j+3} \ra^{\frac{g(T_{j+2})}{4}}\leq
\la k_{j+2} \ra^{\frac{g(T_{j+1})}{2}} .$$

With $\{F_{j}(\cdot)\}$ given explicitly in (\ref{p-pre}), (\ref{p-F1}) and (\ref{F_j+1}), let
\begin{eqnarray}
P(\cdot)&:=& F_0(\cdot)+\tilde P_1 (\cdot )\label{perturbation_AnosovKatok}\\
&=& F_0(\cdot)+ R_{\la k_0,\cdot \ra} F_1(\cdot)R^{-1}_{\la k_0,\cdot \ra} +R_{\la k_0+ k_1,\cdot \ra} F_2(\cdot)R^{-1}_{\la k_0+ k_1,\cdot \ra}  \nonumber\\
   & & + \, \sum_{j\geq 2}R_{\la k_0+ k_1,\cdot \ra} \left(\prod_{n=1}^{j-1} Z_n  R_{\la k_{n+1},\cdot \ra} \right)  F_{j+1}(\cdot) \left(\prod_{n=j-1}^{1} R^{-1}_{\la k_{n+1},\cdot \ra} Z_n^{-1} \right) R^{-1}_{\la k_0+ k_1,\cdot \ra} \ . \nonumber
\end{eqnarray}
 Through (\ref{esti_Fj+1}), we have that
\begin{align*}
\|R_{\la k_0,\cdot \ra} F_1(\cdot)R^{-1}_{\la k_0,\cdot \ra}\|_{r}\leq & \,  16e^{2r |k_0|}\la k_{2}\ra^{\frac58g(T_{1})}< \la k_{2}\ra^{\frac{g(T_{1})}2}, \\
\|R_{\la k_0+ k_1,\cdot \ra} F_2(\cdot)R^{-1}_{\la k_0+ k_1,\cdot \ra}\|_{r}\leq & \,  16e^{2r |k_0+k_1|}\la k_{3}\ra^{\frac58g(T_{2})}< \la k_{3}\ra^{\frac{g(T_{2})}2}.
\end{align*}
Since (\ref{seq_z0}) in Lemma \ref{lemma_seq_z} implies that
\begin{equation}\label{ZZ}
\|Z_n\|, \  \|Z_n^{-1}\|\leq 2\la k_{n+1}\ra^{\frac38 g(T_{n})-\frac12}, \quad 1\leq n \leq j-1,
\end{equation}
 the general term of the sum in (\ref{perturbation_AnosovKatok}) satisfies that
\begin{eqnarray}
& & \left\|R_{\la k_0+ k_1,\cdot \ra} \left(\prod_{n=1}^{j-1} Z_n  R_{\la k_{n+1},\cdot \ra} \right)  F_{j+1}(\cdot) \left(\prod_{n=j-1}^{1}  R^{-1}_{\la k_{n+1},\cdot \ra} Z_n^{-1} \right) R^{-1}_{\la k_0+k_1,\cdot \ra} \right\|_{r} \label{RZRF}\\
&<& 4^{3j-1} e^{2|k_0+k_1|r} \prod_{n=1}^{j-1} e^{2|k_{n+1}|r} \la k_{j+2}\ra^{\frac58g(T_{j+1})} \prod_{n=1}^{j-1}z_n^2 \nonumber\\
&<& 4^{3j-1} e^{2|k_0+k_1|r} \prod_{n=1}^{j-1} e^{2|k_{n+1}|r} \la k_{j+2}  \ra^{\frac58g(T_{j+1})} \la k_{j+1} \ra^{-\frac{g(T_{j})}{1100}} \nonumber\\
&<&\la k_{j+2}\ra^{\frac{g(T_{j+1})}2}.\nonumber
\end{eqnarray}
Noting that (\ref{k0k1_small}) implies $\|F_0\|_{r}<\frac{\varepsilon}2$, we see that the fast decay property of $\left\{\la k_{j+1}\ra^{g(T_{j})}\right\}$ guarantees the convergence of $P$ given in (\ref{perturbation_AnosovKatok}), with
$$\|P\|_{r}\leq \|F_0\|_{r} +  \|\tilde P_1\|_{r} \leq  \frac12 \varepsilon+ \sum_{j\geq 1}\la k_{j+1}\ra^{\frac{g(T_j)}2}<  \frac12 \varepsilon+ \la k_{2}\ra^{\frac{g(T_1)}4}< \varepsilon.$$

Moreover, we have that the error terms $\tilde P_2$ and $\tilde P_{j+2}$ in the preliminary and iterative steps are indeed
\begin{eqnarray*}
\tilde P_2(\cdot)&= &  R_{\la k_1,\cdot \ra} F_2(\cdot)R^{-1}_{\la k_1,\cdot \ra}\\
& &  + \, \sum_{j\geq 2}R_{\la k_1,\cdot \ra} \left(\prod_{n=1}^{j-1} Z_n  R_{\la k_{n+1},\cdot \ra} \right)  F_{j+1}(\cdot) \left(\prod_{n=j-1}^{1} R^{-1}_{\la k_{n+1},\cdot \ra} Z_n^{-1} \right) R^{-1}_{\la k_1,\cdot \ra},
\end{eqnarray*}
and for $j\geq 1$,
$$\tilde P_{j+2}(\cdot)= \sum_{l\geq j+1} \left(\prod_{n=j}^{l-1} Z_n  R_{\la k_{n+1},\cdot \ra} \right)  F_{l+1}(\cdot) \left(\prod_{n=l-1}^{j} R^{-1}_{\la k_{n+1},\cdot \ra} Z_n^{-1} \right).$$
With (\ref{ZZ}) and (\ref{RZRF}), we have
 $$\|\tilde P_2\|_{r}<\la k_3 \ra^{\frac{g(T_2)}{4}} ,\qquad \|\tilde P_{j+2}\|_{r}<\la k_{j+3} \ra^{\frac{g(T_{j+2})}{4}}, \quad j\geq 1, $$
which verifies the hypothesis about the estimates on error terms.

For the sequence of transformations
$$U_1(t):=R_{\la k_0+k_1 \ra t}  , \quad  U_{j+1}(t):= R_{\la k_0+k_1 \ra t} \prod_{n=1}^{j} \left( Z_n  R_{\la k_{n+1} \ra t}\right),$$
we see that $\{\sup_t \|U_j(t)\|\}\to\infty$,
since
$$U_{j+1}(0)= \prod_{n=1}^{j}Z_n =\left(\begin{array}{cc}
\prod_{n=1}^{j}z_n & 0 \\
0 & \prod_{n=1}^{j}z_n^{-1}
\end{array} \right)$$
and the element of the diagonal matrix $Z_n$ satisfies
$$z_n> \frac1{\sqrt{2}}\la k_{n+1} \ra^{\frac38 g(T_n)-\frac12}\to \infty ,\quad n\to \infty.$$
Proposition \ref{prop_AK_ODE} is shown.\qed

\

\section{Oscillatory growth of Sobolev norm -- Proof of Theorem \ref{thm_growth_upper} and \ref{thm_sub_polynome}}\label{sec_pr_thm}

\subsection{Global well-posedness and growth of Sobolev norms}
Consider the quantum Hamiltonian system
${\rm i}\partial_t\psi = \CL(t) \psi$,
where the linear operator $\CL(t)=\CL(t,X,D)$ is a self-adjoint quadric form of $(X,D)$ with time dependent coefficients, i.e.,
$$ \CL(t)=  \frac{-1}2 \left\la \left(\begin{array}{c}
X\\
D
 \end{array}\right) , J L(t) \left(\begin{array}{c}
X\\
D
 \end{array}\right) \right\ra,$$
 with $L(t)\in \sl$ smoothly depending on $t$. We also assume that there is $\beta>0$ such that
$$\sup_t\|L(t)\|\leq\beta.$$
Given $s\geq 0$, let us define the constant
$$\Upsilon_s:= \max \left\{\begin{array}{cc}
                             \|[X^2, \CJ^{\frac{s}{2}}]\CJ^{-\frac{s}{2}}\|_{{\CB}({\CH}^0)}, & \|[XD, \CJ^{\frac{s}{2}}]\CJ^{-\frac{s}{2}}\|_{{\CB}({\CH}^0)} , \\[1mm]
                             \|[DX, \CJ^{\frac{s}{2}}]\CJ^{-\frac{s}{2}}\|_{{\CB}({\CH}^0)} , & \|[D^2, \CJ^{\frac{s}{2}}]\CJ^{-\frac{s}{2}}\|_{{\CB}({\CH}^0)}
                           \end{array}\right\}.$$
We will use in the sequel the following simplified version of Theorem 1.2 of \cite{MR2017}
about the global well-posedness in Sobolev space and the exponential upper bound of Sobolev norm.

\begin{lemma}\label{lemma_wellposed} [Maspero-Robert \cite{MR2017}] Given $s\geq 0$ and $\psi(0)\in {\CH}^s$, we have $\psi(t)\in {\CH}^s$ for every $t\in\R^+$, with
$$\|\psi(t)\|_{{\CH}^s}\leq 2 e^{4\beta \Upsilon_s t}\|\psi(0)\|_{{\CH}^s}.$$
\end{lemma}

\medskip

\subsection{The $o(t^s)-$upper bound of Sobolev norms.}

With the almost reducibility scheme stated in Proposition \ref{prop_upper}, we are going to show in this subsection the upper bound (\ref{upper-thm}) of  $\CH^s-$norm for the solution to Eq. (\ref{eq_original}).

\begin{prop}\label{prop_growth_upper} Consider Eq. (\ref{eq_original}) with $\|P\|_{r}<\varepsilon_*$ as in Proposition \ref{prop_upper} and with
\begin{equation}\label{div_Uj}
\sup_{t\in\R}\{\|U_j(t)\|\}\to \infty,  \quad  j\to\infty.
\end{equation}
Then, for $j$ sufficiently large (depending on $s$),  we have
$$\|\psi(t)\|_{{\CH}^s} \lesssim_s \left\|\psi(0)\right\|_{{\CH}^{s+2}} |\ln\varepsilon_j|^{(4s+4)\tau}\varepsilon_{j}^{\frac{s}{16}}t^s, \quad  \varepsilon_{j-1}^{-\frac{1}{20\tau}-\frac{1}{16}}< t\leq \varepsilon_j^{-\frac{1}{20\tau}-\frac{1}{16}}.$$
\end{prop}

\begin{remark}
Proposition \ref{prop_growth_upper} implies the $o(t^s)-$upper bound of $\|\psi(t) \|_{{\CH}^s}$ under the divergence assumption of $\sup_t\{\|U_j(t)\|\}$. According to Proposition \ref{prop_upper}, the non-reducibility assumption of Eq. (\ref{eq_original}) in Theorem \ref{thm_growth_upper} implies (\ref{div_Uj}), i.e., the non-reducibility of the linear system (\ref{lin_sys_original}). Then $o(t^s)-$upper bound (\ref{upper-thm}) in Theorem \ref{thm_growth_upper} is shown.
In particular, for Eq. (\ref{eq_original}) with the perturbation given in Proposition \ref{prop_lower}, we have the upper bound (\ref{upper-thm}) in Theorem \ref{thm_growth_upper} when $\varepsilon<\varepsilon_*$.
\end{remark}

\proof
For fixed $j\in\N^*$, let us firstly consider the solution to the equation
\begin{equation}\label{eq_j-upper}
{\rm i}\partial_t \psi_j(t,x)=({\CL}_j+{\CP}_{j+1}(t))\psi_j(t,x),\quad \psi_j(0)=\psi_j(0,\cdot)\in {\CH}^{s+2},
\end{equation}
where the linear operator ${\CL}_j$ and ${\CP}_{j+1}(t)$ are defined as in (\ref{L_j}) and (\ref{P_j+1}) in  Proposition \ref{prop_upper}. By the fact that $\|P_{j+1}\|_{\T^d}<\varepsilon_{j+1}$, we see that
\begin{equation}\label{error_s+2-s}
\sup_t\|{\CP}_{j+1}( t )\|_{{\CB}({\CH}^{s+2}\to {\CH}^{s})} \lesssim_s  \varepsilon_{j+1}.
\end{equation}
By Duhamel formula, we have that the solution to Eq. (\ref{eq_j-upper}) has the explicit form
\begin{equation}\label{Duhamel}
\psi_j(t)=e^{-{\rm i}t{\CL}_j}\psi_j(0)-{\rm i}  \int_0^t e^{{\rm i}(t'-t){\CL}_j} {\CP}_{j+1}( t') \psi_j(t')  \, dt',
\end{equation}
where $e^{-{\rm i}t{\CL}_j}$ is the Meta representation of $e^{tL_j}$, i.e., for $L_j=\left(\begin{array}{cc}
             l_j^{11}  & l_j^{02}\\[1mm]
             -l_j^{20} & -l_j^{11}
             \end{array}
\right)$,
$$e^{-{\rm i}t{\CL}_j}={\CM}(e^{tL_j})={\CM}\left(\begin{array}{cc}
             \cos(\varrho_j t) + l_j^{11}\varrho_j^{-1}\sin(\varrho_j t)  & l_j^{02}\varrho_j^{-1}\sin(\varrho_j t) \\[1mm]
             -l_j^{20}\varrho_j^{-1}\sin(\varrho_j t)  & \cos(\varrho_j t)-l_j^{11}\varrho_j^{-1}\sin(\varrho_j t)
             \end{array}
\right)$$
with $\varrho_j:=\sqrt{{\rm det}(L_j)}$. In view of (\ref{l_small_rho}), we have that $\varrho_j>\gamma\cdot |\ln\varepsilon_{j+1}|^{-\tau}$.
Then, according to Proposition \ref{cor_metaplectic_upper} and (\ref{l_small_rho}), we see that
\begin{equation}\label{linear_flow_part}
\|e^{\pm{\rm i}t{\CL}_j}\|_{{\CB}({\CH^s})} \lesssim_s 1+\left(\left|l_j^{20}\right|+\left|l_j^{11}\right|+\left|l_j^{02}\right|\right)^s\left|\frac{\sin(t\varrho_j)}{\varrho_j}\right|^s \lesssim_s  1+\varepsilon_{j}^{\frac{s}{16}} \left|\frac{\sin(t\varrho_j)}{\varrho_j}\right|^s.
\end{equation}
With $\psi_j(0)\in {\CH}^{s+2}$, applying Lemma \ref{lemma_wellposed} to Eq. (\ref{eq_j-upper}) with $\beta=3\varepsilon_j^{\frac{1}{16}}$, we see that
$$\|\psi_j(t')\|_{{\CH}^{s+2}}\leq 2\|\psi_j(0)\|_{{\CH}^{s+2}}\exp\left\{12\varepsilon_j^{\frac{1}{16}} \Upsilon_{s+2} t'\right\},$$
and hence, through (\ref{sequence_epsilon}), for $j$ large enough such that
$$\varepsilon_{j+1}^{\frac14}\varepsilon_j^{-\frac{1}{20\tau}-\frac{1}{16}} |\ln\varepsilon_{j+1}|^{s\tau}<1,\quad \varepsilon_{j}^{\frac1{180\tau}}<\frac{r}{96\Upsilon_{s+2} },$$
we have, by (\ref{error_s+2-s}) and (\ref{linear_flow_part}), for $\varepsilon_{j-1}^{-\frac{1}{20\tau}-\frac{1}{16}}< t\leq \varepsilon_j^{-\frac{1}{20\tau}-\frac{1}{16}}$,
\begin{eqnarray*}
& &\left\|{\rm i}  \int_0^t \left(e^{{\rm i}(t'-t){\CL}_j} {\CP}_{j+1}( t')\right) \psi_j(t')  \, dt'\right\|_{{\CH}^s}\\
&\lesssim_s& \varepsilon_{j+1}\|\psi_j(0)\|_{{\CH}^{s+2}} \, \varepsilon_j^{-\frac{1}{20\tau}-\frac{1}{16}} \left( 1+ \frac{\varepsilon_{j}^{\frac{s}{16}}|\ln\varepsilon_{j+1}|^{s\tau}}{\gamma^{s}}\right)\exp\left\{12\varepsilon_j^{\frac{1}{16}} \Upsilon_{s+2} t\right\}\\
   &\lesssim_s& \varepsilon^{\frac34}_j \|\psi_j(0)\|_{{\CH}^{s+2}} \exp\left\{-\frac38r\varepsilon_j^{-\frac{1}{18\tau}}+12\Upsilon_{s+2}\varepsilon_j^{-\frac{1}{20\tau}}\right\}\\
   &\lesssim_s&\varepsilon^{\frac34}_j \|\psi_j(0)\|_{{\CH}^{s+2}}.
\end{eqnarray*}
Since the ${\CH}^s-$norm of the term $e^{-{\rm i}t{\CL}_j}\psi_j(0)$ can be estimated through (\ref{linear_flow_part}), we have
\begin{eqnarray*}
\|\psi_j(t)\|_{{\CH}^s}&\leq& \|e^{-{\rm i}t{\CL}_j}\|_{{\CB}({\CH^s})} \|\psi_j(0)\|_{{\CH}^s}+ \left\|{\rm i} \int_0^t \left(e^{{\rm i}t'{\CL}_j} {\CP}_{j+1}( t'\right) \psi_j(t')  \, dt'\right\|_{{\CH}^s}\\
&\lesssim_s& \left(1+\varepsilon_{j}^{\frac{s}{16}} \left|\frac{\sin(t\varrho_j)}{\varrho_j}\right|^s\right) \left\|\psi_j(0)\right\|_{{\CH}^{s}} +\varepsilon^{\frac34}_j \|\psi_j(0)\|_{{\CH}^{s+2}}\\
&\lesssim_s& \left\|\psi_j(0)\right\|_{{\CH}^{s+2}}\left(1+ \varepsilon_{j}^{\frac{s}{16}} t^s\right),\quad \varepsilon_{j-1}^{-\frac{1}{20\tau}-\frac{1}{16}}< t\leq \varepsilon_j^{-\frac{1}{20\tau}-\frac{1}{16}}.
\end{eqnarray*}

Let us come back to Eq. (\ref{eq_original}).
Under the divergence assumption (\ref{div_Uj}) of $\{U_j(t)\}$, we have,
according to Proposition \ref{prop_upper},
$$ \psi(t)={\CU}_j(t) \psi_j(t), \quad  {\CU}_j(t) = {\CM}(U_j(t)) \  \  {\rm with}  \  \  \sup_t\|U_j(t)\|\leq  |\ln\varepsilon_j|^{2\tau} $$
such that Eq. (\ref{eq_original}) is conjugated to Eq. (\ref{eq_j-upper}).
By Proposition \ref{prop_elem_meta} and Lemma \ref{lemma_esti_metaplectic}, we have
$$\sup_t\|{\CU}_j(t)\|_{{\CB}({\CH}^s)}\lesssim_s |\ln\varepsilon_j|^{2\tau s}, \quad \sup_t\|{\CU}_j(t)^{-1}\|_{{\CB}({\CH}^{s+2})}\lesssim_s |\ln\varepsilon_j|^{2\tau(s+2)}.$$
Combining with Proposition \ref{prop_growth_upper}, we obtain, for $\varepsilon_{j-1}^{-\frac{1}{20\tau}-\frac{1}{16}}< t\leq \varepsilon_j^{-\frac{1}{20\tau}-\frac{1}{16}}$,
\begin{eqnarray*}
\|\psi(t) \|_{{\CH}^s} &\leq&  \|{\CU}_j(t)\|_{{\CB}({\CH}^s)} \cdot \|\psi_j(t) \|_{{\CH}^s}  \\
&\lesssim_s&|\ln\varepsilon_j|^{2\tau s}\left(1+ \varepsilon_{j}^{\frac{s}{16}} t^s\right)  \cdot \|{\CU}_j(0)^{-1}\|_{{\CB}({\CH}^{s+2})}\left\|\psi(0)\right\|_{{\CH}^{s+2}}\\
&\lesssim_s& \left\|\psi(0)\right\|_{{\CH}^{s+2}} |\ln\varepsilon_j|^{(4s+4)\tau}\varepsilon_{j}^{\frac{s}{16}}t^s.
\end{eqnarray*}
Proposition \ref{prop_growth_upper} is shown.\qed

\subsection{Oscillatory growth of Sobolev norms.}\label{sec_proof_osci}

With the concrete almost reducibility of Eq. (\ref{eq_original}) stated in Proposition \ref{prop_lower}, we are going to show (\ref{limsup-thm}) and (\ref{liminf-thm}) in Theorem \ref{thm_sub_polynome} by constructing suitable sequences of moments $\{T_j\}$ and $\{t_j\}$ as mentioned in (\ref{rem_seq_moments}).
This implies the oscillatory growth of Sobolev norms $\|\psi(t)\|_{\CH^s}$, as well as the optimality of the $o(t^s)-$upper bound obtained in Proposition \ref{prop_growth_upper} when $\varepsilon<\varepsilon_*$.

Fix $s>0$.
For $j\in\N^*$, Eq. (\ref{eq_original}) with $\psi(0)\in {\CH}^{s+2}$ is conjugated to the equation
\begin{equation}\label{eq_j-lower}
{\rm i}\partial_t \psi_j(t,x)=({\CL}_j+{\CP}_{j+1}(t))\psi_j(t,x),\quad \psi_j(0)=\psi_j(0,\cdot)\in {\CH}^{s+2},
\end{equation}
via the time quasi-periodic $L^2-$unitary transformation $\psi(t)={\CU}_j(t)\psi_j(t)$,
where the linear operators ${\CL}_j$, ${\CP}_{j+1}(t)$ and ${\CU}_j(t)$ are defined as in Proposition \ref{prop_lower}, with
\begin{equation}\label{small_P_j+1}
\sup_{t}\|{\CP}_{j+1}( t )\|_{{\CB}({\CH}^{s+2}\to {\CH}^{s})} \lesssim_s  \la k_{j+2}\ra^{\frac{g(T_{j+1})}2}.
\end{equation}
According to Duhamel formula, the solution to Eq. (\ref{eq_j-lower}) is still of the form (\ref{Duhamel}):
$$
\psi_j(t)=e^{-{\rm i}t{\CL}_j}\psi_j(0)-{\rm i} \int_0^t e^{{\rm i}(t'-t){\CL}_j} {\CP}_{j+1}( t') \psi_j(t')  \, dt'.$$
Then, through the conjugation $\psi(t)={\CU}_j(t)\psi_j(t)$, we have
$$\psi(t)={\CU}_j(t)e^{-{\rm i}t{\CL}_j} {\CU}_j(0)^{-1} \psi(0) - {\rm i} \,  {\CU}_j(t)  \int_0^t e^{{\rm i}(t'-t){\CL}_j} {\CP}_{j+1}( t')\psi_j(t')  \, dt'.$$

Recall that $e^{-{\rm i}t {\CL}_j}$ is the Meta representation of $e^{tL_j}$ with
\begin{eqnarray}
 L_j&=&\left(\begin{array}{cc}
              0 & \vp_j+\lm_j \\
              -(\vp_j-\lm_j) & 0
            \end{array} \right), \label{L_j_matrix}\\
e^{tL_j}&=& \left(\begin{array}{cc}
              \cos(\la k_{j+1} \ra t) & (\vp_j+\lm_j)\frac{\sin(\la k_{j+1} \ra t)}{\la k_{j+1} \ra} \\
              -(\vp_j-\lm_j)\frac{\sin(\la k_{j+1} \ra t)}{\la k_{j+1} \ra} & \cos(\la k_{j+1} \ra t)
            \end{array} \right).\label{etL_j}
\end{eqnarray}
Then, according to Proposition \ref{cor_metaplectic_upper}, we have, for any $t\in\R$,
\begin{equation}\label{up_eitLj}
\|e^{\pm{\rm i}t {\CL}_j}\|_{{\CB}({\CH^s})} \lesssim_s 1+\left((\varphi_j+\lm_j)+(\varphi_j-\lm_j)\right)^s \left|\frac{\sin(\la k_{j+1} \ra t)}{\la k_{j+1} \ra}\right|^s
\lesssim_s \la k_{j+1}\ra^{(\frac34g(T_j)-1)s}.
\end{equation}
On the other hand, in view of Lemma \ref{lemma_seq_kT} and \ref{lemma_seq_z}, $U_j(t)$ defined in (\ref{U_j_prop}) satisfies that
$$\sup_t\|U_j(t)\|\leq 4^{j} \prod_{n=1}^{j-1} z_n < 4^j  \la k_{j+1} \ra^{-\frac{g(T_{j})}{2200}} <\la k_{j+1}\ra^{-\frac{g(T_{j})}{2000}},$$
which implies, through Lemma \ref{lemma_esti_metaplectic},
\begin{equation}\label{up_Ujt}
\sup_{t}\|{\CU}_j(t)\|_{{\CB}({\CH}^s)}\lesssim_s \la k_{j+1}\ra^{-\frac{g(T_{j})}{2000}s} ,  \quad \sup_{t}\|{\CU}_j(t)^{-1}\|_{{\CB}({\CH}^{s+2})} \lesssim_s \la k_{j+1}\ra^{-\frac{g(T_{j})}{2000}(s+2)}.
\end{equation}
Recalling that $|\vp_j\pm\lm_j|<2\la k_{j+1}\ra^{\frac34g(T_j)}$, we have,
in view of Lemma \ref{lemma_wellposed},
\begin{eqnarray*}
\|\psi_j(t')\|_{{\CH}^{s+2}}&<& 2 e^{12\la k_{j+1}\ra^{\frac34g(T_j)} \Upsilon_{s+2} t'}\|\psi_j(0)\|_{{\CH}^{s+2}}\\
&\lesssim_s&  e^{12\la k_{j+1}\ra^{\frac34g(T_j)} \Upsilon_{s+2} t'} \la k_{j+1}\ra^{-\frac{g(T_{j})}{2000}(s+2)} \|\psi(0)\|_{{\CH}^{s+2}}.
\end{eqnarray*}
With (\ref{small_P_j+1}), (\ref{up_eitLj}) and (\ref{up_Ujt}), we have, for $j$ large enough (depending on $s$ and $\|\psi(0)\|_{{\CH}^{s+2}}$), for $0< t\leq 4 T_j$,
\begin{eqnarray}
& & \left\|{\CU}_j(t)  \int_0^t e^{{\rm i}(t'-t){\CL}_j} {\CP}_{j+1}( t') \, \psi_j(t')  \, dt'\right\|_{{\CH}^{s}}\label{esti_error-psi}\\
&\lesssim_s& \|\psi(0)\|_{{\CH}^{s+2}} \la k_{j+2}\ra^{\frac{g(T_{j+1})}2} \cdot T_{j}\la k_{j+1}\ra^{\left(\left(\frac34-\frac{1}{1000}\right)g(T_j)-1\right)s-\frac{g(T_j)}{1000}} \exp\left\{48 \Upsilon_{s+2} T_j\la k_{j+1}\ra^{\frac34g(T_j)} \right\} \nonumber\\
&\lesssim_s&\la k_{j+2}\ra^{\frac{g(T_{j+1})}4},\nonumber
\end{eqnarray}
since, by the construction of $\{k_j\}$ in (\ref{seq_k-2}), we have,
\begin{eqnarray*}
& &\la k_{j+2}\ra^{\frac{g(T_{j+1})}{8}}\cdot T_{j}\la k_{j+1}\ra^{\left(\left(\frac34-\frac{1}{1000}\right)g(T_j)-1\right)s-\frac{g(T_j)}{1000}}\\
&<& 9\pi \la k_{j+1}\ra^{\frac{g(T_{j})}{8}}\la k_{j+1}\ra^{\left(\left(\frac34-\frac{1}{1000}\right)g(T_j)-1\right)s-\frac{g(T_j)}{1000}-1}\exp\left\{-\frac14\la  k_{j+1} \ra^{-\left(1+\frac{1}{36\tau}\right)}\right\}\\
&<& \la  k_{j+1} \ra^{-s-1} \exp\left\{-\frac14\la  k_{j+1} \ra^{-\left(1+\frac{1}{36\tau}\right)}\right\} \\
&<&1,
\end{eqnarray*}
and furthermore,
\begin{eqnarray*}
& &\la k_{j+2} \ra^{\frac{g(T_{j+1})}{8}}  \|\psi(0)\|_{{\CH}^{s+2}}  \exp\left\{48 \Upsilon_{s+2} T_j\la k_{j+1}\ra^{\frac34g(T_j)} \right\} \\
&<& \|\psi(0)\|_{{\CH}^{s+2}}  \la k_{j+1}\ra^{\frac{g(T_{j})}{8}}\exp\left\{-\frac14\la k_{j+1}\ra^{-\left(1+\frac{1}{36\tau}\right)}+4000 \Upsilon_{s+2} \la k_{j+1}\ra^{\frac34g(T_j)-1}  \right\} \\
&<&1.
\end{eqnarray*}

The estimate (\ref{esti_error-psi}) means that the second part of $\psi(t)$ with the integral is negligible for $0< t\leq 4 T_j$.
Therefore, to estimate the ${\CH}^s-$norm of the solution $\psi(t)$ to Eq. (\ref{eq_original}) for $0< t\leq 4 T_j$, it is sufficient to focus on the part ${\CU}_j(t)e^{-{\rm i}t{\CL}_j} {\CU}_j^{-1}(0) \psi(0)$ generated by the linear flow.
In view of Proposition \ref{prop_elem_meta}, we see that the linear propagator ${\CU}_j(t)e^{-{\rm i}t{\CL}_j} {\CU}_j^{-1}(0)$ is indeed the Meta representation:
\begin{equation}\label{meta-UeU}
{\CU}_j(t)e^{-{\rm i}t{\CL}_j} {\CU}_j^{-1}(0) = {\CM}\left(U_j(t)\right) {\CM}\left(e^{tL_j}\right) {\CM}\left(U_j^{-1}(0)\right)
= {\CM}\left(U_j(t) e^{tL_j} U_j^{-1}(0)\right),
\end{equation}
where, as in (\ref{U_j_prop}) and (\ref{diagonal_Zj}),
\begin{equation}\label{U_j}
U_j(t)=R_{\la k_0+k_1\ra t}\prod_{n=1}^{j-1} (Z_n R_{\la k_{n+1}\ra t}), \quad  Z_n:= \left(\begin{array}{cc}
                                       z_n & 0 \\
                                       0 & z_n^{-1}
                                     \end{array}\right).
\end{equation}
The oscillatory growth of Sobolev norm will be deduced from the following two propositions.

\begin{prop}\label{prop_limsup}
For non-vanishing $\psi(0)\in {\CH}^s$, we have, for $j\geq 2$,
$$\|{\CU}_j(T_j)e^{-{\rm i}T_j{\CL}_j} {\CU}_j^{-1}(0) \psi(0)\|_{{\CH}^s}\gtrsim_s T_j^{(1-\frac34 g(T_j))s} \|D^s \psi(0)\|_{L^2} .$$
\end{prop}
\begin{remark}\label{rmk_limsup}
Combining (\ref{esti_error-psi}) and Proposition \ref{prop_limsup}, we obtain (\ref{limsup-thm}) since
\begin{equation}\label{eq_rmk_limsup}
\frac{\|\psi(T_j) \|_{{\CH}^s} }{f(T_j)} \gtrsim_s  \frac{T_j^{(1-\frac34 g(T_j))s}\|D^s \psi(0)\|_{L^2}}{T_j^{(1- g(T_j))s}} =T_j^{\frac{g(T_j)s}{4}}\|D^s \psi(0)\|_{L^2} \to \infty, \quad j\to \infty.
\end{equation}
\end{remark}
\proof In view of (\ref{U_j}), we see that
$$U_j(0)=\prod_{n=1}^{j-1} Z_n=\left(\begin{array}{cc}
                       \prod_{n=1}^{j-1}z_n  & 0 \\
                       0 & \prod_{n=1}^{j-1}z_n^{-1}
                     \end{array}
\right)=:\left(\begin{array}{cc}
               {\CZ}_j  & 0 \\
               0 & {\CZ}_j^{-1}
              \end{array}
\right).$$
Then, by applying Lemma \ref{lemma_esti_metaplectic} and \ref{lemma_metaplectic_lower-1} to ${\CU}_j^{-1}(0)={\CM}(U_j^{-1}(0))$, we have
$$\|{\CU}_j^{-1}(0) \psi(0)\|_{{\CH}^s}\lesssim_s {\CZ}_j^{s}\|\psi(0)\|_{{\CH}^s},\quad \|D^s{\CU}_j^{-1}(0) \psi(0)\|_{L^2}= {\CZ}_j^{s} \|D^s  \psi(0)\|_{L^2}.$$
Hence, recalling $\Gamma_s(\cdot)$ defined in Section \ref{sec_Nota},
\begin{equation}\label{Delta_var}
\Gamma_s\left({\CU}_j^{-1}(0) \psi(0)\right)=\frac{\|{\CU}_j^{-1}(0) \psi(0)\|_{{\CH}^s}}{\|D^s{\CU}_j^{-1}(0) \psi(0)\|_{L^2}} \lesssim_s \frac{{\CZ}_j^{s} \|\psi(0)\|_{{\CH}^s}}{{\CZ}_j^{s} \|D^s\psi(0)\|_{L^2}} = \Gamma_s\left(\psi(0)\right).
\end{equation}
In view of (\ref{seqest_2}) in Lemma \ref{lemma_seq_kT}, we have $\sin\left(\la k_{j+1}\ra T_j \right)>\frac34$, then
the fast decay property of the sequence $\{\la k_j\ra\}$ and (\ref{Delta_var}) implies that
$$\frac{(\varphi_j+\lm_j)\sin\left(\la k_{j+1}\ra T_j \right)}{\la k_{j+1}\ra(1+\Gamma_s({\CU}_j^{-1}(0) \psi(0)))^2}\gtrsim_s  \frac{\la k_{j+1}\ra^{\frac34 g(T_j)-1}}{(1+\Gamma_s(\psi(0)))^2} \gtrsim_s \la k_1 \ra^{-\frac14}(1+\Gamma_s(\psi(0)))^{-2}.$$
Since $\la k_1\ra^{\frac18}(1+ \Gamma_s(\psi(0)))$ is sufficiently small, we can apply Proposition \ref{cor_metaplectic_lower} to $e^{T_jL_j}$ and its Meta representation $e^{-{\rm i}T_j{\CL}_j}$, and obtain
\begin{eqnarray}
\|e^{-{\rm i}T_j{\CL}_j}{\CU}_j^{-1}(0) \psi(0)\|_{{\CH}^s}&\gtrsim_{s}& (\varphi_j+\lambda_j)^s \frac{\left|\sin\left(\la k_{j+1} \ra T_j\right)\right|^s}{\la k_{j+1}\ra^s} \left\|D^s{\CU}_j^{-1}(0) \psi(0)\right\|_{L^2} \label{esti_lower_1} \\
&\gtrsim_{s}& {\CZ}^s_j \la k_{j+1} \ra^{(\frac34 g(T_j)-1)s} \|D^s  \psi(0)\|_{L^2}.\nonumber
\end{eqnarray}
By (\ref{seqest_3}) in Lemma \ref{lemma_seq_kT}, we see that $R_{\la k_{n+1}\ra T_j}$, $1\leq n \leq j-1$, is diagonally dominant, with
$$\left|\sin(\la k_{n+1}\ra T_j) (\cos(\la k_{n+1}\ra T_j))^{-1}\right|= |\tan(\la k_{n+1}\ra T_j)|< \frac32\la k_{n+1}\ra^{1-\frac{g(T_{n})}{32}}.$$
Then, according to Lemma \ref{lemma_metaplectic_lower-1}, we see that, for any $u\in{\CH}^s$,
$$ \|{\CM}(R_{\la k_{n+1}\ra T_j})u\|_{{\CH}^s}\gtrsim_{s} \cos^{s}(\la k_{n+1}\ra T_j) \|u\|_{{\CH}^s}. $$
Since, for the diagonal $Z_n$, Lemma \ref{lemma_metaplectic_lower-1} implies that
$\|{\CM}(Z_n)u\|_{{\CH}^s} \gtrsim_{s} z_n^{-s} \|u\|_{{\CH}^s}$, we have
$$\|{\CM}(Z_n R_{\la k_{n+1}\ra T_j})u\|_{{\CH}^s}=\|{\CM}(Z_n){\CM}(R_{\la k_{n+1}\ra T_j})u\|_{{\CH}^s}\gtrsim_{s} z_n^{-s}\cos^s(\la k_{n+1}\ra T_j) \|u\|_{{\CH}^s}.$$
Noting that $R_{\la k_0+k_{1}\ra T_j}=I$ (recalling (\ref{defi_T_j-error})), we have, by Proposition \ref{prop_elem_meta},
$${\CU}_j(T_j)={\CM}(U_j(T_j))={\CM}(Z_1 R_{\la k_{2}\ra T_j})\cdots {\CM}(Z_{j-1} R_{\la k_{j}\ra T_j}),\quad j\geq 2,$$
and hence,  for any $u\in{\CH}^s$,
\begin{equation}\label{esti_lower_2}
\|{\CU}_j(T_j)u\|_{{\CH}^s}   \gtrsim_{s} \left(\prod_{n=1}^{j-1} z_n^{-1}\cos(\la k_{n+1}\ra T_j)\right)^{s} \|u\|_{{\CH}^s}  \gtrsim_{s} {\CZ}_j^{-s} \|u\|_{{\CH}^s},
\end{equation}
where, by (\ref{seqest_3}) in Lemma \ref{lemma_seq_kT}, we have
\begin{eqnarray*}
\prod_{n=1}^{j-1} \cos(\la k_{n+1}\ra T_j)&\geq& \prod_{n=1}^{j-1} \left(1-\dist(\la k_{n+1}\ra T_j,2\pi\Z)^2\right)\\&\geq& \prod_{n=1}^{j-1} \exp\left\{\ln\left(1-\la k_{n+1}\ra^{2-\frac{g(T_{n})}{16}}\right)\right\}\\
&\geq&  \exp\left\{-\frac54\sum_{n=1}^{j-1}\la k_{n+1}\ra^{2-\frac{g(T_{n})}{16}} \right\}\\
&>&\frac{9}{10}.
\end{eqnarray*}
Combining (\ref{esti_lower_1}) and (\ref{esti_lower_2}), we have, for $j\geq 2$,
$$\|{\CU}_j(T_j)e^{-{\rm i}T_j{\CL}_j} {\CU}_j^{-1}(0) \psi(0)\|_{{\CH}^s}\gtrsim_s \la k_{j+1} \ra^{(\frac34 g(T_j)-1)s} \|D^s  \psi(0)\|_{L^2}
  \gtrsim_s  T_j^{(1-\frac34 g(T_j))s}\|D^s  \psi(0)\|_{L^2} .\qed$$

\medskip

\begin{prop}\label{prop_liminf}
For $j\geq 2$,
$\|{\CU}_j(4 T_j)e^{-4{\rm i} T_j{\CL}_j} {\CU}_j^{-1}(0)\|_{{\CB}({\CH}^s)}\lesssim_s 1$.
\end{prop}

\begin{remark}
Combining (\ref{esti_error-psi}) and Proposition \ref{prop_liminf}, we obtain that $\|\psi(4T_j) \|_{{\CH}^s}\lesssim_s 1 $,
which implies (\ref{liminf-thm}) via an infimum limit with the sequence of moments $\{t_j\}=\{4T_j\}$.
\end{remark}

\proof In view of (\ref{meta-UeU}) and Lemma \ref{lemma_esti_metaplectic}, it is sufficient to give a constant upper bound for the $\SL-$matrix $U_j(4 T_j) e^{4T_jL_j} U_j^{-1}(0)$.

\begin{lemma}\label{lemma_est_U_4T} For $U_j(t)=\left(\begin{array}{cc}
          U_j(t)_{11}&  U_j(t)_{12}\\
         U_j(t)_{21}&  U_j(t)_{22}
        \end{array}
\right)$, $j\geq 2$, we have
\begin{equation}\label{est_U_4T}
\left|U_j(4T_j)_{11}\right| , \left|U_j(4T_j)_{21}\right| < \frac54{\CZ}_j , \quad
\left|U_j(4T_j)_{12}\right|, \left|U_j(4T_j)_{22}\right|  <  \frac54 {\CZ}_j^{-1}.
\end{equation}
\end{lemma}
\proof Let $V_n(t):=Z_n R_{\la k_{n+1}\ra t}$, $1\leq n \leq j-1$.
We are going to show that, for $t=4T_j$,
\begin{eqnarray}
\left|\left(\prod_{n=1}^{j-1} V_{n}(4T_j)\right)_{11}\right|, \left|\left(\prod_{n=1}^{j-1} V_{n}(4T_j)\right)_{21}\right|  &<& \left(1+\sum_{n=0}^{j-1}\la k_{n+1}\ra^{\frac{g(T_{n})}{16}}\right)\prod_{n=1}^{j-1} z_{n},\label{prod_V-1} \\
\left|\left(\prod_{n=1}^{j-1} V_{n}(4T_j)\right)_{12}\right| , \left|\left(\prod_{n=1}^{j-1} V_{n}(4T_j)\right)_{22}\right|  &<& \left(1+\sum_{n=0}^{j-1}\la k_{n+1}\ra^{\frac{g(T_n)}{16}}\right) \prod_{n=1}^{j-1} z^{-1}_{n},\label{prod_V-2}
\end{eqnarray}
from which we can deduce (\ref{est_U_4T}) since $\la k_0+ k_1 \ra T_j\in 2\pi\Z$.

 At first, (\ref{prod_V-1}) and (\ref{prod_V-2}) are satisfied for $j=2$ since
\begin{eqnarray*}
 V_1(t)=Z_1 R_{\la k_2\ra t}=\left(\begin{array}{cc}
           z_1\cos(\la k_2\ra t) & z_1\sin(\la k_2\ra t)\\
         -z_1^{-1}\sin(\la k_2\ra t) & z_1^{-1}\cos(\la k_2\ra t)
        \end{array}
\right),
\end{eqnarray*}
and (\ref{seq_z1}) in Lemma \ref{lemma_seq_z} implies that  $|z_1\sin(4\la k_2\ra T_j)|<4z_1^{-1}\la k_{2} \ra^{\frac{5}8 g(T_1)}$.
Assume that $j>2$, and for some $1\leq l \leq j-2$,  we have
\begin{eqnarray}
\left|\left(\prod_{n=1}^{l} V_{n}(4T_j)\right)_{11}\right|, \left|\left(\prod_{n=1}^{l} V_{n}(4T_j)\right)_{21}\right|&<& \left(1+\sum_{n=0}^{l}\la k_{n+1}\ra^{\frac{g(T_{n})}{16}}\right)\prod_{n=1}^{l} z_{n},\label{prod_V_l-1}\\
\left|\left(\prod_{n=1}^{l} V_{n}(4T_j)\right)_{12}\right| , \left|\left(\prod_{n=1}^{l} V_{n}(4T_j)\right)_{22}\right|  &<& \left(1+\sum_{n=0}^{l}\la k_{n+1}\ra^{\frac{g(T_n)}{16}}\right) \prod_{n=1}^{l} z_{n}^{-1} .\label{prod_V_l-2}
\end{eqnarray}
Let us verify the estimates for
\begin{eqnarray*}
\prod_{n=1}^{l+1} V_{n}(t) &=& \prod_{n=1}^{l} V_{n}(t)\left(\begin{array}{cc}
z_{l+1}\cos(\la k_{l+2}\ra t) & z_{l+1}\sin(\la k_{l+2}\ra t) \\[3mm]
-z_{l+1}^{-1}\sin(\la k_{l+2}\ra t) & z_{l+1}^{-1}\cos(\la k_{l+2}\ra t)
\end{array}\right),
\end{eqnarray*}
with the matrix elements given by
\begin{eqnarray*}
\left(\prod_{n=1}^{l+1} V_{n}(t)\right)_{11}&=&  \left(\prod_{n=1}^{l} V_{n}(t)\right)_{11} z_{l+1}\cos(\la k_{l+2}\ra t)- \left(\prod_{n=1}^{l}V_{n}(t)\right)_{12}z^{-1}_{l+1}\sin(\la k_{l+2}\ra t) ,\\
\left(\prod_{n=1}^{l+1} V_{n}(t)\right)_{12}&=& \left(\prod_{n=1}^{l} V_{n}(t)\right)_{11} z_{l+1}\sin(\la k_{l+2}\ra t) + \left(\prod_{n=1}^{l} V_{n}(t)\right)_{12} z^{-1}_{l+1}\cos(\la k_{l+2}\ra t),  \\
\left(\prod_{n=1}^{l+1} V_{n}(t)\right)_{21}&=& \left(\prod_{n=1}^{l} V_{n}(t)\right)_{21} z_{l+1}\cos(\la k_{l+2}\ra t)- \left(\prod_{n=1}^{l} V_{n}(t)\right)_{22} z^{-1}_{l+1}\sin(\la k_{l+2}\ra t), \\
\left(\prod_{n=1}^{l+1} V_{n}(t)\right)_{22}&=& \left(\prod_{n=1}^{l} V_{n}(t)\right)_{21} z_{l+1}\sin(\la k_{l+2}\ra t)+ \left(\prod_{n=1}^{l} V_{n}(t)\right)_{22} z^{-1}_{l+1}\cos(\la k_{l+2}\ra t).
\end{eqnarray*}
Under the assumptions (\ref{prod_V_l-1}) and (\ref{prod_V_l-2}), together with Lemma \ref{lemma_seq_kT} and \ref{lemma_seq_z}, we have, for $t=4T_j$,
\begin{eqnarray*}
 \left|\left(\prod_{n=1}^{l+1} V_{n}(t)\right)_{11}\right|
&\leq& \left|\left(\prod_{n=1}^{l} V_{n}(t)\right)_{11}\right| \left|z_{l+1}\cos(\la k_{l+2}\ra t)\right|
+\left|\left(\prod_{n=1}^{l} V_{n}(t)\right)_{12}\right|\left| z_{l+1}^{-1}\sin(\la k_{l+2}\ra t)\right|\\
&<&\left(1+\sum_{n=0}^{l}\la k_{n+1}\ra^{\frac{g(T_{n})}{16}}\right)\left(\prod_{n=1}^{l+1} z_{n}+4\la k_{l+2} \ra^{1-\frac{g(T_{l+1})}{32}}\prod_{n=1}^{l+1} z^{-1}_{n} \right)  \\
&<&\left(1+\sum_{n=0}^{l+1}\la k_{n+1}\ra^{\frac{g(T_{n})}{16}}\right) \prod_{n=1}^{l+1} z_{n},
\end{eqnarray*}
\begin{eqnarray*}
\left|\left(\prod_{n=1}^{l+1} V_{n}(t)\right)_{12}\right|
&\leq& \left|\left(\prod_{n=1}^{l} V_{n}(t)\right)_{11}\right| \left|z_{l+1}\sin(\la k_{l+2}\ra t)\right| + \left|\left(\prod_{n=1}^{l} V_{n}(t)\right)_{12}\right| \left|z^{-1}_{l+1}\cos(\la k_{l+2}\ra t)\right|\\
&<&\left(1+\sum_{n=0}^{l}\la k_{n+1}\ra^{\frac{g(T_n)}{16}}\right) \left(4z_{l+1}^{-1} \la k_{l+2} \ra^{\frac{5}8 g(T_{l+1})} \prod_{n=1}^{l} z_{n}  + \prod_{n=1}^{l+1} z^{-1}_{n}    \right)\\
&<&\left(1+\sum_{n=0}^{l}\la k_{n+1}\ra^{\frac{g(T_n)}{16}}\right) \left( \frac{\la k_{l+2} \ra^{\frac{g(T_{l+1})}{16}}}2 + 1  \right) \prod_{n=1}^{l+1} z^{-1}_{n}\\
&<&\left(1+\sum_{n=0}^{l+1}\la k_{n+1}\ra^{\frac{g(T_n)}{16}}\right) \prod_{n=1}^{l+1} z^{-1}_{n},
\end{eqnarray*}
and similarly,
\begin{eqnarray*}
\left|\left(\prod_{n=1}^{l+1} V_{n}(t)\right)_{21}\right|
&\leq& \left|\left(\prod_{n=1}^{l} V_{n}(t)\right)_{21}\right| \left|z_{l+1}\cos(\la k_{l+2}\ra t)\right| + \left|\left(\prod_{n=1}^{l} V_{n}(t)\right)_{22}\right| \left|z^{-1}_{l+1}\sin(\la k_{l+2}\ra t)\right|\\
&<&\left(1+\sum_{n=0}^{l+1}\la k_{n+1}\ra^{\frac{g(T_n)}{16}}\right) \prod_{n=1}^{l+1} z_{n},
\end{eqnarray*}
\begin{eqnarray*}
\left|\left(\prod_{n=1}^{l+1} V_{n}(t)\right)_{22}\right| &\leq & \left|\left(\prod_{n=1}^{l} V_{n}(t)\right)_{21} \right| \left|z_{l+1}\sin(\la k_{l+2}\ra t)\right|+ \left|\left(\prod_{n=1}^{l} V_{n}(t)\right)_{22}\right| \left|z^{-1}_{l+1}\cos(\la k_{l+2}\ra t) \right|\\
&< &\left(1+\sum_{n=0}^{l+1}\la k_{n+1}\ra^{\frac{g(T_n)}{16}} \right)\prod_{n=1}^{l+1} z^{-1}_{n}.
\end{eqnarray*}
Hence, (\ref{prod_V-1}) and (\ref{prod_V-2}) are shown.\qed

Now we are ready to prove Proposition \ref{prop_liminf}. Recall that
$$e^{ 4T_j L_j}= \left(\begin{array}{cc}
\cos(4T_j\la k_{j+1}\ra)& (\varphi_j+\lambda_j)\frac{\sin(4T_j\la k_{j+1}\ra)}{\la k_{j+1}\ra} \\[2mm]
-(\varphi_j-\lambda_j)\frac{\sin(4T_j\la k_{j+1}\ra )}{\la k_{j+1}\ra} & \cos(4T_j\la k_{j+1}\ra)
\end{array}\right). $$
By a direct calculation with $U_j^{-1}(0)=\left(\begin{array}{cc}
                       {\CZ}^{-1}_j  & 0 \\
                       0 & {\CZ}_j
                     \end{array}
\right)$, we have
$$ U_j(t) e^{ t L_j}U_j^{-1}(0) = \left(\begin{array}{cc}
 \begin{array}{c}
 {\CZ}^{-1}_j U_j(t)_{11} \cos(\la k_{j+1}\ra t) \\
 - {\CZ}^{-1}_j U_j(t)_{12}\cdot(\varphi_j-\lambda_j)\frac{\sin(\la k_{j+1}\ra t)}{\la k_{j+1}\ra}
 \end{array} &  \begin{array}{c}
{\CZ}_jU_j(t)_{11}\cdot(\varphi_j+\lambda_j)\frac{\sin(\la k_{j+1}\ra  t)}{\la k_{j+1}\ra}   \\
+ {\CZ}_j U_j(t)_{12}\cos(\la k_{j+1}\ra t)
 \end{array}   \\[6mm]
\begin{array}{c}
 {\CZ}^{-1}_j  U_j(t)_{21} \cos(\la k_{j+1}\ra t)  \\
 -  {\CZ}^{-1}_j U_j(t)_{22}\cdot(\varphi_j-\lambda_j)\frac{\sin(\la k_{j+1}\ra t)}{\la k_{j+1}\ra}
 \end{array} &  \begin{array}{c}
  {\CZ}_j   U_j(t)_{21}\cdot(\varphi_j+\lambda_j)\frac{\sin(\la k_{j+1}\ra  t)}{\la k_{j+1}\ra} \\
+ {\CZ}_j U_j(t)_{22}\cos(\la k_{j+1}\ra t)
 \end{array} \end{array}
\right).$$
Then, in view of Lemma \ref{lemma_est_U_4T}, we see that, for $t=4T_j$,
\begin{eqnarray*}
|{\CZ}^{-1}_j U_j(t)_{11}|, |{\CZ}^{-1}_j U_j(t)_{12}| , |{\CZ}^{-1}_j U_j(t)_{21}| , |{\CZ}^{-1}_j U_j(t)_{22}|, |{\CZ}_j U_j(t)_{12}| , |{\CZ}_j U_j(t)_{22}|&\leq&  \frac54,\\
|{\CZ}_j U_j(t)_{11}| , |{\CZ}_j U_j(t)_{21}|   &\leq & \frac54{\CZ}_j^2.
\end{eqnarray*}
Moreover, the matrix elements of $e^{4T_j L_j}$ satisfy that
$$|\varphi_j\pm\lambda_j|\frac{|\sin(4T_j\la k_{j+1}\ra)|}{\la k_{j+1}\ra}< \frac{2\la k_{j+1}\ra^{\frac34g(T_j)}}{\la k_{j+1}\ra} \cdot 4\la k_{j+1}\ra^{1-\frac{g(T_j)}{32}} < 8 \la k_{j+1}\ra^{(\frac34-\frac1{32}) g(T_j)}.$$
Recalling (\ref{seq_z1}) and (\ref{seq_z2}) in Lemma \ref{lemma_seq_z}, we have ${\CZ}_j=\prod_{n=1}^{j-1}z_n< \la k_{j+1} \ra^{-\frac{g(T_{j})}{2200}}$,
$${\CZ}_j^{2}|\varphi_j\pm\lambda_j|\frac{|\sin(4T_j\la k_{j+1}\ra)|}{\la k_{j+1}\ra}<8 \la k_{j+1} \ra^{-\frac{g(T_{j})}{1100}} \la k_{j+1}\ra^{(\frac34-\frac1{32}) g(T_j)}<\la k_{j+1}\ra^{\frac{g(T_j)}2}.$$
Then we see that $\|U_j(4 T_j) e^{4T_jL_j} U_j^{-1}(0) \|< 8$. Proposition \ref{prop_liminf} is shown.\qed

\appendix

\section{Proof of Proposition \ref{prop_elem_meta}}\label{app_pr_prop}

The assertions (i) -- (iii) are obvious by direct computations through (\ref{Meta_a}) and (\ref{Meta_b}).

\medskip
\noindent
Proof of (iv). By direct computations, we see that, for $A\in\SL$,
\begin{itemize}
  \item for $\kappa\in\R$,
  \begin{eqnarray*}
  \CM(A) \, \CM\left(\begin{array}{cc}
               1 & 0 \\
               \kappa & 1
             \end{array}\right) &=& \pm \,  \CM\left(A \, \left(\begin{array}{cc}
               1 & 0 \\
               \kappa & 1
             \end{array}\right)\right),  \\
    \CM\left(\begin{array}{cc}
               1 & 0 \\
               \kappa & 1
             \end{array}\right) \, \CM(A)&=&\pm  \,   \CM\left(\left(\begin{array}{cc}
               1 & 0 \\
               \kappa & 1
             \end{array}\right) \, A\right),
  \end{eqnarray*}
  \item for $\lambda\in\R\setminus\{0\}$,
  \begin{eqnarray*}
  \CM(A) \, \CM\left(\begin{array}{cc}
               \lambda & 0 \\
               0 & \lambda^{-1}
             \end{array}\right) &=&\pm  \,  \CM\left(A \, \left(\begin{array}{cc}
               \lambda & 0 \\
               0 & \lambda^{-1}
             \end{array}\right)\right), \\
     \CM\left(\begin{array}{cc}
               \lambda & 0 \\
               0 & \lambda^{-1}
             \end{array}\right) \, \CM(A) &=& \pm \,  \CM\left(\left(\begin{array}{cc}
               \lambda & 0 \\
               0 & \lambda^{-1}
             \end{array}\right) \, A\right),
  \end{eqnarray*}
\item $\CM(A)\CM(J)=\pm  \,  \CM(AJ)$, $\CM(A)\CM(-J)=\pm  \,  \CM(-AJ)$.
\end{itemize}
Then the assertion (iv) is shown, since for any $A=\left(\begin{array}{cc}
               a & b \\
               c & d
             \end{array}\right)\in\SL$, we have
\begin{itemize}
\item if $a\neq 0$,
\begin{equation}\label{decomp_a}
A=\left(\begin{array}{cc}
               1 & 0 \\
               ca^{-1} & 1
             \end{array}\right)\left(\begin{array}{cc}
               a & 0 \\
               0 & a^{-1}
             \end{array}\right)\left(\begin{array}{cc}
0 &  1 \\
-1 & 0
\end{array}\right)
\left(\begin{array}{cc}
1 &  0 \\
-ba^{-1} & 1
\end{array}\right)
\left(\begin{array}{cc}
0 &  -1 \\
1 & 0
\end{array}\right),
\end{equation}
\item if $b\neq 0$,
\begin{equation}\label{decomp_b}
A=\left(\begin{array}{cc}
               1 & 0 \\
               db^{-1} & 1
             \end{array}\right)\left(\begin{array}{cc}
               b & 0 \\
               0 & b^{-1}
             \end{array}\right)\left(\begin{array}{cc}
0 &  1 \\
-1 & 0
\end{array}\right)
\left(\begin{array}{cc}
1 &  0 \\
ab^{-1} & 1
\end{array}\right). \qed
\end{equation}
\end{itemize}

\medskip
\noindent
Proof of (v). For the particular cases $Y=\pm J$, $\left(\begin{array}{cc}
               1 & 0 \\
               \kappa & 1
             \end{array}\right)$, $\left(\begin{array}{cc}
               \lambda & 0 \\
               0 & \lambda^{-1}
             \end{array}\right)$, the identity (\ref{conj_quadratic}) is shown by direct computations with assertions (i) -- (iii).
Then, for general $Y\in\SL$, (\ref{conj_quadratic}) is deduced through the decompositions (\ref{decomp_a}) and (\ref{decomp_b}) and the assertion (iv).\qed

\medskip
\noindent
Proof of (vi). Given $L\in\sl$, there exists $C_L\in\SL$ such that $L=C_L\tilde L C_L^{-1}$ with
\begin{itemize}
  \item $\tilde L= \left(\begin{array}{cc}
               0 & \theta \\
               -\theta & 0
             \end{array}\right)$ for some $\theta\in\R\setminus\{0\}$ if ${\rm det}(L)>0$
  \item $\tilde L=\left(\begin{array}{cc}
               0 & 0 \\
               \kappa & 0
             \end{array}\right)$ for some $\kappa\in\R$ if ${\rm det}(L)=0$
  \item $\tilde L=\left(\begin{array}{cc}
               l & 0 \\
               0 & -l
             \end{array}\right)$ for some $l\in\R\setminus\{0\}$ if ${\rm det}(L)<0$.
\end{itemize}
By computation with the assertions (i) -- (iii), we see that
$$e^{-{\rm i}t\tilde{\CL}}={\CM}(e^{t\tilde L})=\left\{\begin{array}{ll}
                               {\CM}\left(\begin{array}{cc}
               \cos(\theta t) & \sin(\theta t) \\
               -\sin(\theta t) & \cos(\theta t)
             \end{array}\right),  &  {\rm det}(L)>0 \\[4mm]
                                {\CM}\left(\begin{array}{cc}
               1 & 0 \\
               t\kappa & 1
             \end{array}\right) , &  {\rm det}(L)=0 \\[4mm]
                                {\CM}\left(\begin{array}{cc}
               e^{tl} & 0 \\
               0 & e^{-tl}
             \end{array}\right) , & {\rm det}(L)<0
                              \end{array}
\right.  ,$$
with $\tilde{\CL}$ defined as
$$\tilde{\CL}:=-\frac12 \left\la \left(\begin{array}{c}
X\\
D
 \end{array}\right) , J\tilde L \left(\begin{array}{c}
X\\
D
 \end{array}\right) \right\ra . $$
In view of the assertion (v), we have
\begin{eqnarray*}
{\CM}(C_L) \tilde{\CL} {\CM}(C_L)^{-1}&=& \frac{-1}2 \left\la \left(\begin{array}{c}
X\\
D
 \end{array}\right) , JC_L\tilde L C_L^{-1}\left(\begin{array}{c}
X\\
D
 \end{array}\right) \right\ra \\
&=& \frac{-1}2 \left\la \left(\begin{array}{c}
X\\
D
 \end{array}\right) , JL\left(\begin{array}{c}
X\\
D
 \end{array}\right) \right\ra\\
&=& {\CL}.
\end{eqnarray*}
Since ${\CM}(C_L)$ is unitary, we have
$$e^{-{\rm i}t {\CL}}={\CM}(C_L)e^{-{\rm i}t\tilde{\CL}}{\CM}(C_L)^{-1}={\CM}(C_L){\CM}(e^{t\tilde L}){\CM}(C_L)^{-1}={\CM}(C_L e^{t\tilde L} C_L^{-1})={\CM}(e^{tL}).$$
The assertion (vi) is shown. \qed

\medskip
\noindent
Proof of (vii). Since the two linear systems in (\ref{linear_systems}) are conjugated to each other via the change of variable in (\ref{conj_ode}), we see that
$$L_2(t)=-Y(t)^{-1}\dot{Y}(t)+Y(t)^{-1}L_1(t)Y(t).$$
Hence, by the assertion (v), for every $t\in\R$,
\begin{eqnarray*}
\CL_2(t)&=&\frac{-1}2 \left\la \left(\begin{array}{c}
X\\
D
 \end{array}\right) , J L_2(t) \left(\begin{array}{c}
X\\
D
 \end{array}\right) \right\ra\\
 &=&\frac{1}2 \left\la \left(\begin{array}{c}
X\\
D
 \end{array}\right) , J Y(t)^{-1}\dot{Y}(t) \left(\begin{array}{c}
X\\
D
 \end{array}\right) \right\ra +  \frac{-1}2 \left\la \left(\begin{array}{c}
X\\
D
 \end{array}\right) , J Y(t)^{-1}L_1(t) Y(t) \left(\begin{array}{c}
X\\
D
 \end{array}\right) \right\ra\\
 &=&\frac{1}2 \left\la \left(\begin{array}{c}
X\\
D
 \end{array}\right) , J Y(t)^{-1}\dot{Y}(t) \left(\begin{array}{c}
X\\
D
 \end{array}\right) \right\ra + \CM(Y(t))^{-1}\CL_1(t)\CM(Y(t)).
\end{eqnarray*}
To show the conjugacy between two quadratic quantum Hamiltonians in (\ref{quantum_Hamiltonians}) via the unitary transformation (\ref{conj_pde}),  it is sufficient to show that
\begin{equation}\label{eq_time_deriv}
\CM(Y(t))^{-1} \, {\rm i}\partial_t \CM(Y(t))= \,  \frac{-1}2 \left\la \left(\begin{array}{c}
X\\
D
 \end{array}\right) , J Y(t)^{-1}\dot{Y}(t) \left(\begin{array}{c}
X\\
D
 \end{array}\right) \right\ra,
\end{equation}
where
 for $Y(t)=\left(\begin{array}{cc}
               Y_{11}(t) & Y_{12}(t)\\
               Y_{21}(t) & Y_{22}(t)
             \end{array}\right)$ and $v\in L^2(\R)$, ${\rm i}\partial_t \CM(Y(t))$ is given as
\begin{itemize}
  \item if $Y_{11}\neq 0$, then, with (\ref{Meta_a}),
  \begin{eqnarray}
  {\rm i}\partial_t \CM(Y) v (x) &=&{\rm i} \partial_t\left(\frac{1}{\sqrt{2\pi Y_{11}}}\int_{\R} e^{\frac{\rm i}{2Y_{11}}(Y_{21}x^2+2x\xi-Y_{12}\xi^2)} \widehat v(\xi) d\xi \right) \label{idt_form1}\\
   &=&\frac{-{\rm i}{\dot Y}_{11}}{2Y_{11}}\CM(Y) v (x) -  \frac{Y_{11}\dot{Y}_{21}-\dot{Y}_{11}Y_{21}}{2  Y_{11}^{2}}X^2 \CM(Y) v (x) +   \frac{\dot{Y}_{11}}{Y_{11}^{2}}\CI_1(x)\nonumber\\
   & &+ \, \frac{Y_{11}\dot{Y}_{12}-\dot{Y}_{11}Y_{12}}{2Y_{11}^{2}}\CI_2(x)\nonumber
\end{eqnarray}
with the integral terms $\CI_1(x)$ and $\CI_2(x)$ defined by
\begin{eqnarray*}
\CI_1(x)&:=&\frac{1}{\sqrt{2\pi Y_{11}}} \int_{\R} x\xi e^{\frac{\rm i}{2Y_{11}}\left(Y_{21}x^2+2x\xi-Y_{12}\xi^2\right)} \widehat v(\xi) d\xi,\\
\CI_2(x)&:=&\frac{1}{\sqrt{2\pi Y_{11}}} \int_{\R} \xi^2 e^{\frac{\rm i}{2Y_{11}}\left(Y_{21}x^2+2x\xi-Y_{12}\xi^2\right)} \widehat v(\xi) d\xi.
\end{eqnarray*}
\item if $Y_{11}=0$, then, with (\ref{Meta_b}),
\begin{align}
{\rm i}\partial_t \CM(Y(t)) v (x) &= {\rm i}\frac{d}{dt}\left(\frac{{\rm i}^{\frac12}}{\sqrt{2\pi Y_{12}}}\int_{\R} e^{\frac{\rm i}{2Y_{12}}\left(Y_{22}x^2+2xy\right)}  v(y) dy \right)\label{idt_form2}\\
&=  \frac{-{\rm i}{\dot Y}_{12}}{2  Y_{12}} \CM(Y) v (x) -  \frac{Y_{12}\dot{Y}_{22}-\dot{Y}_{12}Y_{22}}{2  Y_{12}^{2}}  X^2 \CM(Y) v (x) +  \frac{\dot{Y}_{12}}{Y_{12}^{2}}\CJ_1(x)\nonumber
\end{align}
with the integral term $\CJ_1(x)$ defined by
\begin{eqnarray*}
\CJ_1(x):=\frac{{\rm i}^{\frac12}}{\sqrt{2\pi Y_{12}}}\int_{\R} xy e^{\frac{\rm i}{2Y_{12}}\left(Y_{22}x^2+2xy\right)}  v(y) dy.
\end{eqnarray*}
\end{itemize}

We are going to show (\ref{eq_time_deriv}) in the following two cases.

\noindent
Case 1. $Y_{11}\neq 0$.
For $v\in L^2(\R)$, by direct computations on $\CM(Y) v$ with $\CM(Y)$ given by (\ref{Meta_a}), we see that
\begin{eqnarray*}
(XD+DX)\CM(Y)v&=& -{\rm i} \CM(Y)v+\frac{2Y_{21}}{Y_{11}}X^2\CM(Y)v+\frac{2}{Y_{11}}\CI_1 , \\
D^2 \CM(Y)v&=& -\frac{{\rm i}Y_{21}}{Y_{11}}\CM(Y)v +\frac{Y^2_{21}}{Y^2_{11}}X^2\CM(Y)v+\frac{2 Y_{21} }{Y^2_{11}}\CI_1+\frac{1}{Y^2_{11}}\CI_2,
\end{eqnarray*}
which allows us to present the two integral term $\CI_1$ and $\CI_2$ as:
\begin{eqnarray*}
\CI_1&=& \left( \frac{Y_{11}}{2}(XD+DX)+\frac{{\rm i} Y_{11}}{2}  {\rm Id}-Y_{21}X^2\right)\CM(Y)v ,\\
\CI_2&=&  \left(- Y_{11}Y_{21}(XD+DX)+Y^2_{11}D^2+Y^2_{21}X^2\right)\CM(Y)v.
\end{eqnarray*}
Hence, according to (\ref{idt_form1}), we have
\begin{eqnarray*}
 & & {\rm i}\partial_t \CM(Y) v \\
  &=& \frac{- {\rm i}{\dot Y}_{11}}{2Y_{11}} \CM(Y)v -\frac{Y_{11}\dot{Y}_{21}-\dot{Y}_{11}Y_{21}}{2Y^2_{11}}   X^2 \CM(Y)v +  \frac{\dot{Y}_{11}}{Y^2_{11}}\CI_1 +\frac{Y_{11}\dot{Y}_{12}-\dot{Y}_{11}Y_{12}}{2Y_{11}^{2}}\CI_2   \\
   &=& \frac{-1}{2Y^2_{11}}\left[\left(Y_{11}\dot{Y}_{21}+\dot{Y}_{11}Y_{21}-Y_{11}\dot{Y}_{12}Y_{21}^2+\dot{Y}_{11}Y_{12}Y_{21}^2\right)X^2 +Y^2_{11}(\dot{Y}_{11}Y_{12}-Y_{11}\dot{Y}_{12})D^2  \right.\\
   & & \ \ \ \ \ \  \ \ \ \left. + \,  \left( -Y_{11}\dot{Y}_{11}+ Y^2_{11}\dot{Y}_{12}Y_{21}-Y_{11}\dot{Y}_{11}Y_{12}Y_{21}\right) (XD+DX) \right] \CM(Y) v .
\end{eqnarray*}
Since $Y\in \SL$ means that
$$Y_{22}=\frac{1+Y_{12}Y_{21}}{Y_{11}} ,\quad \dot{Y}_{22}=\frac{\dot{Y}_{12}Y_{21}+Y_{12}\dot{Y}_{21}}{Y_{11}}-\frac{\dot{Y}_{11}(1+Y_{12}Y_{21})}{Y^2_{11}},$$
the above formula of $ {\rm i}\partial_t \CM(Y) v$ implies that
$${\rm i}\partial_t \CM(Y)=\frac{-1}2 \left\la \left(\begin{array}{c}
X\\
D
 \end{array}\right) , J \dot{Y} Y^{-1} \left(\begin{array}{c}
X\\
D
 \end{array}\right) \right\ra   \CM(Y). $$
Hence, with the assertion (v), we obtain (\ref{eq_time_deriv}) by seeing that
\begin{eqnarray*}
  \CM(Y)^{-1} \, {\rm i}\partial_t \CM(Y)
 &=& \CM(Y)^{-1}  \, \frac{-1}2 \left\la \left(\begin{array}{c}
X\\
D
 \end{array}\right) , J \dot{Y} Y^{-1} \left(\begin{array}{c}
X\\
D
 \end{array}\right) \right\ra   \CM(Y) \\
&=&\frac{-1}2 \left\la \left(\begin{array}{c}
X\\
D
 \end{array}\right) , J Y^{-1} \dot{Y} Y^{-1} Y \left(\begin{array}{c}
X\\
D
 \end{array}\right) \right\ra   \\
&=&  \frac{-1}2 \left\la \left(\begin{array}{c}
X\\
D
 \end{array}\right) , J Y^{-1}\dot{Y} \left(\begin{array}{c}
X\\
D
 \end{array}\right) \right\ra .
\end{eqnarray*}

\noindent
Case 2. $Y_{11}= 0$. For $v\in L^2(\R)$, by direct computations of $(XD+DX)\CM(Y)v$ with $\CM(Y)$ given by (\ref{Meta_b}), we see that
$$\CJ_1= \left( \frac{Y_{12}}{2}(XD+DX)+\frac{{\rm i} Y_{12}}{2}  {\rm Id}-Y_{22}X^2\right)\CM(Y)v. $$
Hence, according to (\ref{idt_form2}), we have
\begin{eqnarray*}
 {\rm i}\partial_t \CM(Y) v
  &=& \frac{- {\rm i}{\dot Y}_{12}}{2Y_{12}} \CM(Y)v -  \frac{Y_{12}\dot{Y}_{22}-\dot{Y}_{12}Y_{22}}{2  Y_{12}^{2}}  X^2 \CM(Y)v +  \frac{\dot{Y}_{12}}{Y^2_{12}}\CJ_1  \\
   &=& \frac{-1}{2Y^2_{12}}\left[\left(Y_{12}\dot{Y}_{22}+\dot{Y}_{12}Y_{22}\right)X^2  - Y_{12}\dot{Y}_{12}(XD+DX) \right] \CM(Y) v .
\end{eqnarray*}
Since $Y=\left(\begin{array}{cc}
               0 & Y_{12}\\
               Y_{21} & Y_{22}
             \end{array}\right)\in \SL$ means that
$$Y_{21}=\frac{-1}{Y_{12}} ,\quad \dot{Y}_{21}=\frac{\dot{Y}_{12}}{Y^2_{12}},$$
the above formula of $ {\rm i}\partial_t \CM(Y) v$ implies that
$${\rm i}\partial_t \CM(Y)=\frac{-1}2 \left\la \left(\begin{array}{c}
X\\
D
 \end{array}\right) , J \dot{Y} Y^{-1} \left(\begin{array}{c}
X\\
D
 \end{array}\right) \right\ra   \CM(Y). $$
Hence, we obtain (\ref{eq_time_deriv}) as in Case 1. \qed

\section{Proof of Proposition \ref{prop_al_ODE}}\label{App_proof_ar}

The almost reducibility of the quasi-periodic linear system $(\omega, J+P(\cdot))$ was shown by Eliasson \cite{Eli1992}.
Indeed, if $\|P\|_r:=\varepsilon_0$ is small enough (depending on $r,\gamma,\tau,d$), then there exists sequences $\{Y_l\}_{l\in\N^*}\subset C^\omega(2\T^d, {\rm SL}(2,\R))$, $\{A_l\}_{l\in\N^*}\subset {\rm sl}(2,\R)$, and $\{F_l\}_{l\in\N^*}\subset C^\omega(\T^d, {\rm sl}(2,\R))$
 with $\|F_l\|_{\T^d}$ bounded by some $\tilde\varepsilon_l$, which is much smaller than $\varepsilon_0$, such that, with $A_0=J$ and $F_0(\cdot)=P(\cdot)$,
 \begin{equation}\label{conj_KAM_Eliasson}
 \frac{d}{dt}Y_{l}(\omega t)=\left(A_l+F_l(\omega t)\right)Y_{l}(\omega t)-Y_{l}(\omega t) \left(A_{l+1}+F_{l+1}(\omega t)\right).
 \end{equation}
 Note that the above equality means that, via the change of variables $x_l=Y_{l}(\omega t)x_{l+1}$, the linear system $\dot{x}_l=\left(A_l+F_l(\omega t)\right)x_l$ is conjugated to $\dot{x}_{l+1}=\left(A_{l+1}+F_{l+1}(\omega t)\right)x_{l+1}$.
More precisely, at the $l-$th step, for $\pm{\rm i}\xi_l\in\R\cup{\rm i}\R$, two eigenvalues of $A_l$, and
$$N_l:=\frac{2|\ln\tilde\varepsilon_l|}{r_l-r_{l+1}},\qquad r_0=r, \quad r_{l+1}=r_l-\frac{r_0}{2^{l+2}},$$
there are two cases about the construction of the conjugation in the KAM step.
\begin{itemize}
  \item (Non-resonant case) If for every $k\in\Z^d$ with $0<|k|\leq N_l$, we have
  \begin{equation}\label{non_resonant}
  \left|2\xi_l-\la k\ra\right|\geq \tilde\varepsilon_l^{\frac{1}{15}},
  \end{equation}
  then $Y_{l}(\cdot)=e^{\tilde Z_{l}(\cdot)}$ for some $\tilde Z_{l}\in C^{\omega}(2\T^d,{\rm sl}(2,\R))$ with
 $$ \|\tilde Z_{l}\|_{2\T^d}<\tilde\varepsilon_l^{\frac12} ,\quad \|A_{l+1}-A_l\|<\tilde\varepsilon_l^{\frac12}, \quad  \|F_{l+1}\|_{\T^d}<\tilde\varepsilon_{l+1}=\tilde\varepsilon_{l}^{2}.$$

\smallskip

  \item (Resonant case) If for some $n_l\in\Z^d$ with $0<|k_l|\leq N_l$, we have
  \begin{equation}\label{resonant}
  \left|2\xi_l-\la k_l\ra\right|< \tilde\varepsilon_l^{\frac{1}{15}},
  \end{equation}
  then there exists $\tilde Z_{l}\in C^{\omega}(2\T^d,{\rm sl}(2,\R))$ with $\|\tilde Z_{l}\|_{2\T^d}<\tilde\varepsilon_l^{\frac12}$ such that
  $$Y_{l}(\cdot)=C_{A_l}\cdot e^{\tilde Z_{l}(\cdot)}\cdot R_{\frac{\la k_l ,\cdot\ra}{2}},$$
  where $C_{A_l}\in \SL$ such that $A_l=C_{A_l}\left(\begin{array}{cc}
                                  0 & \xi_{l} \\
                                  -\xi_{l}  & 0
                                \end{array}\right) C_{A_l}^{-1}$ satisfying
$$\|C_{A_l}\|\leq 4\pi\sqrt{\frac{\|A_l\|}{\xi_l}}. $$
Note that the Diophantine condition of $\om$ and the resonant condition (\ref{resonant}) imply that
$$\xi_l   > \frac12\left(|\la k_{l}\ra|-\left| 2\xi_{l} - \la k_{l}\ra\right|\right)
>\frac12\left(\frac{\gamma}{|k_{l}|^{\tau}}-\tilde\varepsilon_{l}^{\frac{1}{15}}\right)
>\frac{\gamma}{|\ln\tilde\varepsilon_{l}|^{\tau}}, $$
and hence $\|C_{A_l}\|< |\ln\tilde\varepsilon_{l}|^{\frac34\tau}$. Moreover, we have
  $$\|A_{l+1}\|<\frac{\tilde\varepsilon_l^{\frac{1}{16}}}4,\quad  \|F_{l+1}\|_{\T^d}<\tilde\varepsilon_{l+1}=\tilde\varepsilon_l \, {\rm exp}\left\{-r_{l+1}\tilde\varepsilon_l^{-\frac{1}{18\tau}}\right\}.$$
\end{itemize}
As $l$ goes to $\infty$, the time dependent part $F_{l}$ tends to vanish.
Hence the linear system $(\omega,\, J+P)$ is almost reducible.

Let ${\CR}:=\{l_j\}\subset \N^*$ be the collection of indices of KAM steps where the resonant case occurs, and let us focus on the state just in front of the resonant KAM step. More precisely, for $j\geq 0$, let us define
\begin{equation}\label{quanti_res}
L_j:=A_{l_{j+1}},\quad P_{j+1}:=F_{l_{j+1}},\quad U_j(t)=\prod_{l=0}^{l_{j+1}-1} Y_{l}(\om t).
\end{equation}
We see that, by the description of the resonant case, ${\rm det}(L_j)=\xi^2_{l_{j+1}}$ with
$$\varrho_j  := \xi_{l_{j+1}} >|\ln\tilde\varepsilon_{l_{j+1}}|^{-\frac54\tau}. $$
Moreover, since there are only non-resonant steps between the systems $(\om, A_{l_{j}+1}+F_{l_{j}+1})$ and $(\om, A_{l_{j+1}}+F_{l_{j+1}})$, we have
\begin{equation}\label{esti_Lj}
\|L_j\|=\|A_{l_{j+1}}\|\leq  \|A_{l_{j}+1}\|+\sum_{l=l_{j}+1}^{l_{j+1}}  \|A_{l+1}-A_l\| <\frac{\tilde\varepsilon_{l_j}^{\frac1{16}}}{4}+\sum_{l=l_{j}+1}^{l_{j+1}}\tilde\varepsilon_{l}^{\frac12}<\tilde\varepsilon_{l_j}^{\frac1{16}}.
\end{equation}

In view of \cite{Eli1992}, there are two cases about the almost reducibility.
\begin{itemize}
\item If $\sharp\CR<\infty$, then $(\omega, \, J+P(\cdot))$ is reducible.
For $l_{\bar j}:=\max \CR$, we have
\begin{equation}\label{J_large}
Y_{l}(\cdot)=e^{\tilde Z_{l}(\cdot)} \  \ {\rm with} \   \  \|\tilde Z_{l}\|_{2\T^d}<\tilde\varepsilon^{\frac12}_{l},\qquad \forall \ l\geq l_{\bar j}+1.
\end{equation}
Then, for $j\geq \bar j+1$, let us define
$$L_j:=\lim_{l\to\infty} A_l,\quad P_{j+1}:=0, \quad U_j(t):=U_{l_{\bar j}}(t)\prod_{l=l_{\bar j}}^{\infty}Y_{l}(\om t).$$
Hence there exist $U(t)\in \SL$ for every $t\in\R$ and $L\in\sl$ such that
$$U_j(t)\to \prod_{l=0}^\infty Y_l(\omega t) =: U(t), \quad L_j\to L.$$
Moreover, via the change of variables $y=U(t)w$, the system $\dot{y}=(J+P(\om t))y$ is conjugated to
 $\dot{w}=L w.$

\smallskip

\item If $\sharp\CR=\infty$, then we define the sequence $\{\varepsilon_{j}\}$ by $\varepsilon_{j}:=\tilde\varepsilon_{l_{j}}$. By (\ref{quanti_res}) and (\ref{esti_Lj}), we have
$$\|L_j\|< \varepsilon_{j}^{\frac1{16}},\quad \|P_{j+1}\|_{\T^d}<\tilde\varepsilon_{l_{j+1}}=\varepsilon_{j+1}.$$
According to the construction of the sequence $\{\tilde\varepsilon_{l}\}$, we see that
$$ \varepsilon_{j+1}=\tilde\varepsilon_{l_{j+1}}<\tilde\varepsilon_{l_{j}} {\rm exp}\left\{-r_{l_{j}+1}\tilde\varepsilon_{l_{j}}^{-\frac{1}{18\tau}}\right\}<\varepsilon_j \, {\rm exp}\left\{-\frac{r}{2}\varepsilon_{j}^{-\frac{1}{18\tau}}\right\},$$
which implies that
$$|\ln\varepsilon_{j+1}|> |\ln\varepsilon_{j}|+ \frac{r}{2}\varepsilon_{j}^{-\frac{1}{18\tau}}\geq |\ln\varepsilon_{j}|^{200}.$$
Hence, for $j,k\in\N$ with $k\leq j$, we have
$|\ln\varepsilon_{k}|\leq |\ln\varepsilon_{j}|^{\left(\frac{1}{200}\right)^{j-k}}$. For
$$U_j(t)=\prod_{l=0}^{l_{j+1}-1} Y_{l}(\omega t),$$
we have that, for $l\in\CR$,
$$\|Y_{l}\|_{2\T^d}\leq \|C_{A_l}\|\cdot \|e^{\tilde Z_{l}}\|_{\T^d}\cdot \|R_{\frac{\la n_l ,\cdot\ra}{2}}\|_{2\T^d}\leq6|\ln\tilde\varepsilon_{l}|^{\frac34\tau}<|\ln\tilde\varepsilon_{l}|^{\tau}.$$
Therefore, we have the estimates
$$ \prod_{0\leq l\leq l_{j+1}-1\atop{l\in\CR}}  \|Y_{l}\|_{2\T^d} \leq  \prod_{i=0}^{j} \|Y_{l_i}\|_{2\T^d}\leq \prod_{i=0}^{j}|\ln\tilde\varepsilon_{l_i}|^{\tau} \leq |\ln\tilde\varepsilon_{l_{j}}|^{\tau\sum_{i=0}^j\left(\frac{1}{200}\right)^{j-i}}\leq |\ln\tilde\varepsilon_{l_j}|^{\frac{200}{199}\tau},$$
$$ \prod_{0\leq l\leq l_{j+1}-1\atop{l\not\in\CR}}\|Y_{l}\|_{2\T^d}\leq \prod_{0\leq l\leq l_{j+1}-1}\left(1+\tilde\varepsilon_{l}^{\frac23}\right)<3.$$
Hence, for $j\geq 1$, we have
$$\sup_t\|U_j(t)\|\leq \prod_{l=0}^{l_j-1} \|Y_{l}\|_{2\T^d}< 3 |\ln\tilde\varepsilon_{l_j}|^{\frac{200}{199}\tau}< |\ln\tilde\varepsilon_{l_j}|^{2\tau}=|\ln\varepsilon_{j}|^{2\tau}.$$
\end{itemize}
Proposition \ref{prop_al_ODE} is shown.

\end{document}